\documentclass[a4paper,10pt]{article}
\usepackage{amsmath, amsthm, amssymb}
\usepackage{url}

\usepackage{physics}
\usepackage{graphicx,lipsum}
\graphicspath{ {./downloads/} }

\usepackage[colorlinks,citecolor=red,urlcolor=blue,bookmarks=true,hypertexnames=true]{hyperref}
\usepackage{tikz}
\usetikzlibrary{calc}
\usetikzlibrary{shapes}
\usepackage[autostyle]{csquotes}
\makeatletter

\usepackage[colorlinks]{hyperref}
\usepackage[nameinlink,capitalize]{cleveref}
\newtheorem{theorem}{Theorem}[section]
\newtheorem{corollary}{Corollary}[theorem]

\newtheorem{claim}{Claim}[section]
\newtheorem{lemma}[theorem]{Lemma}
\newtheorem{prop}{Proposition}[section]
\newtheoremstyle{named}{}{}{\itshape}{}{\bfseries}{.}{.5em}{\thmnote{#3's }#1}
\theoremstyle{named}

\newtheorem{remark}{Remark}
\theoremstyle{definition}
\newtheorem{definition}{Definition}[section]

\theoremstyle{remark}

\theoremstyle{conclusion}

\theoremstyle{observation}

\usepackage{braket}
\usepackage{xspace}
\usepackage[margin=1.0in]{geometry}

\usepackage{titlesec}

\usepackage{mathtools}

\DeclarePairedDelimiterX{\inp}[2]{\langle}{\rangle}{#1, #2}
\titleformat{\chapter}
  {\Large\bfseries}
  {}
  {0pt}
  {\huge}
\usepackage{makeidx}  
\usepackage{dsfont}
\usepackage{amsmath,amssymb,amsfonts,amscd}
\usepackage{hyperref}
\usepackage{booktabs} 
\usepackage{appendix}
\usepackage[utf8]{inputenc}
\usetikzlibrary{3d,positioning,calc}
\usepackage{tikz,tikz-cd}
\tikzset{
    invisible/.style={opacity=0},
    visible on/.style={alt={#1{}{invisible}}},
    alt/.code args={<#1>#2#3}{%
      \alt<#1>{\pgfkeysalso{#2}}{\pgfkeysalso{#3}}%
  }
}
\tikzset{
  symbol/.style={
    draw=none,
    every to/.append style={
      edge node={node [sloped, allow upside down, auto=false]{$#1$}}}
  }
}
\usetikzlibrary{cd}
\usepackage{xcolor}
\usepackage{graphicx}
\usepackage{cancel}

\usepackage{fancyhdr} 

\fancypagestyle{firstpage}{%
  \fancyhf{} 
  \fancyfoot[C]{\footnotesize
  }
}

\begin{document}
\thispagestyle{firstpage}

\begin{center}
  \fontsize{13pt}{16pt}\selectfont
  \textsc{\textbf{CONSTRUCTING AND TRANSFORMING STANDARD RESOLUTIONS OF LOCAL MODULES}}
\end{center}

\vspace{0.2cm}

\begin{center}
  \fontsize{10pt}{14pt}\selectfont
  \textsc{HAO XINGBANG}
\end{center}

\vspace{0.4cm}

\begin{center}
  \fontsize{12pt}{14pt}\selectfont
  \textsc{Abstract}
\end{center}

\noindent

We construct standard resolutions for analytic local modules on complex hypersurfaces via standard basis methods, with applications to derived category formulas for local hypersurface blow-ups, revealing new links between standard bases and categorical geometry.

\tableofcontents

\section{Introduction}

\subsection{Main Results}

This article establishes two primary conclusions. First, we expect to construct an acyclic complex over a local ring $(R,\mathfrak m)$ whose terms are direct sums of powers of a prime ideal $P$:
\[
\begin{tikzcd}
\cdots \arrow[r] & \bigoplus^{m_3}_{l=1} P^{\mathrm{a}^{(3)}_l} \arrow[r, "F_2"] &
\bigoplus^{m_2}_{l=1} P^{\mathrm{a}^{(2)}_l} \arrow[r, "F_1"] &
\bigoplus^{m_1}_{l=1} P^{\mathrm{a}^{(1)}_l} \arrow[r, "F_0"] &
\bigoplus^{m_0}_{l=1} P^{\mathrm{a}^{(0)}_l} \arrow[r] & 0 .
\end{tikzcd}
\]
The key property is that the complex remains acyclic after multiplication by $P^r$ for any $r \ge 0$. This condition has a natural geometric interpretation: blowing up the subscheme $Z(P)$ in $X = \operatorname{Spec} R$ yields a variety $X^+$. The complex above then provides a projective resolution of certain coherent sheaf $\mathcal{F}$ on $X^+$.
Via Serre's correspondence,
\[
\bigoplus_{i\ge 0} H^0(X^+,\mathcal F\otimes\mathcal O_{X^+}(i)),
\]
this construction descends to a complex of $R$-modules satisfying the same acyclicity resolution. Although such exact sequences are implicit in geometric settings, explicit algebraic constructions appear to be largely absent from the literature.

Motivated by this geometric picture, we expect that for any regular (analytic) local ring $(R,\mathfrak m)$ and certain $R$-modules $M$, there exists a projective resolution inducing exact sequences with the above property, without modifying the underlying morphisms. We refer to such a resolution as the \emph{standard resolution} of $M$. A natural approach is to study the cokernels of $F_0$ after multiplication by $P^r$, which produce a filtration of $M$. However, the analysis of these filtrations is technically subtle. Moreover, the non-uniqueness of syzygies leads to many resolutions, motivating the study of their parameter spaces.

Inspired by Gr\"obner basis and standard basis techniques, we observe that the construction of projective resolutions depends on several choices, including valuations on $R$, module orderings, and syzygy degrees of freedom. Even when such resolutions are homotopy equivalent, their explicit forms may differ substantially, making the identification of a canonical standard resolution a fundamental problem.

For torsion-free modules $M$, we construct an arbitrary Bourbaki sequence to reduce the problem to finding a standard resolution of an ideal $I$. This reduction avoids complications arising from module orderings in syzygy choices. Applying standard basis techniques with the flag order $(\mathrm{ord}_P,\mathrm{ord}_{\mathfrak m})$, we obtain in Section~\ref{subsec:reg} a projective resolution via degree restrictions. This construction automatically satisfies the desired properties and yields a canonical resolution for this class of modules.

This approach relies crucially on $R$ being a regular local (analytic/formal) ring. For more general local rings, classical syzygy and standard basis algorithms no longer apply. But as least in the hypersurface case $S:=R/(w)$, Eisenbud's work \cite{Eisenbud1980} provides an alternative via matrix factorizations. For a closed subscheme $B \subset Y:=\operatorname{Spec} S$, we can construct a standard resolution of its ideal $I_B$ on $R$.
 This raises the question of how syzygies computed in $R$ relate to those computed in $S$. A partial answer is given by the distinguished triangle
\[
\mathcal O_B\otimes\mathcal O_Y(-Y)[1]\longrightarrow j^*j_*\mathcal O_B\longrightarrow \mathcal O_B ,
\]
which suggests that the complexes $j^*j_*\mathcal O_B$ and its shifts encode the syzygies of $B$ on $Y$. In Section~\ref{subsec:hyp}, we explicitly construct the relevant morphisms using division algorithms and show that the resulting standard resolution in $S$ exhibits the expected periodicity of matrix factorizations.

The author subsequently learned that the above construction is precisely the one first introduced by Shamash in \cite{Shamash1969}. For a general (regular) ring \(R\) and an \(R/(w)\)-module, Shamash's construction provides a method for obtaining an \(R/(w)\)-free resolution from an \(R\)-free resolution. It is now commonly referred to as the Eisenbud--Shamash resolution. We refer the reader to Avramov \cite{Avramov1989} for the first modern treatment of the original argument and to \cite{Avramov2002} for the case of complete intersections. Moreover, the Eisenbud--Shamash algorithm is now implemented in Macaulay2 \cite{Eisenbud2001} and is widely used in computational commutative algebra, and is closely related to dg algebras. We emphasize the connection between this resolution and valuations, and expect that the construction extends naturally to arbitrary complete intersections as \cite{Avramov2002} in the local setting. On the other hand, our results presently do not appear to admit a straightforward extension to general non-local rings.

Second, although our constructions apply only to certain subschemes and torsion-free modules, they already admit a natural interpretation in the framework of derived categories. Consider the following geometric setting:
\begin{itemize}
  \item $R$ is a regular \textbf{formal/analytic/algebraic} local ring, and \(X = \operatorname{Spec} R\);
  \item Let $Z:=Z(P) \subset X$ be a smooth center, and let $X^+ = \mathrm{Bl}_Z X$ denote the blow-up of $X$ along $Z$.
  \item Let $Y \subset X$ be a normal hypersurface containing $Z$ and $Z \subset \mathrm{Sing}(Y)$.
  \item Let $Y^+$ be the strict transform of $Y$ in $X^+$. Its exceptional divisor $E_Y$ is Fano (and flat) over $Z$.
\end{itemize}
We study the relationship between $\mathrm D^b(Y^+)$ and $\mathrm D^b(Y)$. Bondal-Orlov's conjecture, verified in many cases in \cite[Theorem~8.22]{Efimov2020}, suggests an equivalence
\[
\mathrm D^b(Y^+)\simeq \mathrm D^b(Y)/\mathcal K_Z ,
\]
where $\mathcal K_Z$ is a kernel category associated with $Z$. This leads to the question of whether $\mathrm D^b(Y)$ embeds fully faithfully into $\mathrm D^b(Y^+)$, or more generally into a Verdier quotient $\mathrm D^b(Y^+)/\mathcal K$. One may further ask whether there exists a universal such subcategory $\mathcal K_u$ through which all embedding factor. Although flipping may occur, as a first step we construct a computable fully faithful functor
\[
\mathcal R\pi_Y^* \colon \mathrm D^b(Y)\longrightarrow
\mathrm D^b(Y^+)/\{\mathcal K:=\langle i^+_*\mathcal A_{-1}\rangle\},
\]
where $\mathcal A_{-1}$ denotes the residue component in the semi-orthogonal decomposition of the (flat) hypersurface fibration \eqref{equ:SODhyper}. The construction is given in Section~\ref{section:stricttrans}. For $F\in\mathrm D^b(Y)$, we decompose a projective resolution into a torsion-free module part and a perfect complex part. The former is treated using standard resolutions and $P^r$-multiplication, while the latter is handled by derived pullback. These pieces are then glued via adjunction to obtain $\mathcal R\pi_Y^*F$.

\begin{prop}[Prop.~\ref{prop:output}]
The functor $\mathcal R\pi_Y^*$ is well-defined, right admissible, and fully faithful.
\end{prop}

The remaining arguments are formal. Using the semi-orthogonal decomposition of $\mathrm D^b(E_Y)$ in \eqref{equ:SODhyper}, we obtain the following result by the same argument as in the regular case.

\begin{theorem}[Thm.~\ref{thm:output} and Thm.~\ref{thm:algebraic-version}]
If $c-d>0$, then
\[
\mathrm D^b(Y^+)/\mathcal K
=
\left\langle
i^+_*\mathcal O_{E_Y}(-(c-d)+1),\ldots,
i^+_*\mathcal O_{E_Y}(-1),
\operatorname{Im}\mathcal R\pi_Y^*
\right\rangle .
\]
\end{theorem}
The author expects that, for a broad class of cases in our theorems, the methods of \cite[Theorem 8.22]{Efimov2020} yield similar conclusions (at least over algebraically closed fields of characteristic zero), though the approach here remains valuable for its distinct technical perspective. The case $c-d\le 0$ appears to lie beyond the reach of purely derived-category methods and suggests deeper connections with monodromy phenomena, in the spirit of Orlov's work \cite{Orlov2009}. This result can be regarded as a slight generalization of Orlov's blow-up formula for smooth varieties \cite{Orlov1993}, but we postpone a more detailed investigation of this connection to our future work.

\subsection{Notation}
In this paper, we assume that $R$ is an essentially/analytic finitely generated local integral domain over a field $\mathrm{k}$ (we do not restrict the choice of field, though the complex numbers are of particular interest), possibly normal. We use $\widehat{R}$ to denote either the formal or analytic completion of $R$ at its maximal ideal $\mathfrak{m}$, as they exhibit essentially the same homological behavior for our purposes. We write $\mathrm{D}^*(X)$ for the derived category of coherent modules (or coherent sheaves) on a ring (or space) $X$, where $* = b, +, -$ refers, respectively, to the bounded, bounded below, and bounded above derived categories. Throughout, all standard functors such as $\pi_*, \pi^*$, etc., are understood in their derived sense unless otherwise specified.
\section{Construction}

Let $(R,\mathfrak m)$ be a regular local ring with residue field $\mathrm k$. Let $J\subset R$ be an ideal such that $S:=R/J$ is normal and Gorenstein, and let $P\subset R$ be a prime ideal containing $J$ such that $Z:=R/P$ is regular.

\begin{definition}[\cite{Yoshino1990}, Prop.~1.2]
An $S$-module $M$ is called a (maximal) Cohen–Macaulay module if
\[
\mathrm{depth}(M)=\mathrm{dim}(S).
\]
This condition is equivalent to each of the following:
\begin{enumerate}
  \item $\mathrm{Ext}_S^{<d}(S/\mathfrak m,M)=0$;
  \item $\mathrm H^{\neq d}_{\mathfrak m}(M)=0$;
  \item $\mathrm{Ext}_S^{>0}(M,S)=0$,
\end{enumerate}
where $d:=\mathrm{dim}(S)$.
\end{definition}

\begin{lemma}[\cite{Yoshino1990}, Cor.~1.13]
If $M$ is a Cohen–Macaulay $S$-module, then its dual
\[
M^\vee:=\mathrm{Ext}^0_S(M,S)
\]
is also Cohen–Macaulay.
\end{lemma}
Cohen–Macaulay modules play a central role in the study of resolutions over Gorenstein local rings, where they naturally appear as the eventual syzygies.
\begin{lemma}[\cite{Buchweitz2021}, Thm.~5.1.2]\label{Buchlemma}
Let $S$ be a strongly Gorenstein ring\footnote{For local rings, the notions of Gorenstein and strongly Gorenstein coincide.}.
\begin{enumerate}
  \item Any finitely generated $S$-module $F$ admits a short exact sequence
  \[
  0 \longrightarrow P_F \longrightarrow M_F \longrightarrow F \longrightarrow 0,
  \]
  where $P_F$ has finite projective dimension and $M_F$ is Cohen–Macaulay.
  \item Such a presentation is unique up to projective stable equivalence: if
  \[
  0 \longrightarrow P'_F \longrightarrow M'_F \longrightarrow F \longrightarrow 0
  \]
  is another such presentation, then there exist finitely generated projective $S$-modules $Q,Q'$ with
  \[
  M_F\oplus Q \cong M'_F\oplus Q',\qquad
  P_F\oplus Q \cong P'_F\oplus Q'.
  \]
  \item The surjection $p\colon M_F\to F$ is universal among morphisms from Cohen–Macaulay $S$-modules to $F$; that is, any morphism $M'\to F$ with $M'$ Cohen–Macaulay factors uniquely through $p$ in $\mathbf{\underline{mod}}\text{-}S$.
\end{enumerate}
\end{lemma}

The above approximation property admits a derived analogue.

\begin{lemma}[\cite{Orlov2004}, Prop.~1.23]
For any object $F\in\mathrm D^b(S)$, there exists a Cohen–Macaulay module $M_F$ and a common anti-roof in $\mathrm D^b(S)$
\[
\begin{tikzcd}
 & F \arrow[dr,"a"] && M_F \arrow[dl,"b"] \\
 && G
\end{tikzcd}
\]
such that both $\mathrm{Cone}(a)$ and $\mathrm{Cone}(b)$ lie in $\mathrm D^{\mathrm{perf}}(S)$. In particular, the images of $F$ and $M_F$ in the singularity category $\mathrm D_{sg}(S)$ are isomorphic.
\end{lemma}

\begin{proof}
Let $(P^\bullet,d^\bullet)$ be a bounded above minimal free resolution of $F$. For $k\gg0$, the module $M:=\ker d^{-k}$ is Cohen–Macaulay by depth lemma. Setting $G:=M[k+1]$ yields the desired diagram.

Next, consider a right free resolution $(Q^\bullet,e^\bullet)$ of $M$, obtained from a minimal free resolution of $M^\vee$. Then $M_F:=\ker e^k$ is Cohen–Macaulay, providing the required approximation.
\end{proof}

Combining the above results, we obtain the following derived Cohen–Macaulay approximation.

\begin{lemma}\label{CMapprox}
Let $S$ be a strongly Gorenstein ring.
\begin{enumerate}
  \item Any $F\in\mathrm D^b(S)$ fits into a distinguished triangle
  \[
  P_F \longrightarrow F \longrightarrow M_F[n_F] \in \Delta,
  \]
  where $P_F\in\mathrm D^{\mathrm{perf}}(S)$ and $M_F$ is a Cohen–Macaulay module. The objects $F$ and $M_F[n_F]$ have isomorphic images in $\mathrm D_{sg}(S)$.
  \item Such a presentation is unique up to projective stable equivalence in the sense of Lemma~\ref{Buchlemma}.
\end{enumerate}
\end{lemma}

\begin{proof}
Take a bounded above minimal free resolution $P^\bullet$ of $F$. For $k\gg0$, the module $\ker d^{-k}$ is Cohen–Macaulay. Choosing $k$ minimal with this property and setting $M_F:=\ker d^{-k}$ and $n_F:=k+1$ yields the desired triangle.
\end{proof}

\begin{lemma}[\cite{Yoshino1990}, Prop.~1.5]
If $S$ is reduced, then every Cohen–Macaulay $S$-module is reflexive, and hence torsion-free.
\end{lemma}

\begin{lemma}[\cite{BourbakiCommutativeAlgebra}, Chap.~VII, \S4.9, Thm.~6]
Let $S$ be a normal ring and $M$ a finitely generated torsion-free $S$-module.
\begin{enumerate}
  \item There exists a free submodule $L\subset M$ such that $M/L$ is isomorphic to an ideal $I\subset S$, yielding a short exact sequence
  \[
  0\longrightarrow L\longrightarrow M\longrightarrow I\longrightarrow 0.
  \]
  \item Any two such presentations correspond to distinct projections of
  \[
  M\otimes_S K(S)=V_{n-1}\oplus K(S),
  \]
  and differ by a unique element $u\in\mathds P^{n}_{K(S)}$, where $n=\mathrm{rank}(M)$.
  \item In general, the associated ideals need not be isomorphic.
\end{enumerate}
\end{lemma}

\begin{definition}
For a Cohen–Macaulay $S$-module $M$, an ideal $I$ obtained as in the preceding lemma is called a \emph{Bourbaki ideal} of $M$, denoted $I(M)$. The corresponding subscheme $B(M)$ is called the \emph{Bourbaki cycle}. Both depend on the choice of Bourbaki sequence, equivalently on the choice of divisors or of the element $u$.
\end{definition}

\subsection{Standard resolution on regular ring}\label{subsec:reg}

For an ideal $I$ in an analytic regular local ring $\widehat{R}$, a projective resolution of $I$ can be constructed by iteratively taking syzygies until a finite free resolution is obtained. We first consider a regular sequence \( x_{1}, \ldots, x_{c}, x_{c+1}, \ldots, x_{n} \) in \( \widehat{R} \) at the maximal ideal \( \mathfrak{\widehat{m}} \), where \( x_{1}, \ldots, x_{c} \) define the center \( \widehat{P} \). Assign the following valuations:
\[
\operatorname{ord}(x_{1}) = \cdots = \operatorname{ord}(x_{c}) = 1, \quad \operatorname{ord}(x_{c+1}) = \cdots = \operatorname{ord}(x_{n}) = 0.
\]
It gives a well-defined valuation to general formal power series. More generally, we have the following equivalent definition:
\begin{definition}
\begin{enumerate}
    \item Let \( R \) be a Noetherian ring, and let \( P \) be a nontrivial ideal of \( R \). For each \( f \in R \), the order of \( f \) along \( P \) is defined as
    \[
    \mathrm{ord}_{P}(f) := \max \{ n\in \mathds{N} \mid f \in P^{n} \}.
    \]
    If \( f \in P^{n} \) for all \( n \), we set \( \mathrm{ord}_{P}(f) = \infty \). In a local ring, this condition holds if and only if \( f = 0 \).

    \item For each \( f \) in the function field \( K(R) \) of \( R \), the order of \( f \) along \( P \) is defined as
    \[
    \mathrm{ord}_{P}(f) := \mathrm{ord}_{P}(g) - \mathrm{ord}_{P}(h),
    \]
    for any representation \( f = g/h \), where \( g \) and \( h \) are regular functions on \( R \).
\end{enumerate}
\end{definition}

For formal power series, leading terms can be defined in two ways. One depends on the element's order, producing homogeneous polynomials \cite{Hironaka1977}. The other depends on monomial ordering (lexicographical and module orders), producing monomials, and is used in Grauert's division theorem and Buchberger's algorithm \cite{Bayer,DeJongPfister}.

\begin{definition}
Let \( 0 \neq f \in \widehat{R}\simeq\mathrm{k}[[\mathbf{x}]]\) (or \(\mathrm{k}\{\mathbf{x}\}\)) be a nonzero multivariate formal (analytic) power series. The leading term of \( f \) is defined as
\[
L(f) = f \mod \widehat{P}^{\mathrm{ord}_{\widehat{P}}(f) + 1},
\]
and its residue term is defined as \( R(f) := f - L(f) \).
\end{definition}

\begin{definition}
Let \( > \) be a well monomial ordering, and let \( 0 \neq f \in \widehat{R}\simeq\mathrm{k}[[\mathbf{x}]] \) (or \(\mathrm{k}\{\mathbf{x}\}\)) be a multivariate formal (analytic)  power series. Then we may write
\[
f = \sum_{i \geq 1} a_i \mathbf{x}^{\alpha_i},
\]
such that \( a_i \neq 0 \), \( a_i \in \mathrm{k} \), and \( \mathbf{x}^{\alpha_i} < \mathbf{x}^{\alpha_{i+1}} \) for all \( i \). We define:
\begin{enumerate}
\item \( \mathrm{in}_{>}(f) = a_1\mathbf{x}^{\alpha_1} \) as the initial monomial of \( f \);
\item \( \mathrm{tail}_{>}(f) = \sum_{i \geq 2} a_i \mathbf{x}^{\alpha_i} \) as the tail series under the monomial ordering \( > \).
\end{enumerate}
\end{definition}

\begin{definition}
For any ideal \( I \) in the completion \( \widehat{R} \), the initial ideal of \( I \) with respect to the monomial ordering \( > \) is defined as:
\[
\mathrm{in}_{>}(I) := \left\langle \mathrm{in}_{>}(f) \mid f \in I, f \neq 0 \right\rangle.
\]
If \( f_1, \ldots, f_m \) is a system of elements of the ideal \( I \), then it is called a standard basis of \( I \) if
\[
\mathrm{in}_{>}(I) = \left\langle \mathrm{in}_{>}(f_1), \ldots, \mathrm{in}_{>}(f_m) \right\rangle.
\]
A standard basis is called reduced if:
\begin{enumerate}
    \item the coefficient \( \mathrm{c}_{>}(f_i) \)  of \( \mathrm{in}_{>}(f_i) \) in \( f_i \) is \( 1 \) for all \( i \), and
    \item \( \mathrm{in}_{>}(f_i) \) does not divide any monomial occurring in \( \mathrm{in}_{>}(f_j) \) for \( i \neq j \) and \( \mathrm{tail}_{>}(f_i) \).
\end{enumerate}
\end{definition}
The following definition was originally introduced by Hironaka as the \emph{standard basis}. In this context, we refer to it as the \emph{measure basis}.
\begin{definition}
 Under the assumption as above, if \( f_1, \ldots, f_m \) is a system of elements of the ideal \( I \), then it is called a measure basis of \( I \) if
\[
L(I) = \left\langle L(f_1), \ldots, L(f_m) \right\rangle.
\]
A measure basis is called reduced if no monomial of the \(R(f_i)\) belongs to \(L(I)\).
\end{definition}
The existence of a standard basis can be ensured by the Noetherian property of the ring. Indeed, we have the following:
\begin{prop}[\cite{DeJongPfister}, Cor. 7.2.11]
For any fixed lexicographical order on \( x_1, \ldots, x_n \), and for any ideal \( I \) in \(\widehat{R} \), there exists a unique reduced standard basis \(f_1, \ldots, f_m\) of \( I \), such that \( \mathrm{in}_{>}(f_1) < \cdots < \mathrm{in}_{>}(f_m) \) and \( I =\langle f_1, \ldots, f_m\rangle\).
\end{prop}
Here, we use the following hierarchical monomial ordering to define $\mathrm{in}(f_i)$:
\begin{itemize}\label{orderchoice}
    \item First, compare elements by the minimal order of $P$.
    \item Next, apply lexicographical ordering within $P$.
    \item Then, compare by the order with respect to $\mathfrak{m}$, and use lexicographical ordering within $\mathfrak{m}$\footnote{This choice does not affect our conclusions; any well-order would suffice. A Gröbner order is preferable when \(R/P\) is a localization of a polynomial ring, or more generally an algebraic regular local ring, since the resulting transformations remain in an algebraic setting, while such a setting sometimes increases the codimension by introducing complete intersections.}
\end{itemize}
The standard basis provides a method for constructing the minimal projective resolution of \( I \) in \( \widehat{R} \). This method is well-known as the Schreyer's syzygies (or Buchberger algorithm), we refer to  \cite[Lemma 7.2.15]{DeJongPfister}  (or \cite[Sec. 15.5]{EisenbudCA}), it proceeds as follows:

\begin{enumerate}
    \item Start with the standard basis of \( I =\langle f_1, \ldots, f_m\rangle\) and initialize \( G = I \).
    \item Assign an order mark to the standard basis of \( I\), defined as:
\[
{\mathrm{ord}_{\widehat{P}}}(I) = \langle {\mathrm{ord}_{\widehat{P}}}(f_1), \ldots, {\mathrm{ord}_{\widehat{P}}}(f_m) \rangle=\langle \mathrm{a}_{1}, \ldots, \mathrm{a}_{m}\rangle.
\]
    \item For each pair \( (f_i, f_j) \) in \( I \), compute the S(lop)-polynomial:
    \[
    S(f_i, f_j) = \frac{\mathrm{LCM}(\mathrm{in}(f_i), \mathrm{in}(f_j))}{\mathrm{in}(f_i)} \cdot f_i - \frac{\mathrm{LCM}(\mathrm{in}(f_i), \mathrm{in}(f_j))}{\mathrm{in}(f_j)} \cdot f_j=p^{i,j}_{i} f_i + p^{i,j}_{j} f_j,
    \]
    where \( \mathrm{LCM} \) denotes the least common multiple (If \( f_i \) and \( f_j \) come from different module bases, the least common multiple (LCM) is defined to be zero).
    \item \label{Schstep4} Reduce \( S(f_i, f_j) \) modulo \( I \) with respect to the standard basis \( \langle f_1, \ldots, f_m \rangle \) by Grauert's division theorem \cite[Thm. 7.1.9]{DeJongPfister}, The division is expressed as:
\[
S(f_i, f_j) = q^{i,j}_{1} f_1 + q^{i,j}_{2} f_2 + \cdots + q^{i,j}_{m} f_m + r,
\]
where \( q^{i,j}_{k} \) are the quotient polynomials and \( r \) is the remainder. This division is uniquely determined with respect to the given order. Moreover, \( \mathrm{in}(p^{i,j}_{i}f_i) < \mathrm{in}(q^{i,j}_{k} f_k) \) for all \( k \), \( i \), and \( j \), and the remainder \( r \) is zero.
 \item Construct a submodule of \( \widehat{R}^m = \bigoplus_{k=1}^{m} \widehat{R} e_k \) generated by all basis elements defined as follows:
\[
f_{i,j} := p^{i,j}_{i} e_i + p^{i,j}_{j} e_j - \left( q^{i,j}_{1} e_1 + q^{i,j}_{2} e_2 + \cdots + q^{i,j}_{m} e_m \right),
\]
for any \( i \neq j \). This module is precisely the syzygy module of \( I \), denoted as \( \mathrm{Syz}(I) \). It fits into the following exact sequence:
\[
0 \longrightarrow \mathrm{Syz}(I) \longrightarrow \widehat{R}^m \longrightarrow I \longrightarrow 0.
\]
\item Following \cite[Thm.\,15.10]{EisenbudCA}, we equip $\widehat{R}^m$ with the natural module ordering induced from $\widehat{R}$, where $m e_i < n e_j$ holds either when $\mathrm{in}(m f_i) < \mathrm{in}(n f_j)$ \label{compactiblemoduleorder}in $\widehat{R}$, or when these initial terms are equal but $i < j$ in the lexicographical order. The order is completed by defining
\[
\mathrm{ord}_P(e_i) := \mathrm{ord}_P(f_i) \quad \text{for all } i=1,\ldots,m.
\]
\item
Under this natural module order, by the Chain Criterion \cite[Lemma 7.2.18]{DeJongPfister}, the set \( \{f_{i,j}\} \) for all \( i \neq j \) contains a unique subset, denoted as \( \{f^{(1)}_{k}\}_{k=1}^{m_{1}} \), such that \( \{f^{(1)}_{1}, \ldots, f^{(1)}_{m_{1}}\} \) forms a (reduced) standard basis for \( \mathrm{Syz}(I) \) \cite[Thm. 15.10]{EisenbudCA}.
\item Therefore, we obtain a new standard basis of \( \mathrm{Syz}(I) = \langle f^{(1)}_1, \ldots, f^{(1)}_{m_{1}} \rangle \) and  there is also a compatible order valuation:
\[
{\mathrm{ord}_{\widehat{P}}}(\mathrm{Syz}(I)) = \langle {\mathrm{ord}_{\widehat{P}}}(f^{(1)}_1), \ldots, {\mathrm{ord}_{\widehat{P}}}(f^{(1)}_{m_{1}}) \rangle = \langle \mathrm{a}^{(1)}_1, \mathrm{a}^{(1)}_2, \ldots, \mathrm{a}^{(1)}_{m_1} \rangle.
\]
then initialize \( G = \mathrm{Syz}(I) \).
    \item Repeat steps (3)–(8) until the syzygy reduce to zero.
\end{enumerate}
By analyzing all possible initial terms in \(\mathrm{Syz}^k(I)\) and the successive elimination of variables along syzygy paths by module order choice in (\ref{compactiblemoduleorder}), one obtains the termination of the syzygy algorithm and we have a finite projective resolution of $I$ by \cite[Cor. 15.11]{EisenbudCA}. for example:
\[
\begin{tikzcd}
0 \arrow[r] & \widehat{R}^{m_k} \arrow[r, "F_k"] & \widehat{R}^{m_{k-1}} \arrow[r, "F_{k-1}"] & \cdots \arrow[r, "F_2"] & \widehat{R}^{m_1} \arrow[r, "F_1"] & \widehat{R}^{m_0} \arrow[r, "F_0"] & \widehat{R} \arrow[r] & \widehat{R}/I \arrow[r] & 0
\end{tikzcd}
\]
we notice that even a choice of reduced standard basis, this process does not necessarily produce a reduced projective resolution, for example the ideal $I = (x^2+y^2, xy, y^3)$ with the ordering $\mathrm{ord}(x) = \mathrm{ord}(y) = 1$ and $x < y$.\\

For \( 0 \leq i \leq k \), the matrix \( F_i \) is defined as:
\[
F_i = \begin{pmatrix}
f^{(i)}_1 & f^{(i)}_2 & \cdots & f^{(i)}_{m_i}
\end{pmatrix},
\]
where \( f^{(i)}_j \) denotes the \( j \)-th column of \( F_i \), and the orders assigned to each column of \( F_i \) are given by:
\[
{\mathrm{ord}_{\widehat{P}}}(F_i) = \langle \mathrm{a}^{(i)}_1, \mathrm{a}^{(i)}_2, \ldots, \mathrm{a}^{(i)}_{m_i} \rangle,
\]
where \( \mathrm{a}^{(i)}_j \) represents the order assigned to the \( j \)-th column \( f^{(i)}_j \) of \( F_i \).
Let \( F_i = \{ f^{(i)}_{m,n} \} \) be a matrix representing \( F_i \), where \( f^{(i)}_{m,n} \) denotes the element in the \( m \)-th row and \( n \)-th column. The \textit{leading matrix} \( L(F_i) \) is defined as:
\[
L(F_i) = \left\{ L(f^{(i)}_{m,n}) \right\}_{(1,1) \leq (m,n) \leq (m_{i-1},m_{i})}
\]
where each entry \( L(f^{(i)}_{m,n}) \) is given by:
\[
L(f^{(i)}_{m,n}) = f^{(i)}_{m,n} \mod \widehat{P}^{\mathrm{a}^{(i)}_{n} - \mathrm{a}^{(i-1)}_{m} + 1}.
\]
similarly, we have its column representation:
\[
L(F_i) = \begin{pmatrix}
L(f^{(i)}_1) & L(f^{(i)}_2) & \cdots & L(f^{(i)}_{m_i})
\end{pmatrix},
\]
where \( L(f^{(i)}_j) \) denotes the \( j \)-th column of \( L(F_i) \).

\begin{lemma}
Under the natural module order on $\widehat{R}^{m_{i-1}}$ induced by the construction, the leading terms
\[
L(f^{(i)}_1), \dots, L(f^{(i)}_{m_i})
\]
form a standard basis of the module they generate, denoted $L(\mathrm{Syz}^{i}(I))$.
\end{lemma}

\begin{proof}
By construction, $\operatorname{in}(L(f^{(i)}_j)) = \operatorname{in}(f^{(i)}_j)$ for all $j$. Hence $L(f^{(i)}_1),\dots,L(f^{(i)}_{m_i})$ generate the initial module and form a standard basis.
\end{proof}

\begin{lemma}\label{lemmaleadingregular}
The sequence
\[
\begin{tikzcd}
0 \arrow[r] & \widehat{R}^{m_k} \arrow[r, "L(F_k)"] & \widehat{R}^{m_{k-1}} \arrow[r, "L(F_{k-1})"] & \cdots \arrow[r, "L(F_1)"] & \widehat{R}^{m_0} \arrow[r, "L(F_0)"] & \widehat{R} \arrow[r] & \widehat{R}/L(I) \arrow[r] & 0
\end{tikzcd}
\]
is exact and coincides with the sequence obtained by applying Schreyer's syzygies to $L(I)$.
\end{lemma}

\begin{proof}
Let $L(I) = \langle L(f_1),\dots,L(f_m) \rangle$ be the standard basis.
The Schreyer construction for syzygies depends only on initial terms:
\[
S(L(f_i), L(f_j)) = p_i^{i,j} L(f_i) + p_j^{i,j} L(f_j),
\]
and the leading term of $S(f_i,f_j)$ with respect to $\operatorname{ord}_{\widehat{P}}$ is identified with $S(L(f_i),L(f_j))$.
By induction, all syzygy matrices $F_i$ have leading matrices $L(F_i)$ forming an exact sequence.
All maps are homogeneous; the sequence is graded and agrees with that of $L(I)$.
\end{proof}

Under the above assumptions, we have the following lemma:

\begin{lemma}\label{lemmagradedexactsequence}
For any integer \( r \), we have the following exact sequence:
\[
\begin{tikzcd}
0 \arrow[r] & \bigoplus^{m_{k}}_{l=1} \widehat{P}^{r-\mathrm{a}^{(k)}_l} \arrow[r, "F_k"] & \cdots \arrow[r, "F_2"] & \bigoplus^{m_{1}}_{l=1} \widehat{P}^{r-\mathrm{a}^{(1)}_l} \arrow[r, "F_1"] & \bigoplus^{m_0}_{l=1} \widehat{P}^{r-\mathrm{a}_l} \arrow[r, "F_0"] & \widehat{P}^{r} \arrow[r] & \widehat{P}^{\,r}/(I\cap\widehat{P}^{\,r}) \arrow[r] & 0
\end{tikzcd}
\]
Here we set $\widehat{P}^{\,i}:=\widehat{R}$ for all $i\leq 0$.
\end{lemma}
\begin{proof}
Clearly,
\[
\mathrm{Im}\big(F_2|_{\bigoplus_l \widehat{P}^{\,r-\mathrm{a}^{(2)}_l}}\big) \subset \mathrm{Ker}(F_1) \cap \bigoplus_l \widehat{P}^{\,r-\mathrm{a}^{(1)}_l}.
\]
Conversely, let $t \in \widehat{R}^{m_2}$ satisfy $F_2(t) \in \bigoplus_l \widehat{P}^{\,r-\mathrm{a}^{(1)}_l}$.
Let $[t]_k$ be the homogeneous component of minimal order $k$. Then using lemma \ref{lemmaleadingregular}
\[
L(F_2)([t]_k) = 0 \implies [t]_k = L(F_3)([s]_k) \text{ for some } [s]_k \in \widehat{R}^{m_3}.
\]
Set $t' := t - F_3([s]_k)$. Then $F_2(t) = F_2(t')$ and $\mathrm{ord}(t') > \mathrm{ord}(t)$. Iterating finitely yields $t^\infty \in \bigoplus_l \widehat{P}^{\,r-\mathrm{a}^{(2)}_l}$ with $F_2(t) = F_2(t^\infty)$.
Hence
\[
\mathrm{Im}\big(F_2|_{\bigoplus_l \widehat{P}^{\,r-\mathrm{a}^{(2)}_l}}\big) = \mathrm{Ker}(F_1) \cap \bigoplus_l \widehat{P}^{\,r-\mathrm{a}^{(1)}_l}.
\]
The exactness of the other terms follows similarly.
\end{proof}

\paragraph{Freedom of resolutions.}
Finally, we note the flexibility in constructing such resolutions. Let
\[
\widehat{R} = \mathrm{k}[[x_1,\dots,x_n]] \cong \mathrm{k}[[x'_1,\dots,x'_n]],
\]
with induced isomorphism \(\alpha\). A resolution constructed via \(\mathrm{k}[[x'_i]]\) is mapped to \(\mathrm{k}[[x_i]]\) by \(\alpha\), and compared with the resolution built directly in \(\mathrm{k}[[x_i]]\). They are homotopy equivalent, yielding a commutative diagram of resolutions:

\[
\begin{tikzcd}[row sep=normal, column sep=normal]
0 \arrow[r] &
\widehat{R}^{m_k} \arrow[r, "F_k"] \arrow[d, "\phi_k"] &
\cdots \arrow[r, "F_2"] &
\widehat{R}^{m_1} \arrow[r, "F_1"] \arrow[d, "\phi_1"] &
\widehat{R}^{m_0} \arrow[r, "F_0"] \arrow[d, "\phi_0"] &
\widehat{R} \arrow[r] \arrow[d, "\mathrm{id}"] &
\widehat{R}/I \arrow[r] \arrow[d, "\mathrm{id}"] &
0 \\
0 \arrow[r] &
\widehat{R}^{n_k} \arrow[r, "G_k"] \arrow[d, "\varphi_k"] &
\cdots \arrow[r, "G_2"] &
\widehat{R}^{n_1} \arrow[r, "G_1"] \arrow[d, "\varphi_1"] &
\widehat{R}^{n_0} \arrow[r, "G_0"] \arrow[d, "\varphi_0"] &
\widehat{R} \arrow[r] \arrow[d, "\mathrm{id}"] &
\widehat{R}/I \arrow[r] \arrow[d, "\mathrm{id}"] &
0 \\
0 \arrow[r] &
\widehat{R}^{m_k} \arrow[r, "F_k"] &
\cdots \arrow[r, "F_2"] &
\widehat{R}^{m_1} \arrow[r, "F_1"] &
\widehat{R}^{m_0} \arrow[r, "F_0"] &
\widehat{R} \arrow[r] &
\widehat{R}/I \arrow[r] &
0
\end{tikzcd}
\]
where there exist homotopy maps $h_i: \widehat{R}^{m_i} \to \widehat{R}^{m_{i+1}}$ and $k_i: \widehat{R}^{n_i} \to \widehat{R}^{n_{i+1}}$ such that:
\[
\varphi_i \circ \phi_i - \mathrm{id} = F_{i+1} \circ h_i + h_{i-1} \circ F_i
\]
\[
\phi_i \circ \varphi_i - \mathrm{id} = G_{i+1} \circ k_i + k_{i-1} \circ G_i
\]

\begin{lemma}\label{gradedfiterexact}
We claim that for all integers $r$, there are maps $\phi$ and $\varphi$ induce the following homotopy-commutative diagram via $h^\infty$ and $k^\infty$, for any resolution constructed in Lemma~\ref{lemmagradedexactsequence}, under different choices of coordinates.
\[
\begin{tikzcd}[row sep=normal, column sep=normal]
0 \arrow[r] &
\bigoplus\limits_{l=1}^{m_k} \widehat{P}^{r-\mathrm{a}^{(k)}_l} \arrow[r, "F_k"] \arrow[d, "\phi_k"] &
\cdots \arrow[r, "F_2"] &
\bigoplus\limits_{l=1}^{m_1} \widehat{P}^{r-\mathrm{a}^{(1)}_l} \arrow[r, "F_1"] \arrow[d, "\phi_1"] &
\bigoplus\limits_{l=1}^{m_0} \widehat{P}^{r-\mathrm{a}_l} \arrow[r, "F_0"] \arrow[d, "\phi_0"] &
\widehat{P}^{r} \arrow[r] \arrow[d, "\mathrm{id}"] &
\widehat{P}^{\,r}/(I\cap\widehat{P}^{\,r}) \arrow[r] \arrow[d, "\mathrm{id}"] &
0 \\
0 \arrow[r] &
\bigoplus\limits_{l=1}^{n_k} \widehat{P}^{r-\mathrm{b}^{(k)}_l} \arrow[r, "G_k"] \arrow[d, "\varphi_k"] &
\cdots \arrow[r, "G_2"] &
\bigoplus\limits_{l=1}^{n_1} \widehat{P}^{r-\mathrm{b}^{(1)}_l} \arrow[r, "G_1"] \arrow[d, "\varphi_1"] &
\bigoplus\limits_{l=1}^{n_0} \widehat{P}^{r-\mathrm{b}_l} \arrow[r, "G_0"] \arrow[d, "\varphi_0"] &
\widehat{P}^{r} \arrow[r] \arrow[d, "\mathrm{id}"] &
\widehat{P}^{\,r}/(I\cap\widehat{P}^{\,r}) \arrow[r] \arrow[d, "\mathrm{id}"] &
0 \\
0 \arrow[r] &
\bigoplus\limits_{l=1}^{m_k} \widehat{P}^{r-\mathrm{a}^{(k)}_l} \arrow[r, "F_k"] &
\cdots \arrow[r, "F_2"] &
\bigoplus\limits_{l=1}^{m_1} \widehat{P}^{r-\mathrm{a}^{(1)}_l} \arrow[r, "F_1"] &
\bigoplus\limits_{l=1}^{m_0} \widehat{P}^{r-\mathrm{a}_l} \arrow[r, "F_0"] &
\widehat{P}^{r} \arrow[r] &
\widehat{P}^{\,r}/(I\cap\widehat{P}^{\,r}) \arrow[r] &
0
\end{tikzcd}
\]
\end{lemma}
\begin{proof}
To demonstrate the well-definedness of the morphism
\[
\phi_0: \bigoplus_{l=1}^{m_0} \widehat{P}^{r - a_l} \longrightarrow \bigoplus_{l=1}^{n_0} \widehat{P}^{r - b_l},
\]
we consider the relation
\[
\phi_0 \circ G_0 = F_0,
\]
where \( G_0 \) is a matrix whose columns form a standard basis in the presentation \( k[[x'_1, \ldots, x'_n]] \). Since \( G_0 \) consists of standard basis elements, each column of \( F_0 \) can be expressed as a combination of these basis elements. This implies that the entries \( \phi_{0, i,j} \) of the matrix \( \phi_0 \) satisfy the inequality
\[
{\mathrm{ord}_{\widehat{P}}}(\phi_{0, i,j}) \geq a_i - b_j,
\]
for all indices \( i \) and \( j \). Therefore, the morphism \( \phi_0 \) is well-defined between the respective graded components.\\

To demonstrate the existence of homotopy, we assume \( h_{-2} = h_{-1} = 0=h^{\infty}_{-1} = h^{\infty}_{-2} \), and we have:
\[
\varphi_0 \circ \phi_0 - \mathrm{id} = F_1 \circ h_0,
\]
For any \( m \in \bigoplus\limits_{l=1}^{m_0} \widehat{P}^{r-\mathrm{a}_l} \) by construction of \(\phi_0\) and \(\varphi_0\):
\[
(\varphi_0 \circ \phi_0 - \mathrm{id})(m) = F_1 \circ h_0(m) \in \bigoplus\limits_{l=1}^{m_0} \widehat{P}^{r-\mathrm{a}_l}.
\]
noticing the exact sequence under presentation \( k[[x_1, \ldots, x_n]] \) and \( F_0 \circ F_1 = 0 \) , there exists \( [h_0(m)]^{\infty} \) in \(\bigoplus\limits_{l=1}^{m_1} \widehat{P}^{r-\mathrm{a}^{(1)}_l}\) as in the argument of Lemma~\ref{lemmagradedexactsequence}, such that:

\[
F_1 \circ h_0(m) = F_1 \circ [h_0(m)]^{\infty},
\]
which defines a map:
\[
h_0^{\infty} \colon m \mapsto [h_0(m)]^{\infty}.
\]
This map is a module homomorphism, and we have:
\[
\varphi_0 \circ \phi_0 - \mathrm{id} = F_1 \circ h_0^{\infty},
\]
when restricted to \( \bigoplus\limits_{l=1}^{m_0} \widehat{P}^{r-\mathrm{a}_l} \). The other terms follow similarly.
\end{proof}
A particular case occurs when different lexicographical orderings of \(x_1,\dots,x_n\) are used. For the reduced resolution, although a homotopy equivalence exists, it may fail to induce a well-behaved graded version; for example, the ideal
\(I = (x^2 + y^3, xy^2, y^5)\) with \(\mathrm{ord}(x)=\mathrm{ord}(y)=1\) and \(x<y\).
Resolutions constructed as above on \(R\) (or \(\widehat{R}\)), up to (graded) homotopy, are called \textbf{standard resolutions} of \(I\) (or \(R/I\)) on \(R\) (or \(\widehat{R}\)).

\subsection{Standard resolution on hypersurface}\label{subsec:hyp}
Let \( (R, \mathfrak{m}) \) be a regular local ring, and let \( J \) be an ideal such that \( S := R/J \) is normal Gorenstein. Suppose \( P \) is a prime ideal in \( R \) containing \( J \) such that \( Z := \operatorname{Spec}(R/P) \) is regular. We first consider the case where \( J \) defines a hypersurface \(w\). In general, we assume that $\operatorname{ord}_\mathfrak{m}(w) > 1$.
Furthermore, let \( I \) be an arbitrary ideal in \( R \) such that \( I \supseteq J \). Geometrically, we have the following commutative diagram:
\[
\begin{tikzcd}
\operatorname{Spec} R/P  \arrow[d, "k_1"] \arrow[dr, "k"] & \\
\operatorname{Spec} R/J \arrow[r, "j"] & \operatorname{Spec} R  \\
\operatorname{Spec} R/I \arrow[u, "i_1"] \arrow[ur, "i"] &
\end{tikzcd}
\hspace{1cm}\text{or equivalently}\hspace{1cm}
\begin{tikzcd}
Z \arrow[d, "k_1"] \arrow[dr, "k"] & \\
Y \arrow[r, "j"] & X  \\
B \arrow[u, "i_1"] \arrow[ur, "i"] &
\end{tikzcd}
\]
We consider two types of excess distinguished triangles.
From the exact sequence
\[
0 \longrightarrow \mathcal{O}_{X}(-Y) \xrightarrow{\cdot w} \mathcal{O}_{X} \longrightarrow \mathcal{O}_{Y} \longrightarrow 0,
\]

\begin{lemma}[\cite{HuybrechtsFM}, Cor. 11.4]
If \( Y \) is smooth, then for any \( F \in \mathrm{D}^b(Y) \), there exists a distinguished triangle:
\[
F \otimes \mathcal{O}_Y(-Y)[1] \longrightarrow j^{*}j_{*}F \longrightarrow F \quad \in \Delta.
\]
Moreover, applying this distinguished triangle to any other distinguished triangle results in a commutative diagram in the triangulated category.
\end{lemma}

This property need not hold for a singular hypersurface \(Y\); here we restrict to the basic case, for which it holds for any scheme \(X\) and any Cartier hypersurface \(Y\).

\begin{lemma}[\cite{HuybrechtsFM}, Cor. 11.2]
For any \( Y \) and any \( F \in \mathrm{Coh}(Y) \), there exists a unique distinguished triangle in \( \mathrm{D}^b(Y) \):
\[
F \otimes \mathcal{O}_Y(-Y)[1] \longrightarrow j^{*}j_{*}F \longrightarrow F \quad \in \Delta.
\]
\end{lemma}

\begin{proof}
We aim to prove that
\[
\mathcal{H}^\ell (j^* j_* F) = F (-Y)[1] \oplus F.
\]
Since \(j_*\) commutes with taking cohomology of complexes, this is equivalent to proving:
\[
j_* \mathcal{H}^\ell (j^* j_* F) = \mathcal{H}^\ell (j_* j^* j_* F) = j_* F (-Y)[1] \oplus j_* F.
\]
This follows immediately, as the functorial properties of \( j_* j^* j_* \) induce this direct sum decomposition.
\end{proof}

\begin{corollary}
For any positive integer \( r \) and \( F \in \mathrm{Coh}(Y) \), we have the following left filtration:

\[
\begin{tikzcd}[column sep=small, row sep=small]
 & F \arrow[rr] && F(-Y)[2] \arrow[dl] \arrow[rr] && F(-2Y)[4] \arrow[dl] \\
 & & j^{*}j_{*}F[1] \arrow[ul, dashed, "\Delta"] && j^{*}j_{*}F(-Y)[3] \arrow[ul, dashed, "\Delta"] & \cdots
\end{tikzcd}
\]

\[
\begin{tikzcd}[column sep=small, row sep=small]
 & F(-(r-1)Y)[2(r-1)] \arrow[rr] && F(-rY)[2r] \arrow[dl] \\
 & & j^{*}j_{*}F(-(r-1)Y)[2r-1] \arrow[ul, dashed, "\Delta"]
\end{tikzcd}
\]

and the right Postnikov filtration of \( F \):
\[
\begin{tikzcd}[column sep=small, row sep=small]
 & F(rY)[-2r] \arrow[rr] && F((r-1)Y)[-2(r-1)] \arrow[dl] \\
 & & j^{*}j_{*}F((r-1)Y)[-(2r-1)] \arrow[ul, dashed, "\Delta"] & \cdots
\end{tikzcd}
\]

\[
\begin{tikzcd}[column sep=small, row sep=small]
 & F(2Y)[-4] \arrow[rr] && F(Y)[-2] \arrow[dl] \arrow[rr] && F \arrow[dl] \\
 & & j^{*}j_{*}F(Y)[-3] \arrow[ul, dashed, "\Delta"] && j^{*}j_{*}F[-1] \arrow[ul, dashed, "\Delta"]
\end{tikzcd}
\]
\end{corollary}

\begin{corollary}\label{excesshyper}
By applying the octahedral axiom in the triangulated category to the above filtrations, we obtain the following distinguished triangle:
\[
F \longrightarrow F(-rY)[2r] \longrightarrow F_r \longrightarrow F[1],
\]
where \( F_r \) is constructed via a series of extensions from \( j^{*}j_{*}F[1] \) to \( j^{*}j_{*}F(-(r-1)Y)[2r-1] \). For example, for \( F_2 \), we have the following distinguished triangle:
\[
j^{*}j_{*}F[1] \longrightarrow F_2 \longrightarrow j^{*}j_{*}F(-Y)[3] \longrightarrow j^{*}j_{*}F[2].
\]
and more generally, \( F_r \) admits the following filtration:
\[
\begin{tikzcd}[column sep=small, row sep=small]
 & 0 \arrow[rr] && F_1 \arrow[dl] \arrow[rr] && F_2 \arrow[dl] \\
 & & j^{*}j_{*}F[1] \arrow[ul, dashed, "\Delta"] && j^{*}j_{*}F(-Y)[3] \arrow[ul, dashed, "\Delta"] &
\end{tikzcd}\cdots
\begin{tikzcd}[column sep=small, row sep=small]
 & F_{r-1} \arrow[rr] && F_r \arrow[dl] \\
 & & j^{*}j_{*}F(-(r-1)Y)[2r-1] \arrow[ul, dashed, "\Delta"]
\end{tikzcd}
\]
\end{corollary}

We have the following simple cohomology calculation:
\begin{lemma}\label{lemmaexccheck}
\[
\mathrm{Ext}^1 \left( j^{*}j_{*} F(- (n-1) Y) [2n-1], j^{*}j_{*} F (-(m-1) Y) [2m-1] \right)
=
\begin{cases}
\mathrm{Hom} \left( F,  F \right) & n = m+1 \\
0 & n > m+1
\end{cases}
\]
\end{lemma}

\begin{proof}
We have:
\[
\begin{aligned}
&\mathrm{Ext}^1 \left( j^{*}j_{*} F(- (n-1) Y) [2n-1], j^{*}j_{*} F (-(m-1) Y) [2m-1] \right) \\
&\quad = \mathrm{Ext}^{2m-2n+1} \left( j^{*}j_{*} F(- (n-1) Y), j^{*}j_{*} F (-(m-1) Y) \right) \\
&\quad = \mathrm{Ext}^{2m-2n+1} \left( j_{*} F(- (n-1) Y), j_{*}j^{*}j_{*} F (-(m-1) Y) \right) \\
&\quad = \mathrm{Ext}^{2m-2n+1} \left( j_{*} F(- (n-1) Y), j_{*} F (-(m-1) Y) \oplus j_{*} F (-mY) [1] \right).
\end{aligned}
\]
If we require \( n \geq m+1 \), then the only non-trivial case is \( n = m+1 \). In this case, the above equality becomes:
\[
\mathrm{Ext}^{-1} \left( j_{*} F(-mY), j_{*} F (-(m-1) Y) \oplus j_{*} F (-mY) [1] \right)
= \mathrm{Hom} \left( j_{*} F, j_{*} F \right)= \mathrm{Hom} \left( F, F \right).
\]
\end{proof}

Therefore, based on simple inductive calculations and observing the symmetry in the extensions at each step of the filtration for \( F_k \), we have:
\begin{corollary}\label{corprim}
For any positive integer $r$ and $F\in\operatorname{Coh}(Y)$, the filtration extension of $F_r$ is completely determined by the classes
\[
\operatorname{Ext}^{-1}\bigl(j^{*}j_{*}F(-mY),j^{*}j_{*}F(-(m-1)Y)\bigr),
\qquad 1\le m\le r.
\]
Moreover, there exists a primitive class $s(F):=\operatorname{id}_{F}\in\operatorname{Hom}(F,F)$, unique up to scaling, which determines the extension.
For example, for $F_2$ there is a distinguished triangle
\[
j^{*}j_{*}F[1]\to F_2\to j^{*}j_{*}F(-Y)[3]\to j^{*}j_{*}F[2].
\]
The corresponding extension class is in
\[
\operatorname{Ext}^{-1}\bigl(j^{*}j_{*}F(-Y),j^{*}j_{*}F\bigr)
\cong
\operatorname{Hom}(F,F),
\]
and is represented by $\operatorname{id}_{F}$. In particular, this primitive class cannot be written as the composition of a nontrivial element of
$\operatorname{Ext}^{-1}\bigl(j^{*}j_{*}F(-Y),j^{*}j_{*}F\bigr)$
with a nontrivial element of
$\operatorname{Hom}\bigl(j^{*}j_{*}F,j^{*}j_{*}F\bigr)\cong\operatorname{Hom}(F,F)$.
\end{corollary}
\begin{proof}
According to the construction, our extension is given by the composition by \cite[Cor. 11.4]{HuybrechtsFM}:
\[
j^{*}j_{*}F(-mY)[1] \longrightarrow F(-mY)[1] \longrightarrow j^{*}j_{*}F(-(m-1)Y),
\]
where the first map corresponds to the identity map in \(\mathrm{End}(j_{*}F(-mY))\) via trivial duality, and the second map corresponds to the identity map in \(\mathrm{End}(j_{*}F(-(m-1)Y))\) via Grothendieck–Verdier duality. Therefore, based on the previous cohomology calculations, their composition remains trivial under these correspondences.
\end{proof}

We consider the following simple example: let  \(\mathcal{O}_{B} \) be the structure sheaf of an arbitrary non-trivial cycle \( B \) on \( Y \). Specifically, we consider the standard resolution of \( j_{*}\mathcal{O}_{B} \) over \( \mathcal{O}_{X} \):
\[
\begin{tikzcd}
0 \arrow[r] & \mathcal{O}_{X}^{m_k} \arrow[r, "F_k"] & \mathcal{O}_{X}^{m_{k-1}} \arrow[r, "F_{k-1}"] & \cdots \arrow[r, "F_2"] & \mathcal{O}_{X}^{m_1} \arrow[r, "F_1"] & \mathcal{O}_{X}^{m_0} \arrow[r, "F_0"] & \mathcal{O}_{X} \arrow[r] & j_{*}\mathcal{O}_{B} \arrow[r] & 0
\end{tikzcd}
\]
Then, pulling back along \( j \), we obtain a quasi-isomorphism of its projective presentation:
\[
\begin{tikzcd}
0 \arrow[r] & \mathcal{O}_{Y}^{m_k} \arrow[r, "F_k"] & \mathcal{O}_{Y}^{m_{k-1}} \arrow[r, "F_{k-1}"] & \cdots \arrow[r, "F_2"] & \mathcal{O}_{Y}^{m_1} \arrow[r, "F_1"] & \mathcal{O}_{Y}^{m_0} \arrow[r, "F_0"] & \mathcal{O}_{Y} \arrow[r,"\cong"] & j^*j_*\mathcal{O}_{B} \arrow[r] & 0
\end{tikzcd}
\]
Consider the extension given by \( s(\mathcal{O}_{B}) \) in \( \mathrm{Ext}^{-1}(j^{*}j_{*}\mathcal{O}_{B}(-Y), j^{*}j_{*}\mathcal{O}_{B})= \mathcal{O}_{B}\), as constructed above. If we have the corresponding relative projective resolution, e.g., \(Y\) is affine and hence \( \mathrm{Ext}^{>0}(\mathcal{O}_{Y}, \mathcal{O}_{Y}) = 0\), it induces a map between the corresponding resolutions:
\[
\begin{tikzcd}[column sep=small, row sep=small]
\mathcal{O}_{Y}^{m_k}(-Y) \arrow[r] & \cdots \arrow[r] \arrow[d] & \mathcal{O}_{Y}^{m_1}(-Y) \arrow[r] \arrow[d] & \mathcal{O}_{Y}^{m_0}(-Y) \arrow[r] \arrow[d] & \mathcal{O}_{Y}(-Y) \arrow[r,"\cong"] \arrow[d] & j^*j_*\mathcal{O}_{B}(-Y)  \\
& \mathcal{O}_{Y}^{m_k} \arrow[r] & \cdots \arrow[r] & \mathcal{O}_{Y}^{m_1} \arrow[r] & \mathcal{O}_{Y}^{m_0} \arrow[r] & \mathcal{O}_{Y} \arrow[r,"\cong"] & j^*j_*\mathcal{O}_{B}
\end{tikzcd}
\]
Moreover, the mapping cone or totalization of such a map provides a projective presentation of \( \mathcal{O}_{B,2} \).

\[
\begin{tikzcd}
\mathcal{O}_{Y}^{m_k}(-Y) \arrow[r, "F_{k+2,2}"] & \mathcal{O}_{Y}^{m_{k-1}}(-Y) \arrow[r, "F_{k+1,2}"] & \cdots\arrow[r, "F_{2,2}"] & \mathcal{O}_{Y}^{m_1}\oplus \mathcal{O}_{Y}(-Y) \arrow[r, "F_{1,2}"] & \mathcal{O}_{Y}^{m_0} \arrow[r, "F_{0,2}"] & \mathcal{O}_{Y} \arrow[r, "\cong"] & \mathcal{O}_{B,2}
\end{tikzcd}
\]
Similarly, we repeat the above steps while considering
\( \mathrm{Ext}^{2}(j^{*}j_{*}\mathcal{O}_{B}(-2Y), \mathcal{O}_{B,2})= \mathcal{O}_{B}\),
and continue this process for any sufficiently large \( r \). We obtain the commutative diagram:
\[
\begin{tikzcd}[column sep=small, row sep=small]
 & \vdots \arrow[d]& & & \vdots \arrow[d]  \\
\mathcal{O}_{Y}^{m_k}(-2Y) \arrow[r] & \cdots \arrow[r] \arrow[d] & \mathcal{O}_{Y}^{m_0}(-2Y) \arrow[r] \arrow[d] & \mathcal{O}_{Y}(-2Y) \arrow[ddl,dashed]\arrow[r,"\cong"] \arrow[d] & j^*j_*\mathcal{O}_{B}(-2Y) \arrow[d] \\
& \mathcal{O}_{Y}^{m_k}(-Y) \arrow[r] & \cdots \arrow[r] \arrow[d] & \mathcal{O}_{Y}^{m_0}(-Y) \arrow[r] \arrow[d] & \mathcal{O}_{Y}(-Y) \arrow[r,"\cong"] \arrow[d] & j^*j_*\mathcal{O}_{B}(-Y) \arrow[d] \\
& & \mathcal{O}_{Y}^{m_k} \arrow[r] & \cdots \arrow[r] & \mathcal{O}_{Y}^{m_0} \arrow[r] & \mathcal{O}_{Y} \arrow[r,"\cong"] & j^*j_*\mathcal{O}_{B}
\end{tikzcd}
\]
Then, by performing iterated totalization, we obtain a projective presentation for any \( \mathcal{O}_{B,r} \):
\[
\begin{tikzcd}
 \cdots \arrow[r, "F_{i+1,r}"] & \bigoplus_{j \mid -1 \leq i-2j \leq i} \mathcal{O}_{Y}^{m_{i-2j}}(-jY) \arrow[r, "F_{i,r}"] & \cdots \arrow[r, "F_{1,r}"] & \mathcal{O}_{Y}^{m_0} \arrow[r, "F_{0,r}"] & \mathcal{O}_{Y} \arrow[r, "\cong"] & \mathcal{O}_{B,r}
\end{tikzcd}
\]
We can easily observe that if \( r \) is sufficiently large, the rank of each term in this exact sequence increases in 2-periodic as the degree decreases, and stabilizes after degree \(-(k-1)\) (if we define the degree of \(\mathcal{O}_{Y}^{m_0} \) to be zero). The exact sequence forms a 2-periodic complex from degree \(-k\) up to the term of degree \(-k-r\), and after degree \(-k-r-1\) up to \(-2k-r+1\), the rank of each term decreases as the degree decreases in 2-periodic, eventually vanishing. In particular, we have the following description for the term of degree \(-(k-1)\):
\[
V_{k-1}:=\bigoplus_{j \mid -1 \leq k-2j-2 \leq k-2} \mathcal{O}_{Y}^{m_{k-2j-2}}(-jY).
\]
and term of degree \(-k\)
\[
V_{k}:=\bigoplus_{j \mid -1 \leq k-2j-1 \leq k-1} \mathcal{O}_{Y}^{m_{k-2j-1}}(-jY).
\]
Thus, we conclude that \( V_{k-1} \) and \( V_{k} \) are locally free sheaves of the same rank. Furthermore, we have that for \( 1 \leq 2t \leq r-1 \),
\[
V_{k+2t-1} = V_{k-1}(-tY), \quad \text{and} \quad V_{k+2t} = V_{k}(-tY),
\]
in these stable 2-periodic complexes, neither the terms nor the maps between them change as \( r \) increases.\\

Let us now explain how to construct these induced maps. This construction may be viewed as an Eisenbud--Shamash resolution with respect to a standard basis; compare with \cite[Theorem~3.1.3]{Avramov1998}. To do this, we  return to the analytic local ring cases i.e. on $\widehat{X}$:
\begin{enumerate}
  \item First we consider \( j^*j_*\mathcal{O}_{\widehat{B}}(-\widehat{Y})[1] \) to \( j^*j_*\mathcal{O}_{\widehat{B}} \) for the \( {\mathcal{O}_{\widehat{B}}} \)-family.
\[
\begin{tikzcd}[column sep=normal, row sep=normal]
& \cdots \arrow[r, "F_{3}"] & \mathcal{O}_{\widehat{Y}}^{m_2}(-\widehat{Y}) \arrow[r, "F_{2}"] \arrow[d, "-K_3"] & \mathcal{O}_{\widehat{Y}}^{m_1}(-\widehat{Y}) \arrow[r, "F_{1}"] \arrow[d, "K_2"] & \mathcal{O}_{\widehat{Y}}^{m_0}(-\widehat{Y}) \arrow[r, "F_{0}"] \arrow[d, "-K_1"] & \mathcal{O}_{\widehat{Y}}(-\widehat{Y}) \arrow[d, "K_0"] \\
& & \cdots \arrow[r, "F_{3}"] & \mathcal{O}_{\widehat{Y}}^{m_2} \arrow[r, "F_{2}"] & \mathcal{O}_{\widehat{Y}}^{m_1} \arrow[r, "F_{1}"] & \mathcal{O}_{\widehat{Y}}^{m_0} \arrow[r, "F_{0}"] & \mathcal{O}_{\widehat{Y}}
\end{tikzcd}
\]
Since \( \widehat{B} \) is contained in \( \widehat{Y} \), and note that
  \[
  F_0 = \begin{pmatrix}
  f^{(0)}_1 & f^{(0)}_2 & \cdots & f^{(0)}_{m_0}
  \end{pmatrix}
  \]
  is a standard basis for \( I \) in \( \widehat{X} \), and \( w \) is an element of \( I \). By Grauert's division theorem \cite[Thm. 7.1.9]{DeJongPfister}, there exists a unique division expression:
  \[
w = \sum_{l=1}^{m_0} f^{(0)}_l \cdot k^{(0)}_l + r^{(0)},
  \]
  where the remainder \( r^{(0)} = 0 \), and the following condition on the initial terms holds:
  \[
  \mathrm{in}(t_0 \cdot w) \leq \mathrm{in}(f^{(0)}_l \cdot k^{(0)}_l) \quad \text{for any } l.
  \]
  and \(F_0 \circ K_0 =  w\).

  \item Next, observing that on \( \widehat{R} \), the following equality holds:
  \[
  F_0 \circ (w-K_0 \circ F_0) = F_0 \circ w -F_0 \circ K_0 \circ F_0 = F_0 \circ w - w \circ F_0 = 0.
  \]
  This implies that \( w - K_0 \circ F_0 \) lies in the kernel of \( F_0 \), meaning its column vectors belong to \( \mathrm{Syz}(I) \).

  Consider
  \[
  F_1 = \begin{pmatrix}
  f^{(1)}_1 & f^{(1)}_2 & \cdots & f^{(1)}_{m_1}
  \end{pmatrix},
  \]
  a standard basis for \( \mathrm{Syz}(I) \) in \( \widehat{R} \) under the natural module order, where each \( f^{(1)}_i \) is an \( m_0 \times m_1 \) matrix. Applying the division theorem to each column of the right-hand side matrix yields:
  \[
  \sum_{l=1}^{m_1} f^{(1)}_l \cdot k^{(1)}_l + r^{(1)} = w - K_0 \circ F_0,
  \]
  with \( r^{(1)} = 0 \). Additionally, for each column, the initial terms satisfy:
  \[
  \mathrm{in}(w - K_0 \circ F_0) \leq \mathrm{in}(f^{(1)}_l \cdot k^{(1)}_l) \quad \text{for any } l.
  \]
  Thus, we obtain  \(F_1 \circ K_1 = w - K_0 \circ F_0\).
  \item Similarly, we consider the following equality on \( \widehat{R} \):
  \[
  F_1 \circ (w - K_1 \circ F_1) = F_1 \circ w - (F_1 \circ K_1) \circ F_1 = F_1 \circ w - (w - K_0 \circ F_0) \circ F_1 = 0.
  \]
  This shows that \( w - K_1 \circ F_1 \) lies in the kernel of \( F_1 \), meaning its column vectors belong to \( \mathrm{Syz}^2(I) \).

Recall that
  \[
  F_2 = \begin{pmatrix}
  f^{(2)}_1 & f^{(2)}_2 & \cdots & f^{(2)}_{m_2}
  \end{pmatrix},
  \]
  is a standard basis for \( \mathrm{Syz}^2(I) \) in \( \widehat{R} \) under the natural module order. Applying the division theorem similarly, we obtain:
  \[
  F_{2} \circ K_{2} = w - K_1 \circ F_1.
  \]

\item Repeating the above steps, we obtain for any \( i \), on \( \widehat{R} \):
\[
K_i \circ F_i = w - F_{i+1} \circ K_{i+1}.
\]
Considering the natural projection of such maps to \( \widehat{Y} \), a simple inspection shows that the morphism is not homotopy to zero which induces a non-trivial map:
\[
\begin{tikzcd}
j^*j_*{\mathcal{O}_{\widehat{B}}}(-\widehat{Y})[1] \arrow[r, "s({\mathcal{O}_{\widehat{B}}})"] & j^*j_*\mathcal{O}_{\widehat{B}}
\end{tikzcd}
\]
corresponds to the primitive extension given by the excess distinguished triangle in \ref{corprim}.

\item For any element \( t \) of \( {\mathcal{O}_{\widehat{B}}} \), it induces a map:
\[
t: j^*j_*\mathcal{O}_{\widehat{B}} \longrightarrow j^*j_*\mathcal{O}_{\widehat{B}}.
\]
We consider the following composition:
\[
\begin{tikzcd}
j^*j_*\mathcal{O}_{\widehat{B}}(-\widehat{Y})[1] \arrow[r,"s({\mathcal{O}_{\widehat{B}}})"] & j^*j_*\mathcal{O}_{\widehat{B}} \arrow[r, "t"] & j^*j_*\mathcal{O}_{\widehat{B}}.
\end{tikzcd}
\]
Based on cohomology calculations and primitive, we obtain all the mappings in \ref{corprim}.

\item Second, we consider \( j^*j_*\mathcal{O}_{\widehat{B}}(-2\widehat{Y})[2] \) to  \(\mathcal{O}_{\widehat{B},2}\) which is a \( {\mathcal{O}_{\widehat{B}}} \)-family.

By symmetry, we construct the $K^{(1)}_i$ matrices similar to the previous steps, then observe that:
\[
F_1 \circ (K_1 \circ K_0) = (f - K_0 \circ F_0)\circ K_0 = fK_0 - K_0f = 0.
\]
Since \( F_1 \) admits a natural standard basis, there exists a matrix \( K_0^{(1)} \) such that
\[
K_1 \circ K_0 + F_2 \circ K_0^{(1)}=0.
\]

\item Next, we proceed by induction to show that for any non-negative integer \( 0 \leq i \leq k-3 \), there exist matrices \( K_i^{(1)} \) such that
\[
K_{i+1} \circ K_i + F_{i+2} \circ K_i^{(1)} + K_{i-1}^{(1)} \circ F_{i-1} = 0.
\]
The computation shows:
\begin{align*}
F_{i+1} \circ ( K_{i+1} \circ K_i+ K_{i-1}^{(1)} \circ F_{i-1})
&= (w - K_i \circ F_i) \circ K_i +(F_{i+1} \circ K_{i-1}^{(1)} )\circ F_{i-1}\\
&= w \circ K_i - K_i \circ (F_i \circ K_i)+(F_{i+1} \circ K_{i-1}^{(1)} )\circ F_{i-1} \\
&= w \circ K_i - K_i \circ (w - K_{i-1} \circ F_{i-1})+(F_{i+1} \circ K_{i-1}^{(1)} )\circ F_{i-1} \\
&= K_i \circ K_{i-1} \circ F_{i-1}+(F_{i+1} \circ K_{i-1}^{(1)} )\circ F_{i-1}\\
&= K_i \circ K_{i-1} \circ F_{i-1}- (K_{i} \circ K_{i-1} + K_{i-2}^{(1)} \circ F_{i-2})\circ F_{i-1}\\
&= 0
\end{align*}
The last step uses the induction hypothesis at step \( i-1 \). Since \( \mathrm{Im}\,F_{i+2} \) admits a natural standard basis, by construction via division we can define a matrix \( K_i^{(1)} \) such that
\[
K_{i+1} \circ K_i + F_{i+2} \circ K_i^{(1)} + K_{i-1}^{(1)} \circ F_{i-1} = 0.
\]

\[
\begin{tikzcd}[column sep=small, row sep=small]
\cdots \arrow[r, "F_{3}"] & \mathcal{O}_{\widehat{X}}^{m_2}(-2\widehat{Y}) \arrow[r, "F_{2}"] \arrow[d, "K_3"] & \mathcal{O}_{\widehat{X}}^{m_1}(-2\widehat{Y}) \arrow[r, "F_{1}"] \arrow[d, "-K_2"] & \mathcal{O}_{\widehat{X}}^{m_0}(-2\widehat{Y}) \arrow[r, "F_{0}"] \arrow[d, "K_1"]\arrow[ddl,dashed, "K^{(1)}_1"{pos=.6}] & \mathcal{O}_{\widehat{X}}(-2\widehat{Y}) \arrow[d, "-K_0"]\arrow[ddl,dashed, "-K^{(1)}_0"{pos=.6}] \\
& \cdots \arrow[r, "F_{3}"] & \mathcal{O}_{\widehat{X}}^{m_2}(-\widehat{Y}) \arrow[r, "F_{2}"] \arrow[d, "-K_3"] & \mathcal{O}_{\widehat{X}}^{m_1}(-\widehat{Y}) \arrow[r, "F_{1}"] \arrow[d, "K_2"] & \mathcal{O}_{\widehat{X}}^{m_0}(-\widehat{Y}) \arrow[r, "F_{0}"] \arrow[d, "-K_1"] & \mathcal{O}_{\widehat{X}}(-\widehat{Y}) \arrow[d, "K_0"] \\
& & \cdots \arrow[r, "F_{3}"] & \mathcal{O}_{\widehat{X}}^{m_2} \arrow[r, "F_{2}"] & \mathcal{O}_{\widehat{X}}^{m_1} \arrow[r, "F_{1}"] & \mathcal{O}_{\widehat{X}}^{m_0} \arrow[r, "F_{0}"] & \mathcal{O}_{\widehat{X}}
\end{tikzcd}
\]
It is not difficult to see that the above commutative diagram is not homotopy trivial since \(F_0\circ K_0=w\), this construction provides the primitive extension, which together with its scalings determines all the morphisms from $j^*j_*\mathcal{O}_{\widehat{B}}(-2\widehat{Y})[2]$ to $\mathcal{O}_{B,2}$.

\item Thirdly, we consider \( j^*j_*\mathcal{O}_{\widehat{B}}(-3\widehat{Y})[3] \) to  \(\mathcal{O}_{B,3}\) which is also a \( {\mathcal{O}_{\widehat{B}}} \)-family. Similarly, we have constructed matrices $K_i$, $K^{(1)}_i$, $K^{(2)}_i$ and the following commutative diagram:
\[
\begin{tikzcd}[column sep=small, row sep=large]
\cdots \arrow[r, "F_{3}"] & \mathcal{O}_{\widehat{X}}^{m_2}(-3\widehat{Y}) \arrow[r, "F_{2}"] \arrow[d, "-K_3"] & \mathcal{O}_{\widehat{X}}^{m_1}(-3\widehat{Y})\arrow[ddl,dashed] \arrow[r, "F_{1}"] \arrow[d, "K_2"] & \mathcal{O}_{\widehat{X}}^{m_0}(-3\widehat{Y}) \arrow[dddll,dashed] \arrow[r, "F_{0}"] \arrow[d, "-K_1"]\arrow[ddl,dashed, "-K^{(1)}_1"{pos=.6}] & \mathcal{O}_{\widehat{X}}(-3\widehat{Y}) \arrow[d, "K_0"]\arrow[ddl,dashed, "K^{(1)}_0"{pos=.6}] \arrow[dddll,dashed, "K^{(2)}_0"{above,pos=.1}]&\\
&\cdots\arrow[d] \arrow[r, "F_{3}"] & \mathcal{O}_{\widehat{X}}^{m_2}(-2\widehat{Y}) \arrow[ddl,dashed]\arrow[r, "F_{2}"] \arrow[d, "K_3"] & \mathcal{O}_{\widehat{X}}^{m_1}(-2\widehat{Y}) \arrow[ddl,dashed]\arrow[r, "F_{1}"] \arrow[d, "-K_2"] & \mathcal{O}_{\widehat{X}}^{m_0}(-2\widehat{Y}) \arrow[r, "F_{0}"] \arrow[d, "K_1"]\arrow[ddl,dashed, "K^{(1)}_1"{pos=.6}] & \mathcal{O}_{\widehat{X}}(-2\widehat{Y}) \arrow[d, "-K_0"]\arrow[ddl,dashed, "-K^{(1)}_0"{pos=.6}] \\
&\cdots\arrow[d] \arrow[r, "F_{4}"]& \mathcal{O}_{\widehat{X}}^{m_3}\arrow[r, "F_{3}"]\arrow[d, "K_{4}"] & \mathcal{O}_{\widehat{X}}^{m_2}(-\widehat{Y}) \arrow[r, "F_{2}"] \arrow[d, "-K_3"] & \mathcal{O}_{\widehat{X}}^{m_1}(-\widehat{Y}) \arrow[r, "F_{1}"] \arrow[d, "K_2"] & \mathcal{O}_{\widehat{X}}^{m_0}(-\widehat{Y}) \arrow[r, "F_{0}"] \arrow[d, "-K_1"] & \mathcal{O}_{\widehat{X}}(-\widehat{Y}) \arrow[d, "K_0"] \\
&\cdots\arrow[r, "F_{5}"] &\mathcal{O}_{\widehat{X}}^{m_4}\arrow[r, "F_{4}"] & \mathcal{O}_{\widehat{X}}^{m_3}\arrow[r, "F_{3}"] & \mathcal{O}_{\widehat{X}}^{m_2} \arrow[r, "F_{2}"] & \mathcal{O}_{\widehat{X}}^{m_1} \arrow[r, "F_{1}"] & \mathcal{O}_{\widehat{X}}^{m_0} \arrow[r, "F_{0}"] & \mathcal{O}_{\widehat{X}}
\end{tikzcd}
\]
The constructed matrices satisfy the compatibility condition:
\[
K_{i+1}^{(1)} \circ K_i + K_{i+3} \circ K_i^{(1)} + F_{i+4} \circ K_i^{(2)} + K_{i-1}^{(2)} \circ F_{i-1} = 0
\]
\[
K_{i+1} \circ K_i + F_{i+2} \circ K_i^{(1)} + K_{i-1}^{(1)} \circ F_{i-1} = 0.
\]
for all relevant indices $i$.
\item By induction, we prove that for any non-negative integer $a$ and integer $i$, we can construct matrices $K_i^{(a)}$ (with the convention $K_i^{(0)} := K_i$) satisfying:
\[
\sum_{\substack{a_1+a_2=a-1}} K_{i+2a_2+1}^{(a_1)} \circ K_i^{(a_2)} + F_{i+2a} \circ K_i^{(a)} + K_{i-1}^{(a)} \circ F_{i-1} = 0
\]
where the sum is taken over all pairs of non-negative integers $(a_1,a_2)$ with $a_1+a_2=a-1$.

Assume we have constructed $K_{i'}^{(a')}$ for all pairs $(a',i')$ satisfying either:
\begin{enumerate}
\item $a' < a$ (lower order cases), or
\item $a' = a$ and $i' < i$ (same order but lower index cases)
\end{enumerate}

such that they fulfill the previous conditions. Then for the pair $(a,i)$, there exists $K_i^{(a)}$ satisfying the recurrence relation, since
\begin{align*}
F_{i+2a-1} \circ K_{i-1}^{(a)} \circ F_{i-1} &=
-\sum_{\substack{a_1+a_2=a-1}} \Big( K_{i+2a_2}^{(a_1)} \circ K_{i-1}^{(a_2)} + K_{i-2}^{(a)} \circ F_{i-2} \Big) \circ F_{i-1}\\
&=-\sum_{\substack{a_1+a_2=a-1}} K_{i+2a_2}^{(a_1)} \circ K_{i-1}^{(a_2)} \circ F_{i-1}
\end{align*}
and
\begin{align*}
F_{i+2a-1} \circ \left( \sum_{\substack{a_1+a_2=a-1}} K_{i+2a_2+1}^{(a_1)} \circ K_i^{(a_2)} \right)
&= \sum_{\substack{a_1+a_2=a-1}} \left(F_{i+2a+1} \circ K_{i+2a_2+1}^{(a_1)} \right)\circ K_i^{(a_2)} \\
&= \sum_{\substack{a_1+a_2=a-1}} \left(F_{i+2a_1} \circ K_{i}^{(a_1)} \right)_{2a_2+1}\circ K_i^{(a_2)}
\end{align*}
\begin{align*}
&= \sum_{\substack{a_1+a_2=a-1}} \left( - \sum_{\substack{a_{11}+a_{12}=a_1-1}} K_{i+2a_{12}+1}^{(a_{11})} \circ K_i^{(a_{12})} - K_{i-1}^{(a_1)} \circ F_{i-1} \right)_{2a_2+1} \circ K_i^{(a_2)} \\
&= -\sum_{\substack{a_1+a_2=a-1\\a_{11}+a_{12}=a_1-1}} K_{i+2a_{12}+1+2a_2+1}^{(a_{11})} \circ K_{i+2a_2+1}^{(a_{12})} \circ K_i^{(a_2)} \\
&\quad - \sum_{\substack{a_1+a_2=a-1}} K_{i+2a_2}^{(a_1)} \circ F_{i+2a_2} \circ K_i^{(a_2)}
\end{align*}
expanding the last term using recursion condition we have
\begin{align*}
&- \sum_{\substack{a_1+a_2=a-1}} K_{i+2a_2}^{(a_1)} \circ \left(F_{i+2a_2} \circ K_i^{(a_2)}\right) \\
&= - \sum_{\substack{a_1+a_2=a-1}} K_{i+2a_2}^{(a_1)} \circ \left(
- \sum_{\substack{a_{21}+a_{22}=a_2-1}} K_{i+2a_{22}+1}^{(a_{21})} \circ K_i^{(a_{22})}
- K_{i-1}^{(a_2)} \circ F_{i-1}
\right) \\
&= \sum_{\substack{a_1+a_2=a-1 \\ a_{21}+a_{22}=a_2-1}} K_{i+2a_2}^{(a_1)} \circ K_{i+2a_{22}+1}^{(a_{21})} \circ K_i^{(a_{22})} \\
&\quad + \sum_{\substack{a_1+a_2=a-1}} K_{i+2a_2}^{(a_1)} \circ K_{i-1}^{(a_2)} \circ F_{i-1}
\end{align*}
noticing
\begin{align*}
&\sum_{\substack{a_1+a_2=a-1 \\ a_{11}+a_{12}=a_1-1}}
K_{i+2a_{12}+2a_2+2}^{(a_{11})} \circ K_{i+2a_2+1}^{(a_{12})} \circ K_i^{(a_2)} \\
&= \sum_{\substack{a_1+a_2=a-1 \\ a_{21}+a_{22}=a_2-1}}
K_{i+2a_2}^{(a_1)} \circ K_{i+2a_{22}+1}^{(a_{21})} \circ K_i^{(a_{22})} = \sum_{\substack{b+c+d=a-2}}
K_{i+2c+2d+2}^{(b)} \circ K_{i+2d+1}^{(c)} \circ K_i^{(d)}
\end{align*}
So combine all above, the constructed matrices satisfy under induction assuption:
\[
F_{i+2a-1} \circ \left( \sum_{\substack{a_1 + a_2 = a-1}} K_{i+2a_2+1}^{(a_1)} \circ K_i^{(a_2)} + K_{i-1}^{(a)} \circ F_{i-1}\right) = 0
\]
This shows that \( \sum_{\substack{a_1 + a_2 = a-1}} K_{i+2a_2+1}^{(a_1)} \circ K_i^{(a_2)} + K_{i-1}^{(a)} \circ F_{i-1}\) lies in the image of \(F_{i+2a}\), recall that
  \[
  F_{i+2a} = \begin{pmatrix}
  f^{(i+2a)}_1 & f^{(i+2a)}_2 & \cdots & f^{(i+2a)}_{m_{i+2a}}
  \end{pmatrix},
  \]
  is a standard basis under the natural module order. Applying the division theorem, we obtain matrix $K_i^{(a)}$ such that:
\[
\sum_{\substack{a_1+a_2=a-1}} K_{i+2a_2+1}^{(a_1)} \circ K_i^{(a_2)} + F_{i+2a} \circ K_i^{(a)} + K_{i-1}^{(a)} \circ F_{i-1} = 0
\]
\item
Since the projective resolution over \(\widehat{R}\) is finite, there exists \(a_0\) such that for all \(a \ge a_0\) and all indices \(i\), \(K_i^{(a)} = 0\). Hence the inductive construction terminates in finitely many steps, producing only finitely many nonzero matrices. This construction relies solely on division and precedes cohomology computation.

\item We consider the following block matrices:
\[
A:= \begin{pmatrix}
K_0 & F_1 & & & \\
K_0^{(1)} & K_2 & F_3 & & \\
K_0^{(2)} & K_2^{(1)} & K_4 & \ddots & \\
\vdots & \vdots & \ddots & \ddots & F_{2k'-1} \\
K_0^{(k')} & K_2^{(k'-1)} & \cdots & K_{2k'-2}^{(1)} & K_{2k'}
\end{pmatrix}, \quad
B: =\begin{pmatrix}
F_0 &        &        &        &        \\
K_1 & F_2   &        &        &        \\
K_1^{(1)} & K_3 & F_4    &        &        \\
\vdots & \vdots & \ddots & \ddots &        \\
K_1^{(k'-1)} & K_3^{(k'-2)} & \cdots & K_{2k'-1} & F_{2k'}
\end{pmatrix}
\]
From the preceding discussion, we can directly verify through computation that for any $k'$ satisfying $2k' \leq k$, the following identity holds:
\[
A \circ B = B \circ A = w \cdot \mathrm{id}
\]
\item Using the cohomology lemma \ref{excesshyper}, we obtain a completely symmetric right projective resolution of $\mathcal{O}_{\widehat{B}}$ over $\widehat{Y}$.
\end{enumerate}
In particular, the above construction provides a unique matrix factorization $(A,B)$ representation in \( \widehat{X} \) for the stable periodic part of \( \mathcal{O}_{\widehat{B},r} \):
\[
\begin{tikzcd}[column sep=small, row sep=small]
& & \vdots \arrow[d] & & \vdots \arrow[d] \\
\cdots \arrow[r, "F_{3}"] & \mathcal{O}_{\widehat{X}}^{m_2}(-2\widehat{Y}) \arrow[r, "F_{2}"] \arrow[d, "K_3"] & \mathcal{O}_{\widehat{X}}^{m_1}(-2\widehat{Y}) \arrow[r, "F_{1}"] \arrow[d, "-K_2"] & \mathcal{O}_{\widehat{X}}^{m_0}(-2\widehat{Y}) \arrow[ddl,dashed]\arrow[r, "F_{0}"] \arrow[d, "K_1"] & \mathcal{O}_{\widehat{X}}(-2\widehat{Y}) \arrow[d, "-K_0"]  \arrow[ddl,dashed]\\
& \cdots \arrow[r, "F_{3}"] & \mathcal{O}_{\widehat{X}}^{m_2}(-\widehat{Y}) \arrow[r, "F_{2}"] \arrow[d, "-K_3"] & \mathcal{O}_{\widehat{X}}^{m_1}(-\widehat{Y}) \arrow[r, "F_{1}"] \arrow[d, "K_2"] & \mathcal{O}_{\widehat{X}}^{m_0}(-\widehat{Y}) \arrow[r, "F_{0}"] \arrow[d, "-K_1"] & \mathcal{O}_{\widehat{X}}(-\widehat{Y}) \arrow[d, "K_0"] \\
& & \cdots \arrow[r, "F_{3}"] & \mathcal{O}_{\widehat{X}}^{m_2} \arrow[r, "F_{2}"] & \mathcal{O}_{\widehat{X}}^{m_1} \arrow[r, "F_{1}"] & \mathcal{O}_{\widehat{X}}^{m_0} \arrow[r, "F_{0}"] & \mathcal{O}_{\widehat{X}}
\end{tikzcd}
\]
and the submatrices of $A$ and $B$ provide a complete representation of all double factorization totalization morphisms. We therefore obtain a projective resolution of $\mathcal{O}_{\widehat{B}}$ over $\widehat{Y}$ end with a matrix factorization $(A,B)$:
\[
\begin{tikzcd}[column sep=normal]
\cdots \arrow[r] &
\mathbf{V}_{k+2} \arrow[r, "\mathbf{F}_{k+2}"] &
\mathbf{V}_{k+1} \arrow[r, "\mathbf{F}_{k+1}"] &
\mathbf{V}_k \arrow[r, "\mathbf{F}_{k}"] &
\cdots \arrow[r, "\mathbf{F}_1"] &
\mathbf{V}_0 \arrow[r, "\mathbf{F}_0"] &
\mathcal{O}_{\widehat{Y}} \arrow[r, "\mathbf{F}_{-1}"] &
\mathcal{O}_{\widehat{B}} \arrow[r] &
0
\end{tikzcd}
\]
Additionally, we recall that the exact sequence for the resolution on \( \widehat{X} \) is determined by its leading terms:
\[
\begin{tikzcd}
0 \arrow[r] & \mathcal{O}_{\widehat{X}}^{m_k} \arrow[r, "L(F_k)"] & \mathcal{O}_{\widehat{X}}^{m_{k-1}} \arrow[r, "L(F_{k-1})"] & \cdots \arrow[r, "L(F_2)"] & \mathcal{O}_{\widehat{X}}^{m_1} \arrow[r, "L(F_1)"] & \mathcal{O}_{\widehat{X}}^{m_0} \arrow[r, "L(F_0)"] & \mathcal{O}_{\widehat{X}} \arrow[r] & j_{*}\mathcal{O}_{L(\widehat{B})} \arrow[r] & 0
\end{tikzcd}
\]
and the following exact sequence associated with the leading variety $L(\widehat{Y})$; defined by the leading terms of the equations of $\widehat{Y}$:
\[
\begin{tikzcd}
0 \arrow[r] & \mathcal{O}_{\widehat{X}} \arrow[r, "L(w)"] & \mathcal{O}_{\widehat{X}} \arrow[r] & \mathcal{O}_{L(\widehat{Y})} \arrow[r] & 0
\end{tikzcd}
\]
where we may assume that \( d := {\mathrm{ord}_{\widehat{P}}}(w) = {\mathrm{ord}_{\widehat{P}}}(L(w)) \). By the properties of the standard basis, we have \( L(w) \) contained in \( L(I) \). Therefore, in a completely analogous manner, we can repeat the previous algorithm for the above exact sequence and obtain the following leading term matrix factorization:
\[
\begin{tikzcd}[column sep=normal, row sep=normal]
& \cdots \arrow[r, "L(F_{3})"] & \mathcal{O}_{\widehat{X}}^{m_2}(-\widehat{Y}) \arrow[r, "L(F_{2})"] \arrow[d, "-L(K_3)"] & \mathcal{O}_{\widehat{X}}^{m_1}(-\widehat{Y}) \arrow[r, "L(F_{1})"] \arrow[d, "L(K_2)"] & \mathcal{O}_{\widehat{X}}^{m_0}(-\widehat{Y}) \arrow[r, "L(F_{0})"] \arrow[d, "-L(K_1)"] & \mathcal{O}_{\widehat{X}}(-\widehat{Y}) \arrow[d, "L(K_0)"] \\
& & \cdots \arrow[r, "L(F_{3})"] & \mathcal{O}_{\widehat{X}}^{m_2} \arrow[r, "L(F_{2})"] & \mathcal{O}_{\widehat{X}}^{m_1} \arrow[r, "L(F_{1})"] & \mathcal{O}_{\widehat{X}}^{m_0} \arrow[r, "L(F_{0})"] & \mathcal{O}_{\widehat{X}}
\end{tikzcd}
\]
Consistently, the above construction provides a unique \( L(w) \) matrix factorization in \( \widehat{X} \) for the stable periodic part of \( \mathcal{O}_{L(\widehat{B}),r} \):
\[
\begin{tikzcd}[column sep=normal, row sep=normal]
& & \vdots \arrow[d] & & \vdots \arrow[d] \\
\cdots \arrow[r, "L(F_{3})"] & \mathcal{O}_{\widehat{X}}^{m_2}(-2\widehat{Y}) \arrow[r, "L(F_{2})"] \arrow[d, "L(K_3)"] & \mathcal{O}_{\widehat{X}}^{m_1}(-2\widehat{Y}) \arrow[r, "L(F_{1})"] \arrow[d, "-L(K_2)"] & \mathcal{O}_{\widehat{X}}^{m_0}(-2\widehat{Y}) \arrow[r, "L(F_{0})"] \arrow[d, "L(K_1)"]\arrow[ddl,dashed, "L(K^{(1)}_1)",pos=0.6]  & \mathcal{O}_{\widehat{X}}(-2\widehat{Y}) \arrow[d, "-L(K_0)"] \arrow[ddl,dashed, "-L(K^{(1)}_0)",pos=0.7] \\
& \cdots \arrow[r, "L(F_{3})"] & \mathcal{O}_{\widehat{X}}^{m_2}(-\widehat{Y}) \arrow[r, "L(F_{2})"] \arrow[d, "-L(K_3)"] & \mathcal{O}_{\widehat{X}}^{m_1}(-\widehat{Y}) \arrow[r, "L(F_{1})"] \arrow[d, "L(K_2)"] & \mathcal{O}_{\widehat{X}}^{m_0}(-\widehat{Y}) \arrow[r, "L(F_{0})"] \arrow[d, "-L(K_1)"] & \mathcal{O}_{\widehat{X}}(-\widehat{Y}) \arrow[d, "L(K_0)"] \\
& & \cdots \arrow[r, "L(F_{3})"] & \mathcal{O}_{\widehat{X}}^{m_2} \arrow[r, "L(F_{2})"] & \mathcal{O}_{\widehat{X}}^{m_1} \arrow[r, "L(F_{1})"] & \mathcal{O}_{\widehat{X}}^{m_0} \arrow[r, "L(F_{0})"] & \mathcal{O}_{\widehat{X}}
\end{tikzcd}
\]
Moreover, the entries of all the mapping matrices are homogeneous, and under the natural projection to \( \widehat{Y} \), it gives rise to an acyclic double complex end with matrix factorization $(L(A),L(B))$
\[
\begin{tikzcd}[column sep=normal]
\cdots \arrow[r] &
\mathbf{V}_{k+2} \arrow[r, "L(\mathbf{F}_{k+2})"] &
\mathbf{V}_{k+1} \arrow[r, "L(\mathbf{F}_{k+1})"] &
\mathbf{V}_k \arrow[r, "L(\mathbf{F}_{k})"] &
\cdots \arrow[r, "L(\mathbf{F}_1)"] &
\mathbf{V}_0 \arrow[r, "L(\mathbf{F}_0)"] &
\mathcal{O}_{\widehat{Y}} \arrow[r, "L(\mathbf{F}_{-1})"] &
\mathcal{O}_{L(\widehat{B})} \arrow[r] &
0
\end{tikzcd}
\]
Next, following the usual approach, we consider the following commutative diagram. Our goal is to prove that it provides a well-defined acyclic double complex on \( \widehat{Y} \) (matrix factorization on \( \widehat{X} \)):
\[
\begin{tikzcd}[column sep=small, row sep=small]
 &\vdots \arrow[d] & \vdots \arrow[d] \\
 \bigoplus^{m_{1}}_{l=1} \widehat{P}^{r-\mathrm{a}^{(1)}_l-2d}(-2\widehat{Y}) \arrow[r, "F_{1}"] \arrow[d, "-K_2"] & \bigoplus^{m_{0}}_{l=1} \widehat{P}^{r-\mathrm{a}^{(0)}_l-2d}(-2\widehat{Y}) \arrow[r, "F_{0}"] \arrow[d, "K_1"] &  \widehat{P}^{r-2d}(-2\widehat{Y}) \arrow[d, "-K_0"] \arrow[ddl, "-K^{(1)}_0",dashed,pos=0.6]\\
 \cdots\arrow[r, "F_{2}"] & \bigoplus^{m_{1}}_{l=1} \widehat{P}^{r-\mathrm{a}^{(1)}_l-d}(-\widehat{Y}) \arrow[r, "F_{1}"] \arrow[d, "K_2"] & \bigoplus^{m_{0}}_{l=1} \widehat{P}^{r-\mathrm{a}^{(0)}_l-d}(-\widehat{Y}) \arrow[r, "F_{0}"] \arrow[d, "-K_1"] &  \widehat{P}^{r-d}(-\widehat{Y}) \arrow[d, "K_0"] \\
 \cdots & \bigoplus^{m_{2}}_{l=1} \widehat{P}^{r-\mathrm{a}^{(2)}_l} \arrow[r, "F_{2}"] & \bigoplus^{m_{1}}_{l=1} \widehat{P}^{r-\mathrm{a}^{(1)}_l} \arrow[r, "F_{1}"] & \bigoplus^{m_{0}}_{l=1} \widehat{P}^{r-\mathrm{a}^{(0)}_l} \arrow[r, "F_{0}"] & \widehat{P}^{r}
\end{tikzcd}
\]
we denote the complex also as follows:
\[
\begin{tikzcd}[column sep=normal]
\cdots \arrow[r] &
\mathbf{V}^r_{4} \arrow[r, "\mathbf{F}_{4}"] &
\mathbf{V}^r_{3} \arrow[r, "\mathbf{F}_{3}"] &
\mathbf{V}^r_2 \arrow[r, "\mathbf{F}_{2}"] &
\mathbf{V}^r_1  \arrow[r, "\mathbf{F}_1"] &
\mathbf{V}^r_0 \arrow[r, "\mathbf{F}_0"] &
\widehat{P}^r\mathcal{O}_{\widehat{Y}} \arrow[r, "\mathbf{F}_{-1}"] &
\widehat{P}^r\mathcal{O}_{\widehat{B}(M)} \arrow[r] &
0
\end{tikzcd}
\]

\begin{lemma}\label{lemmaleadinghper}
Consider any element \( t \) in \( \mathcal{O}^{m_i}_{\widehat{Y}}(-j\widehat{Y}) \) such that:
\[
(F_i(t), K_{i+1}(t)) \in \left( \bigoplus^{m_{i+1}}_{l=1} \widehat{P}^{r-\mathrm{a}^{(i)}_l - jd}(-j\widehat{Y}), \bigoplus^{m_{i-1}}_{l=1} \widehat{P}^{r-\mathrm{a}^{(i-1)}_l - (j-1)d}(-(j-1)\widehat{Y}) \right).
\]
Then, there exists an element \( t^\infty \) in \( \bigoplus^{m_{i}}_{l=1} \widehat{P}^{r-\mathrm{a}^{(i)}_l - jd}(-j\widehat{Y}) \) such that:
\[
(F_i(t), K_{i+1}(t)) = (F_i(t^\infty), K_{i+1}(t^\infty)).
\]
where \( i \) and \( j \) are arbitrary integers within the definition.
\end{lemma}

\begin{proof}
If \( t \) is an element satisfying the requirement, then we consider
\[
t \mod \bigoplus^{m_{i}}_{l=1} \widehat{P}^{r-\mathrm{a}^{(i)}_l - jd - k}(-j\widehat{Y})
\]
for any \( k \geq 1 \), and we assume \( k \) is the largest value for which this modulo is nonzero. Denote this element by \( [t]_k \). From the modulo equality, it is straightforward to deduce that:
\[
(L(F_i), L(K_{i+1}))([t]_k) = 0.
\]
there exists a homogeneous element \( [s]_{k} \) in
\[
\left( \bigoplus^{m_{i+1}}_{l=1} \widehat{P}^{r-\mathrm{a}^{(i)}_l - jd}(-(j+1)\widehat{Y}), \bigoplus^{m_{i-1}}_{l=1} \widehat{P}^{r-\mathrm{a}^{(i-1)}_l - (j-1)d}(-j\widehat{Y}) \right)
\]
such that:
\[
[t]_k = (L(K_{i}),L(F_{i+1}))([s]_{k}).
\]
Now, consider the new element:
\[
t' := t - (K_{i},F_{i+1})([s]_{k}).
\]
We know that:
\[
(F_{i},K_{i+1})(t) = (F_{i},K_{i+1})(t - (K_{i},F_{i+1}))([s]_{k}) = (F_{i},K_{i+1})(t').
\]
Moreover, if \( (K_{i},F_{i+1}) = L((K_{i},F_{i+1})) + R((K_{i},F_{i+1})) \) and \( t = [t]_k + \langle t \rangle_{k} \), then we have:
\[
t' = \langle t \rangle_{k} - R((K_{i},F_{i+1}))([s]_{k}).
\]
By construction, the order of \( t' \) is bigger than the order of \( t \) under the natural module ordering. We now replace \( t \) with \( t' \) and repeat the above algorithm. Since there are only finitely many corresponding orders, this algorithm will terminate after a finite number of steps. We output an element \( t^{\infty} \) satisfying:
\[
(F_{i},K_{i+1})(t) = (F_{i},K_{i+1})(t^{\infty}),
\]
and \( t^{\infty} \) belongs to \( \bigoplus^{m_{i}}_{l=1} \widehat{P}^{r-\mathrm{a}^{(i)}_l - jd}(-j\widehat{Y}) \).
\end{proof}

\begin{corollary}
For any integers \( r \) and \( s \), the above double complex is acyclic, and it is quasi-isomorphic to \( \widehat{P}^r \cdot \mathcal{O}_{\widehat{B}, s} \) on \( \widehat{Y} \).
\end{corollary}

\begin{proof}
Since the proof is identical to the previous one, we provide a brief description here. First, by the definition of the leading order for each term, the map is well-defined between different ideals, e.g.,
\[
\mathrm{Im}\, \mathbf{F}_i |_{\mathbf{V}^r_i} \subset \mathbf{V}^r_{i+1}.
\]
Moreover, since the original double complex is acyclic, we have
\[
\mathrm{Im}\, \mathbf{F}_i |_{\mathbf{V}^r_i} \subset \mathrm{Ker}\, \mathbf{F}_{i+1} \cap \mathbf{V}^r_{i+1}.
\]
Additionally, by a similar argument as Lemma \ref{lemmaleadinghper}, for any integer $i$ and element $t$, if there is
\[
\mathbf{F}_i(t) \in \mathbf{V}^r_{i+1},
\]
then there exists an element $t^\infty \in \mathbf{V}^r_i$ satisfying \(\mathbf{F}_i(t) = \mathbf{F}_i(t^\infty)\),
so we also have
\[
\mathrm{Ker}\, \mathbf{F}_{i+1} \cap \mathbf{V}^r_{i+1} = \mathrm{Im}\, \mathbf{F}_i \cap \mathbf{V}^r_{i+1} \subset \mathrm{Im}\, \mathbf{F}_i |_{\mathbf{V}^r_i},
\]
where $\mathbf{V}^r_{i}$ is the term and $\mathbf{F}_i$ is the map of degree $-i$ in the totalization of the above complex.
\end{proof}
So for sufficiently large \( s \), we get the following standard exact resolution of \(\widehat{P}^{r}\mathcal{O}_{\widehat{B}, s}\) in \( \widehat{S} \):
\[
\begin{tikzcd}[column sep=small]
 \cdots & \bigoplus^{m_{2}}_{l=1} \widehat{P}^{r-\mathrm{a}^{(2)}_l} \oplus \bigoplus^{m_{0}}_{l=1} \widehat{P}^{r-\mathrm{a}^{(0)}_l-d}(-\widehat{Y}) \arrow[r, "\mathbf{F}_{2}"] & \bigoplus^{m_{1}}_{l=1} \widehat{P}^{r-\mathrm{a}^{(1)}_l} \oplus \widehat{P}^{r-d}(-\widehat{Y}) \arrow[r, "\mathbf{F}_{1}"] & \bigoplus^{m_{0}}_{l=1} \widehat{P}^{r-\mathrm{a}^{(0)}_l}
\end{tikzcd}
\]
\[
\begin{tikzcd}[column sep=small]
\bigoplus^{m_{0}}_{l=1} \widehat{P}^{r-\mathrm{a}^{(0)}_l} \arrow[r, "\mathbf{F}_{0}"] & \widehat{P}^{r} \arrow[r, "\mathbf{F}_{-1}"] & \widehat{P}^r\mathcal{O}_{\widehat{B}} \arrow[r] & 0
\end{tikzcd}
\]
In particular, we have the following description for the term of degree \(-(k-1)\):
\[
\mathbf{V}^r_{k-1}:=\bigoplus_{j \mid -1 \leq k-2j-2 \leq k-2}\bigoplus^{m_{k-2j-2}}_{l=1} \widehat{P}^{r-\mathrm{a}^{(k-2j-2)}_l - jd}(-j\widehat{Y}).
\]
and term of degree \(-k\)
\[
\mathbf{V}^r_{k}:=\bigoplus_{j \mid -1 \leq k-2j-1 \leq k-1}\bigoplus^{m_{k-2j-1}}_{l=1} \widehat{P}^{r-\mathrm{a}^{(k-2j-1)}_l - jd}(-j\widehat{Y}).
\]
and also for \( 1 \leq 2t \leq s-1 \),
\[
\mathbf{V}^r_{k+2t-1} = \widehat{P}^{-td} \cdot \mathbf{V}^r_{k-1}(-t\widehat{Y}), \quad \text{and} \quad \mathbf{V}^r_{k+2t} = \widehat{P}^{-td} \cdot \mathbf{V}^r_{k}(-t\widehat{Y}),
\]
when \( t \) also is sufficiently large, all the exponents of \( \widehat{P} \) become negative. As a result, the sequence gradually degenerates into the previously given stable period-two complex of \( \mathbf{V}_{k+2t-1} \) and \( \mathbf{V}_{k+2t} \).\\

\begin{remark}
We note the flexibility in our construction. Suppose \(\widehat{S}\) has two presentations:
\[
\widehat{S} \simeq \mathrm{k}[[x_1,\dots,x_n]]/(f) \simeq \mathrm{k}[[x'_1,\dots,x'_n]]/(f'), \quad
\widehat{R} \simeq \mathrm{k}[[x_1,\dots,x_n]] \simeq \mathrm{k}[[x'_1,\dots,x'_n]],
\]
with \(\alpha: \mathrm{k}[[x'_i]] \xrightarrow{\cong} \mathrm{k}[[x_i]]\) such that \(\alpha(f') = u \cdot f\) for a unit \(u\) \cite[Def.~9.1.1]{DeJongPfister}. For all \(r\), the maps \(\phi,\varphi,h^\infty,k^\infty\) induce a homotopy commutative diagram relating the resolutions constructed above from \(f\) and \(\alpha(f')\) (cf. Lemma~\ref{gradedfiterexact}), compatible with the division construction. Corresponding direct sums induce homotopies between the two presentations on \(\widehat{Y}\).
\end{remark}

Even if a reduced standard resolution of \(j_*\mathcal{O}_{\widehat{B}}\) in \(\widehat{R}\) exists, the resolution over \(\widehat{S}\) may not be reduced. Elementary operations on \(\mathbf{V}^r_i\) remove rows/columns outside \(\mathfrak{m}\) without changing leading degrees; by periodic stabilization it suffices to apply this to two consecutive matrices \((\mathbf{F}_k,\mathbf{F}_{k+1})\). Hence, for sufficiently large \(s\) and any \(r\), there exists a reduced resolution of \(\widehat{P}^r \mathcal{O}_{\widehat{B}}\) on \(\widehat{Y}\) depending only on a lexicographic order, satisfying the periodicity
\[
\mathbf{V}^r_{k+2t-1} = \widehat{P}^{-td} \mathbf{V}^r_{k-1}(-t\widehat{Y}), \quad
\mathbf{V}^r_{k+2t} = \widehat{P}^{-td} \mathbf{V}^r_k(-t\widehat{Y}), \quad
\mathbf{V}^{r+j}_i = \widehat{P}^j \mathbf{V}^r_i.
\]
So we call the resolutions constructed above (without reduction) up to homotopy the \textbf{standard resolutions} of \(\mathcal{O}_{\widehat{B}}\) on \(\widehat{Y}\), and the reduced version as a resolution with the \emph{asymptotically periodic} property.\\

For any Cohen–Macaulay module \( M \) on \( \widehat{Y} \) and any Bourbaki sequence associated with \( M \) on \( \widehat{Y} \):
\[
0 \longrightarrow L \longrightarrow M \longrightarrow I(M) \longrightarrow 0,
\]
we examine the corresponding Bourbaki cycle \( \widehat{B}(M) \) and obtain the following standard resolution for \( j_*\mathcal{O}_{\widehat{B}(M)} \) on \( \widehat{X} \):
\[
\begin{tikzcd}
\cdots \arrow[r] & \mathcal{O}_{\widehat{X}}^{m_2} \arrow[r, "F_2"] & \mathcal{O}_{\widehat{X}}^{m_1} \arrow[r, "F_1"] & \mathcal{O}_{\widehat{X}}^{m_0} \arrow[r, "F_0"] & \mathcal{O}_{\widehat{X}} \arrow[r] & j_{*}\mathcal{O}_{\widehat{B}(M)} \arrow[r] & 0
\end{tikzcd}
\]
For all integers \( r \), we have the following standard resolutions:
\[
\begin{tikzcd}
\cdots \arrow[r] & \bigoplus^{m_{2}}_{l=1} \widehat{P}^{r-\mathrm{a}^{(2)}_l} \arrow[r, "F_2"] & \bigoplus^{m_{1}}_{l=1} \widehat{P}^{r-\mathrm{a}^{(1)}_l} \arrow[r, "F_1"] & \bigoplus^{m_0}_{l=1} \widehat{P}^{r-\mathrm{a}^{(0)}_l} \arrow[r, "F_0"] & \widehat{P}^{r} \arrow[r, "F_{-1}"] & \widehat{P}^{r}\mathcal{O}_{\widehat{B}(M)} \arrow[r] & 0
\end{tikzcd}
\]
and then standard resolutions of \( \mathcal{O}_{\widehat{B}} \) on \( \widehat{Y} \) for all integers \( r \):
\[
\begin{tikzcd}[column sep=normal]
\cdots \arrow[r] &
\mathbf{V}^r_{4} \arrow[r, "\mathbf{F}_{4}"] &
\mathbf{V}^r_{3} \arrow[r, "\mathbf{F}_{3}"] &
\mathbf{V}^r_2 \arrow[r, "\mathbf{F}_{2}"] &
\mathbf{V}^r_1  \arrow[r, "\mathbf{F}_1"] &
\mathbf{V}^r_0 \arrow[r, "\mathbf{F}_0"] &
\widehat{P}^r \arrow[r, "\mathbf{F}_{-1}"] &
\widehat{P}^r \mathcal{O}_{\widehat{B}(M)} \arrow[r] &
0
\end{tikzcd}
\]
It is straightforward to see that, by using the periodicity property, we can extend the above exact sequences to obtain the standard version of the Bourbaki sequences for any integer \(r\):
\[
0 \longrightarrow \mathbf{L}^r \longrightarrow \mathbf{M}^r \longrightarrow \widehat{P}^r \longrightarrow \widehat{P}^r \mathcal{O}_{\widehat{B}(M)} \longrightarrow 0,
\]
\( \mathbf{M}^r \) has the following \textbf{standard resolution}:
\[
\begin{tikzcd}[column sep=normal]
\cdots \arrow[r] &
\mathbf{U}^r_4 \arrow[r, "\mathbf{\bar{F}}_4"] &
\mathbf{U}^r_3 \arrow[r, "\mathbf{\bar{F}}_3"] &
\mathbf{U}^r_2 \arrow[r, "\mathbf{\bar{F}}_2"] &
\mathbf{U}^r_1 \arrow[r, "\mathbf{\bar{F}}_1"] &
\mathbf{U}^r_0 \arrow[r, "\mathbf{\bar{F}}_0"] &
\mathbf{M}^r \arrow[r] &
0
\end{tikzcd}
\]
and \( \mathbf{L}^r \) also admits a standard resolution of the form:
\[
0 \longrightarrow \mathbf{W}^r_{k} \longrightarrow \cdots \longrightarrow \mathbf{W}^r_{1} \longrightarrow \mathbf{L}^r \to 0
\]
and where we set for any integer \( t \),
\[
\mathbf{U}^r_{k+2t-1} := \widehat{P}^{-td} \cdot \mathbf{V}^r_{k-1}(-t\widehat{Y}), \quad \text{and} \quad \mathbf{U}^r_{k+2t} := \widehat{P}^{-td} \cdot \mathbf{V}^r_{k}(-t\widehat{Y}).
\]
and the module \( \mathbf{W}_i^r \) as the quotient:
\[
\mathbf{W}_i^r := \mathbf{U}_i^r / \mathbf{V}_i^r,
\]
for any integer \( i \), and the quotient is taken via the natural embedding of the submodules. Accordingly, we naturally extend or restrict the map \( \mathbf{F}^r_i \) to the new resolutions.\\

Moreover, defining \( \widehat{A} := \mathbf{\bar{F}}_i \) for \( i \equiv 0 \pmod{2} \) and \( \widehat{B} := \mathbf{\bar{F}}_i \) for \( i \equiv 1 \pmod{2} \), we obtain an infinite extension of the exact sequence on \( Y \), given by
\[
\begin{tikzcd}[column sep=normal]
\cdots \arrow[r] &
\mathbf{U}^r_{2} \arrow[r, "\widehat{A}"] &
\mathbf{U}^r_{1} \arrow[r, "\widehat{B}"] &
\mathbf{U}^r_{0} \arrow[r, "\widehat{A}"] &
\mathbf{U}^r_{-1} \arrow[r, "\widehat{B}"] &
\mathbf{U}^r_{-2} \arrow[r, "\widehat{A}"] &
\cdots
\end{tikzcd}
\]
Thus, in this formulation, we have
\[
\mathbf{M}^r = \operatorname{Coker} \widehat{A} \big|_{\mathbf{U}^r_0},
\]
Based on the previous construction, for a matrix \( \widehat{A} \) (or \( \widehat{B} \)), we have its leading matrix \( L(\widehat{A}) \) (and \( L(\widehat{B}) \)). Similarly, we define their leading cokernels on \(L(\widehat{Y})\) as follows:
\[
L(\mathbf{M}^r): = \operatorname{Coker} L(\widehat{A}) \big|_{\mathbf{U}^r_0},
\]
noticing this formulation allows some freedom in the choice of the Bourbaki sequence.\\

\subsection{Measure resolution on hypersurface}\label{subsec:mea}
Let \((R,\mathfrak m)\) be a regular local ring, \(J=(w)\) a hypersurface in \(R\). By a valuation\footnote{Although we mainly consider nonnegative valuations on \(R\), the following result remains valid for valuations taking negative values as well.} on \(R\) we mean a map
\[
val:R\to\mathbb Z\cup\{\infty\},
\qquad
val(0)=\infty,
\]
satisfying \(val(fg)=val(f)+val(g)\) and \(val(f+g)\ge\min\{val(f),val(g)\}\). Since matrix factorizations only depend on elements on \(R\), we shall not consider valuations on \(S\). For \(a\in\mathbb Z\), define
\[
P_a:=\{\,f\in R\mid val(f)\ge a\,\}.
\]
We also consider its extension \(P_aS\subset S\). By abuse of notation, we shall denote \(P_aS\) by \(P_a\) whenever no confusion can arise. Later, we specialize step by step to weighted monomial valuations. We consider the following property and refer to \cite[Chap.~7]{Yoshino1990} for background on matrix factorizations.
\begin{definition}

\begin{enumerate}

    \item A rank \( n \) matrix factorization \( (A,B) \) of an Cohen–Macaulay module \( M \) on \( S \) is called \emph{asymptotically periodic} presentation under valuation \(val\) if there exists a pair of submodules \( (U_{1},U_{2}) \) of free modules such that:
    \begin{enumerate}
        \item \( U_{1} = \bigoplus_{i=1}^{n} P_{a_{i}} \) and \( U_{2} = \bigoplus_{j=1}^{n} P_{b_{j}} \), where \( a_{i}, b_{j} > 0 \) are positive integers.
        \item The shifted modules are defined by \( U_{1}(d) := \bigoplus_{i=1}^{n} P_{a_{i}+d} \) and \( U_{2}(d) := \bigoplus_{j=1}^{n} P_{b_{j}+d} \).
        \item The maps \( A: U_1 \to U_2 \) and \( B: U_2 \to U_1(d) \) are well-defined on \(R\), where \( d := val(w) \).
    \end{enumerate}

    \item For a morphism \( A \colon (U_{1},U_{2}) \to (U_{2},U_{1}(d)) \), the \emph{cone of the image of \( A \) restricted to \( U_{1} \)} is defined as:
    \[
    \mathrm{Cone}(A|_{U_{1}}) := \bigcap_{\substack{\{c_{i}\} \\ \mathrm{Im}\, A|_{U_{1}} \,\subset\, \bigoplus_{i} P_{c_{i}}}} \bigoplus_{i} P_{c_{i}},
    \]
    which is the minimal module of the form \( \bigoplus_{i} P_{c_{i}} \) containing \( \mathrm{Im}\, A|_{U_{1}} \).
\end{enumerate}
\end{definition}

\begin{lemma}
  If \( (A,B) \) is an asymptotic periodic presentation of a Cohen–Macaulay module \( M \) on \( S\), then we have
  \[
  U_{2} = \mathrm{Cone}(A|_ {U_{1}}) \quad \text{and} \quad U_{1}(d) = \mathrm{Cone}(B|_ {U_{2}}).
  \]
\end{lemma}

\begin{proof}
Let \(\operatorname{Cone}(X)\) denote the minimal module of the form \(\bigoplus_i P_{c_i}\) containing \(X\subset R^n\), relative to the chosen basis of the matrix factorization. In particular, for \(X\subset Y\subset R^n\), one has
\[
\operatorname{Cone}(X)\subset \operatorname{Cone}(Y).
\]
By our definition, we have \( \mathrm{Cone}(A|_ {U_{1}}) \subset U_{2} \) and
\[
\mathrm{Cone}(B|_ {\mathrm{Cone}(A|_ {U_{1}})}) \subset \mathrm{Cone}(B|_ {U_{2}}) \subset U_{1}(d).
\]
Since \(val(w)=d\) and we have
$w\cdot U_1\subset B\bigl(A(U_1)\bigr)\subset B\bigl(\operatorname{Cone}(A(U_1))\bigr)$. Therefore,
\[
U_1(d)
\subset
\operatorname{Cone}(w\cdot U_1)
\subset
\mathrm{Cone}(B|_ {(A|_ {U_{1}})})
\subset
\mathrm{Cone}(B|_ {\mathrm{Cone}(A|_ {U_{1}})}).
\]
This shows that all the inclusions above are in fact isomorphisms.
\end{proof}
The following lemma is trivial, but it illustrates that the above definition constitutes a purely numerical structure.
\begin{lemma}
If \( (A,B) \) is an asymptotic periodic presentation of a Cohen–Macaulay module \( M \) on \( S \) such that our definition holds under submodules \( (U_{1},U_{2}) := (\bigoplus_{i=1}^{n} P_{a_{i}},\bigoplus_{j=1}^{n} P_{b_{j}}) \), then our definition also holds under submodules \( (U_{1}(k),U_{2}(k)) := (\bigoplus_{i=1}^{n} P_{a_{i}+k},\bigoplus_{j=1}^{n} P_{b_{j}+k}) \) for any integer \( k \).
\end{lemma}

\begin{proof}
This is entirely determined by the valuation of the matrix entries. If the matrix pair \( (A,B) \) has entries
\[
(A_{ij})_{i,j}, (B_{mn})_{m,n}
\]
with valuation
\[
val((A_{ij})_{i,j},(B_{mn})_{m,n}) = ((a_{ij})_{i,j},(b_{mn})_{m,n}),
\]
then our asymptotic periodicity condition for \( (U_{1},U_{2}) \) is given by the inequalities:
\[
\min_{j=1,\dots,n} \{a_{ij} + a_{j} \} \geq b_{i}
\quad
\min_{j=1,\dots,n} \{b_{ij} + b_{j} \} \geq a_{i} + d,
\]
for all \( i \) from \( 1 \) to \( n \), where \( d := val(w) \). It is clear that if the integer sequences \( (a_{i})_{i},(b_{j})_{j} \) satisfy these conditions, then \( (a_{i}+k)_{i},(b_{j}+k)_{j} \) also satisfy the same conditions for any integer \( k \). This corresponds to submodules \( (U_{1}(k),U_{2}(k)) := (\bigoplus_{i=1}^{n} P_{a_{i}+k},\bigoplus_{j=1}^{n} P_{b_{j}+k}) \).
\end{proof}

\begin{lemma}
Given an asymptotically periodic representation $(A,B)$ on $R$, there exists an induced family of long exact sequences on $S$ parameterized by $k\in\mathbb{Z}$:
\[
\cdots \to U_{1}(k) \xrightarrow{A} U_{2}(k) \xrightarrow{B} U_{1}(k+d) \xrightarrow{A} U_{2}(k+d) \xrightarrow{B} U_{1}(k+2d) \to \cdots
\]
\end{lemma}

\begin{proof}
Assuming $k=0$, we need to prove $\mathrm{Im}\, A|_{U_1} = \mathrm{Ker}\, B \cap U_2$.\\

First, the inclusion $\mathrm{Im}\, A|_{U_1} \subset \mathrm{Im}\, A \cap U_2 = \mathrm{Ker}\, B \cap U_2$ follows immediately from the definitions. For the reverse inclusion, consider any $v \in U_2$ such that $B(v) = 0$ in $S$. Since $(A,B)$ forms a matrix factorization of $w$, in $R^n$ we have element $u$ and the identity $w \cdot u = B\circ A(u) = B(v)$, noticing $B(v)$ is in $U_1(d)$ by assumption and $w$ is an order $d$ divisor, this implies $u$ is in $U_1$. Therefore, $v$ belongs to $\mathrm{Im}\, A|_{U_1}$, completing the proof.
\end{proof}

There clearly exist examples where certain presentations fail to be asymptotically periodic. Consider $R = \mathrm{k}[x,y]_{\langle x,y \rangle}$ with $w = x^{a} + y^{b}$ for integers $b \geq a > 2$, and $P := \langle x,y \rangle$ with $val:=\mathrm{ord}_{P}$. We have two equivalent matrix factorizations:

\[
\begin{pmatrix}
x^{a-1} & y^{b-1} \\
y & -x \\
\end{pmatrix},
\begin{pmatrix}
x & y^{b-1} \\
y & -x^{a-1} \\
\end{pmatrix}
\]
and
\[
\begin{pmatrix}
x^{a-1}+y & y^{b-1}-x \\
y & -x \\
\end{pmatrix},
\begin{pmatrix}
x & y^{b-1}-x \\
y & -x^{a-1}-y \\
\end{pmatrix}.
\]
The first representation is asymptotically periodic, while the second is not.
\begin{prop}
Let \( M \) be a Cohen–Macaulay module over  the analytic/formal local hypersurface ring. Suppose that \( val \) denotes the order with respect to a regular center \( P \). Then, any reduced standard resolution constructed in section \ref{subsec:hyp} provides an \emph{asymptotically periodic} representation of \( M \).
\end{prop}

\begin{lemma}\label{lemma1asy}
Let $\mathrm{k}$ be an algebraically closed field of characteristic zero.
If $(A,B)$ is an asymptotic periodic representation of a Cohen–Macaulay module on $R/w$ with respect to a valuation $val$, then the matrix factorization
\[
(A',B') :=
\left(
\begin{pmatrix}
A & u \\
v & -B
\end{pmatrix},
\begin{pmatrix}
B & u \\
v & -A
\end{pmatrix}
\right)
\]
representing the corresponding Cohen–Macaulay module on $R[u,v]_{\mathfrak{m}'}/(w+uv)$ is also asymptotically periodic with respect to any valuation $val'$ of $R[u,v]_{\mathfrak{m}'}$ such that its restriction to $R$ is $val$ and
\[
val'(w+uv) = \min\{val'(w), val'(u) + val'(v)\}.
\]
\end{lemma}

\begin{proof}
If $(A,B)$ is asymptotically periodic with submodules $(U_1, U_2)$, then there is an exact sequence
\[
\cdots \xrightarrow{A} U_2 \xrightarrow{B} U_1(d) \xrightarrow{A} U_2(d) \xrightarrow{B} \cdots
\]
over $R/w$, where $d = val(w)$.

Let $val'(u) = a$, $val'(v) = b$, and $val'(w+uv) = d' = \min\{d, a+b\}$. Consider the sequence
\[
U_1 \oplus U_2(k) \xrightarrow{A'} U_2(\min\{0, a+k\}) \oplus U_1(\min\{b, d+k\}) \xrightarrow{B'} U_1(\min\{d, a+b, a+d+k\}) \oplus U_2(\min\{b, d+k, a+b+k\}).
\]
This sequence is well-defined over $R[u,v]_{\mathfrak{m}'}/(w+uv)$ for any integer $k$. Choosing $k \le b - d'$, we obtain
\[
U_1 \oplus U_2(k) \xrightarrow{A'} U_2(\min\{0, a+k\}) \oplus U_1(\min\{b, d+k\}) \xrightarrow{B'} U_1(d') \oplus U_2(k+d'),
\]
which is asymptotically periodic with respect to $val'(w+uv)$, and hence exact over $R[u,v]_{\mathfrak{m}'}/(w+uv)$.
\end{proof}

\begin{lemma}\label{lemma2asy}
If $(A,-A)$ is an asymptotic periodic representation of certain Cohen–Macaulay module on $R/w$ with respect to $val$, then the matrix factorization
$(u+A,u-A)$
representing the corresponding Cohen–Macaulay module on $R[u]_{\mathfrak{m}'}/w+u^2$ is also asymptotic period with respect to any $val'$ of $R[u]_{\mathfrak{m}'}$ such that its restriction on $R$ is $val$ and $val'(w+u^2)=\min\{val'(w), val'(u^2)\}$.
\end{lemma}

\begin{proof}
If \( (A, -A) \) is asymptotically periodic with submodules \( (U_1, U_2) \), so we have an exact sequence
\[
\cdots \xrightarrow{A} U_2 \xrightarrow{-A} U_1(d) \xrightarrow{A} U_2(d) \xrightarrow{-A} \cdots
\]
on \( R/w \). Assume \( val'(u) = a \) and \( val'(w + u^2) = d' =\min\{d, 2a\} \). Then the sequence
\[
U_1 \cup U_2(k) \xrightarrow{u + A} U_2(\min\{0, a + k\}) \cup U_1(\min\{a, d + k\}) \xrightarrow{u - A} U_1(\min\{d, 2a, a + d + k\}) \cup U_2(\min\{a, d + k, 2a + k\})
\]
is well-defined over \( R[u]_{\mathfrak{m}'} / (w + u^2) \) for any integer \( k \). So if we choose \( d' - d - a \leq k \leq a - d' \), which is a nonempty interval, we obtain
\[
U_1 \cup U_2(k) \xrightarrow{u + A} U_2(\min\{0, a + k\}) \cup U_1(\min\{a, d + k\}) \xrightarrow{u - A} U_1(d') \cup U_2(k + d')
\]
which is asymptotically periodic  of \(val'(w+u^2)\) on $R[u]_{\mathfrak{m}'}/w+u^2$ hence exact.
\end{proof}

\begin{lemma}\label{lemma3asy}
For one dimensional or two dimensional simple singularity e.g. $(R,w)\simeq(\mathrm{k}[x,y]_\mathfrak{m},w)$ or $(R,w)=(\mathrm{k}[x,y,z]_\mathfrak{m},w)$  where $w$ is any simple polynomial and $\mathrm{k}$ be an algebraically closed field of characteristic zero. Then their Cohen–Macaulay modules admit, for all its representations given in the \cite[Chap.9]{Yoshino1990} and \cite[Chap.9, Sec.4]{LeuschkeWiegand2012}, a choice of submodules such that it becomes asymptotically periodic under any valuation  such that \( val(f) \) is equal to the minimal valuation of all monomials appearing in \( f \).
\end{lemma}

\begin{proof}[Outline of the proof]
If the matrix pair \( (A, B) \) has entries
\[
(A_{ij})_{i,j}, \quad (B_{mn})_{m,n}
\]
with valuations
\[
val\big((A_{ij})_{i,j}, (B_{mn})_{m,n}\big) = \big((a_{ij})_{i,j}, (b_{mn})_{m,n}\big),
\]
then the existence of  asymptotic periodicity presentation is equivalent to the existence of integers \( a_i \) and \( b_i \) such that the following inequalities hold:
\[
\min_{j = 1, \dots, n} \{ a_{ij} + a_j \} \geq b_i,
\quad
\min_{j = 1, \dots, n} \{ b_{ij} + b_j \} \geq a_i + d,
\]
for all \( i = 1, \dots, n \), where \( d := val(w) \). We verify the matrix factorizations \( (A, B) \) for all simple singularities appearing in the literature. In each case, we obtain a non-empty solution set to the system of inequalities described above.
\end{proof}
Combining all the Lemma \ref{lemma1asy}, \ref{lemma2asy}, \ref{lemma3asy} and \cite[Chap.12]{Yoshino1990}, we arrive at the following conclusion.
\begin{corollary}
Let $(R, w) \simeq (\mathrm{k}[x_0, \ldots, x_n]_{\mathfrak{m}}, w)$ be a local ring defined by a simple singularity, where $w$ is a simple polynomial and $\mathrm{k}$ be an algebraically closed field of characteristic zero. Then, under any choice of valuation  such that \( val(f) \) is equal to the minimal valuation of all monomials appearing in \( f \), every Cohen–Macaulay module over $(R, w)$ admits an asymptotically periodic presentation.
\end{corollary}

\begin{lemma}
Let $S = R/(w)$ with a valuation $val$. Given any short exact sequence of Cohen–Macaulay modules:
\[
0 \to M \to M'' \to M' \to 0
\]
where $M$ (resp. $M'$) admits a matrix factorization $(A,B)$ (resp. $(A',B')$), then by \cite[Rem.7.8]{Yoshino1990} $M''$ has an induced matrix factorization $(A'',B'')$ with:
\[
A'' := \begin{pmatrix}
A & E \\
0 & A'
\end{pmatrix}, \quad
B'' := \begin{pmatrix}
B & F \\
0 & B'
\end{pmatrix}
\]
Let $(U_1,U_2)$ and $(U'_1,U'_2)$ be asymptotic periodic presentations of $(A,B)$ and $(A',B')$ respectively. Then there exists an integer $k_0 \gg 0$ such that for all $k \geq k_0$, the direct sum
\[
(U_1 \oplus U'_1(k), U_2 \oplus U'_2(k))
\]
forms an asymptotic periodic presentation of $(A'',B'')$.

\end{lemma}
\begin{proof}
Suppose we are given matrices with respective entries
\[
A = (A_{ij}), \quad B = (B_{mn}), \quad A' = (A'_{ij}), \quad B' = (B'_{mn}), \quad E = (E_{ij}), \quad F = (F_{mn}),
\]
and their valuations are
\[
val(A_{ij}) = a_{ij}, \quad val(B_{mn}) = b_{mn}, \quad
val(A'_{ij}) = a'_{ij}, \quad val(B'_{mn}) = b'_{mn}, \quad
val(E_{ij}) = e_{ij}, \quad val(F_{mn}) = f_{mn}.
\]

Let the associated graded modules be
\[
(U_1, U_2) := \left( \bigoplus_{i=1}^{n} P_{a_i}, \bigoplus_{j=1}^{n} P_{b_j} \right), \quad
(U_1', U_2') := \left( \bigoplus_{i=1}^{n} P_{a'_i}, \bigoplus_{j=1}^{n} P_{b'_j} \right),
\]
then by assumption the following inequalities hold:
\[
\min_j \{ a_{ij} + a_j \} \geq b_i, \qquad
\min_j \{ b_{ij} + b_j \} \geq a_i + d,
\]
\[
\min_j \{ a'_{ij} + a'_j \} \geq b'_i, \qquad
\min_j \{ b'_{ij} + b'_j \} \geq a'_i + d.
\]
So there exists \( k \gg 0 \) such that
\[
\min_j \{ a_{ij} + a_j,  e_{ij} + a'_j+k \} \geq b_i, \qquad
\min_j \{ b_{ij} + b_j,  f_{ij} + b'_j+k \} \geq a_i + d,
\]
\[
\min_j \{ a'_{ij} + a'_j +k\} \geq b'_i+k, \qquad
\min_j \{ b'_{ij} + b'_j +k\} \geq a'_i + d+k.
\]
\end{proof}

\section{Transformation}

\subsection{Blow-up}\label{section:stricttrans}
 Let \( (R, m) \) be a regular \textbf{analytic/formal} \footnote{For the case when \(R\) is \textbf{algebraic}, we refer to the discussion preceding Theorem~\ref{thm:algebraic-version}. First, we need to handle the formal case.}\label{expla}
 local ring with residue field \(\mathrm{k}\), and let \( J = \langle w \rangle \) be a hypersurface. Suppose \( P \) is a prime ideal in \( R \) and we require a classic constraint on \( P \) and \( J \):

\begin{definition}\label{permissiblecenter}
A prime ideal \( P \) in \( R \) is called a \textbf{permissible center} for \( J \), if it satisfies the following conditions:
\begin{enumerate}
    \item[(i)] \( P \supseteq J \) and \( R/P \) is regular,
    \item[(ii)] \( \mathrm{gr}_{P/J}(R/J):= \bigoplus_{i \geq 0} (P/J)^i/(P/J)^{i+1} \) is flat over \( R/P \).
\end{enumerate}
\end{definition}
Set
\(X = \operatorname{Spec} R\),  \(Y = \operatorname{Spec} R/J\) and \(Z = \operatorname{Spec}(R/P)\). The second condition states that the order of \( w \) along the center \( Z \) remains invariant. In the case where \( S \) is a hypersurface, this condition is equivalent to requiring that \( Z \) is contained in \( \mathrm{Sing}(S) \) (or \( \mathrm{Reg}(S) \)). We consider the global monoidal transformation of \(R\) along \(P\),
\[
\mathrm{Bl}_P R := \mathcal{R}(P)
  := \bigoplus_{i \ge 0} P^i ,
\qquad
X^+ := \operatorname{Proj} \mathrm{Bl}_P R.
\]
The morphism \( X^+ \to X \) is projective, and its exceptional divisor
is a projective bundle over \( Z \).
Additionally, we consider the strict transform of \( Y \) (or any cycle \( B \)).
This is described by
\[
\mathrm{Bl}_P S
  := \bigoplus_{i \ge 0} P^i / J\cap P^i ,
\qquad
Y^+ := \operatorname{Proj} \mathrm{Bl}_P S ,
\]
we denote the corresponding strict transform by \( Y^+ \) or \( B^+ \). These constructions fit into the following commutative diagram:
\[
\begin{tikzcd}
E_{Y} \arrow[r, hookrightarrow, "i^+"]\arrow[d, "\pi_{E_{Y}}"] & Y^+ \arrow[r, hookrightarrow, "j^+"]\arrow[d, "\pi_{Y}"] & X^+ \arrow[r, hookleftarrow, "i^+"]\arrow[d, "\pi"] & E \arrow[d, "\pi_{E}"]\\
Z \arrow[r, hookrightarrow, "i"] & Y \arrow[r, hookrightarrow, "j"] & X \arrow[r, hookleftarrow,"k"]& Z
\end{tikzcd}
\]
In particular, since \( Z \) is permissible, \( E_{Y} \) is a flat hyper-fibration over the center \( Z \).\\

According to Serre's construction, we have the following abelian category isomorphism:
\[
\Pi: \mathrm{Gr}\big(\mathcal{R}(P)\big) \big/ \mathrm{Tors}\big(\mathcal{R}(P)\big) \simeq \mathrm{Coh}(X^+) : \Gamma,
\]
by \cite[Thm. 3.4.4]{EGAII}, where
\[
\Gamma(M) := \bigoplus_{i \geq 0} \mathrm{H}^{0}(X^+, M(i)).
\]
As before, we set $P^i=R$ for all $i\le 0$. Since like $k\ge 1$, we have the graded description
\[
\mathcal{O}_{X^+}(-k)\longleftrightarrow
\underbrace{R\oplus\cdots\oplus R}_{\text{$k$ copies created by Serre torsion}}
\oplus R\oplus P\oplus P^2\oplus\cdots .
\]
so for any integer \( k \), we have:
\[
\mathcal{O}_{X^+}(-k) := \Pi\big(\mathcal{R}(P)(-k)\big) = \Pi\big(\bigoplus_{i \geq 0} P^{i-k}\big),
\]
from this and the direct sum of the following exact sequence for any integer \( r \):
\[
0 \longrightarrow P^{r+1} \longrightarrow P^{r} \longrightarrow P^{r}/P^{r+1} \longrightarrow 0.
\]
it follows that:
\[
\mathcal{O}_{X^+}(-k) \simeq \mathcal{O}_{X^+}(kE).
\]
Therefore, if we have standard resolution of \(P^r\mathcal{O}_{B}\) on \(R\) as in Section \ref{subsec:reg}:
\[
\begin{tikzcd}
0 \arrow[r] & \bigoplus^{m_{k}}_{l=1} P^{r-\mathrm{a}^{(k)}_l} \arrow[r, "F_k"] & \cdots \arrow[r, "F_2"] & \bigoplus^{m_{1}}_{l=1} P^{r-\mathrm{a}^{(1)}_l} \arrow[r, "F_1"] & \bigoplus^{m_0}_{l=1} P^{r-\mathrm{a}_l} \arrow[r, "F_0"] & P^{r} \arrow[r] & P^r\mathcal{O}_{B} \arrow[r] & 0
\end{tikzcd}
\]
by summing the exact sequences for all \( r \), we have a resolution for \( j^+_*\mathcal{O}_{B^+} \) on \( X^+ \):
\[
\begin{tikzcd}
\cdots \arrow[r] & \bigoplus^{m_{2}}_{l=1} \mathcal{O}_{X^+}(\mathrm{a}^{(2)}_l E) \arrow[r, "F^+_2"] & \bigoplus^{m_{1}}_{l=1} \mathcal{O}_{X^+}(\mathrm{a}^{(1)}_l E) \arrow[r, "F^+_1"] & \bigoplus^{m_0}_{l=1} \mathcal{O}_{X^+}(\mathrm{a}^{(0)}_l E) \arrow[r, "F^+_0"] & \mathcal{O}_{X^+} \arrow[r] & j^+_*\mathcal{O}_{B^+} \arrow[r] & 0
\end{tikzcd}
\]
\begin{definition}
The above resolution of \( j^+_*\mathcal{O}_{B^+} \) is standard on \( X^+ \).
\end{definition}
We have an analogous picture for \( Y \) and \( Y^+ \), induced by the Rees algebra construction:
\[
\Pi: \mathrm{Gr}\big(\mathcal{R}(P/J)\big) \big/ \mathrm{Tors}\big(\mathcal{R}(P/J)\big) \xrightarrow{\sim} \mathrm{Coh}(Y^+):\Gamma,
\]
and there is also a resolution for \( j^+_*\mathcal{O}_{Y^+} \) on \( X^+ \):
\[
\begin{tikzcd}
0 \arrow[r] & \mathcal{O}_{X^+}({\mathrm{ord}_{P}} w \cdot E) \arrow[r] & \mathcal{O}_{X^+} \arrow[r] & j^+_*\mathcal{O}_{Y^+} \arrow[r] & 0
\end{tikzcd}
\]
\begin{lemma} \label{ratlemma}
Suppose \( Z \) has codimension \( c \) in \(X\) and \({\mathrm{ord}_{P}} w \) is \(d\) with \( c > d \). Then we have
\[
\pi_{Y*} \mathcal{O}_{Y^+} \simeq \mathcal{O}_Y.
\]
\end{lemma}

\begin{proof}
According to the exact sequence above, it suffices to show that \( \pi_{*}\mathcal{O}_{X^+}(dE) = \mathcal{O}_X \). This follows by induction using the short exact sequences
\[
0 \longrightarrow \mathcal{O}_{X^+}((j-1)E) \longrightarrow \mathcal{O}_{X^+}(jE) \longrightarrow \mathcal{O}_E(-j) \longrightarrow 0
\]
for \( 0 \leq j \leq d \), and the fact that \( \pi_{E*} \mathcal{O}_E(-j) = 0 \) for all \( 1 \leq j \leq c-1 \). This proves the claim.
\end{proof}
If we consider the standard resolution of \( P^r \mathcal{O}_{B} \) on \(Y\) as in Section \ref{subsec:hyp}:
\[
\begin{tikzcd}[column sep=normal]
\cdots \arrow[r] &
\mathbf{V}^r_{3} \arrow[r, "\mathbf{F}_{3}"] &
\mathbf{V}^r_2 \arrow[r, "\mathbf{F}_{2}"] &
\mathbf{V}^r_1  \arrow[r, "\mathbf{F}_1"] &
\mathbf{V}^r_0 \arrow[r, "\mathbf{F}_0"] &
P^r \arrow[r, "\mathbf{F}_{-1}"] &
P^r \mathcal{O}_{B} \arrow[r] &
0
\end{tikzcd}
\]
we then construct its direct sum resolution for arbitrary \( r \) as resolution  for \( \mathcal{O}_{B^+} \) on \( Y^+ \):
\[
\begin{tikzcd}[column sep=normal]
\cdots \arrow[r] &
\mathcal{V}_{3} \arrow[r, "\mathcal{F}_{3}"] &
\mathcal{V}_2 \arrow[r, "\mathcal{F}_{2}"] &
\mathcal{V}_1  \arrow[r, "\mathcal{F}_1"] &
\mathcal{V}_0 \arrow[r, "\mathcal{F}_0"] &
\mathcal{O}_{X^+}  \arrow[r, "\mathcal{F}_{-1}"] &
\mathcal{O}_{B^+} \arrow[r] &
0
\end{tikzcd}
\]
\begin{definition}
The above resolution of \( \mathcal{O}_{B^+} \) is standard on \( Y^+ \) .
\end{definition}

If \( M \) is a Cohen–Macaulay module over \( S \) with Bourbaki cycle \(B\), we consider its standard resolution:
\[
\begin{tikzcd}[column sep=normal]
\cdots \arrow[r] &
\mathbf{U}^r_4 \arrow[r, "\mathbf{\bar{F}}_4"] &
\mathbf{U}^r_3 \arrow[r, "\mathbf{\bar{F}}_3"] &
\mathbf{U}^r_2 \arrow[r, "\mathbf{\bar{F}}_2"] &
\mathbf{U}^r_1 \arrow[r, "\mathbf{\bar{F}}_1"] &
\mathbf{U}^r_0 \arrow[r, "\mathbf{\bar{F}}_0"] &
\mathbf{M}^r \arrow[r] &
0
\end{tikzcd}
\]
Similarly, we construct its direct sum resolution for arbitrary \( r \), which gives rise to an acyclic complex over \( Y^+ \). This complex is isomorphic to a sheaf over \( Y^+ \), which we denote by \( \mathcal{M} \):
\[
\begin{tikzcd}[column sep=normal]
\cdots \arrow[r] &
\mathcal{U}_4 \arrow[r, "\mathcal{\bar{F}}_4"] &
\mathcal{U}_3 \arrow[r, "\mathcal{\bar{F}}_3"] &
\mathcal{U}_2 \arrow[r, "\mathcal{\bar{F}}_2"] &
\mathcal{U}_1 \arrow[r, "\mathcal{\bar{F}}_1"] &
\mathcal{U}_0 \arrow[r, "\mathcal{\bar{F}}_0"] &
\mathcal{M} \arrow[r] &
0
\end{tikzcd}
\]
\begin{definition}
The above resolution of \( \mathcal{M}\) is standard on \( Y^+ \) .
\end{definition}

Next, we consider the normal cone of \( P \) that is \( \operatorname{gr}_{P/J}(R/J) \). By definition, we have the following exact sequence:
\[
0 \longrightarrow \operatorname{gr}_P(J,R)=L(w) \cdot \operatorname{gr}_P(R) \longrightarrow \operatorname{gr}_P(R) \longrightarrow \operatorname{gr}_{P/J}(R/J) \longrightarrow 0.
\]
According to the construction, we define
\[
E := \operatorname{Proj}_Z \operatorname{gr}_P(R), \quad E_{Y} := \operatorname{Proj}_Z \operatorname{gr}_{P/J}(R/J).
\]
There are similar isomorphisms of abelian categories :
\[
\Pi:  \mathrm{QGr}\big(\operatorname{gr}_P(R)\big) :=\mathrm{Gr}\big(\operatorname{gr}_P(R)\big) \big/ \mathrm{Tors}\big(\operatorname{gr}_P(R)\big) \simeq \mathrm{Coh}(E) : \Gamma,
\]
or
\[
\Pi:  \mathrm{QGr}\big( \operatorname{gr}_{P}(S)\big) :=\mathrm{Gr}\big( \operatorname{gr}_{P/J}(R/J)\big) \big/ \mathrm{Tors}\big( \operatorname{gr}_{P/J}(R/J)\big) \simeq \mathrm{Coh}(E_Y) : \Gamma.
\]
Furthermore, the above equivalence induces an equivalence between the corresponding derived categories:
\[
\mathrm{D^b}\left(\mathrm{QGr}(\operatorname{gr}_{P}(R))\right) \simeq \mathrm{D^b}(\mathrm{Coh}(E)),
\]
or
\[
\mathrm{D^b}\left(\mathrm{QGr}(\operatorname{gr}_{P}(S))\right) \simeq \mathrm{D^b}(\mathrm{Coh}(E_Y)).
\]
We assume that the dimension of \( R \) is \( n \), the codimension of \( R/P \) in \( R \) is \( c \), and the order of \( J \) at \( P \), denoted \( \mathrm{ord}_P(J) \), is \( d \). Consequently, the dimension of \( S \) is \( n - 1 \), and the codimension of \( S/P \) in \( S \) is \( c - 1 \). Therefore, the relative Gorenstein parameter \( a(E_Y) \) of \( E_Y \) viewed as a flat hypersurface fibration over \(Z\) is \( c - d \).

\begin{prop}[\cite{Orlov2009} {Thm. 2.5}, \cite{BurkeStevenson2015} {Thm. 6.4}]\label{orlovlemma}
If \( a(E_S) \geq 0 \), then the bounded derived category
\[
\mathrm{D}^b\left(\mathrm{QGr}\left(\operatorname{gr}_P(S)\right)\right)
\]
admits a semi-orthogonal decomposition of the form :
\[
\mathrm{D}^b\left(\mathrm{QGr}\left(\operatorname{gr}_P(S)\right)\right) =
\left\langle
\mathrm{D}^b(\mathrm{Coh}Z) \otimes_{\mathcal{O}_Z} \mathcal{A}(i - a + 1),\
\dots,\
\mathrm{D}^b(\mathrm{Coh}Z) \otimes_{\mathcal{O}_Z} \mathcal{A}(i),\
\mathcal{A}_i
\right\rangle
\]\label{equ:SODhyper}
\noindent where \( \mathcal{A}(i) \) denotes the \( i \)-th shift of the graded \( \mathcal{O}_Z \)-algebra \( \operatorname{gr}_P(S) \) composed with projection to \( \mathrm{QGr}\left(\operatorname{gr}_P(S)\right) \).
\end{prop}

\begin{remark}
We briefly recall the history of the above theorem. Orlov \cite{Orlov2009} first established it for graded algebras \(A\) over an arbitrary field \(k\) satisfying \(A_0=k\). Later, following the same philosophy, Burke and Stevenson \cite{BurkeStevenson2015} proved a substantially more general version for \(\mathbb{N}\)-graded noetherian rings
\[
A=\bigoplus_{i\ge 0} A_i,
\]
where \(A_0\) is any ring of finite global dimension. In particular, their result does not require \(A\) to be commutative. The latter setting relies on the notion of an Artin--Schelter Gorenstein ring. Namely, one assumes
\[
\operatorname{inj.dim}_A A < \infty \text{ and } \mathbf{R}\!\operatorname{Hom}_{\mathrm{gr}\,A}(A_0, A) \simeq A_0[n](a) \text{ for some } n, a \in \mathbb{Z}.
\]
In our situation,
\[
A=A_0[x_1,\ldots,x_c]/(L(w))
\]
with its natural grading and \(A_0\) a regular local ring. These conditions are automatically satisfied. Moreover, \(A\) has an associated Artin--Schelter index \(a\). When \(a>0\), after base change to the generic point of \(\operatorname{Spec}A_0\), this integer coincides with the Fano index of the corresponding generic hypersurface fiber. In fact, the normally flat assumption in Definition \ref{permissiblecenter} is not needed.
\end{remark}
In particular, we have a description of the residual category \( \mathcal{A}_i \). We consider the derived category of graded \( \mathcal{A} \)-modules whose homogeneous components are concentrated in degrees \( \geq -i \):
\[
\mathrm{D}^b\left(\mathrm{Gr}_{\geq -i}\,\mathcal{A}\right) = \mathrm{D}^b\left(\mathrm{Gr}_{\geq -i}\left(\operatorname{gr}_P(S)\right)\right),
\]
and we obtain the following weak semi-orthogonal decomposition of \( \mathrm{D}^b\left(\mathrm{Gr}_{\geq -i}\,\mathcal{A}\right) \) by \cite[Lemma 2.4]{Orlov2009}:
\[
\mathrm{D}^b\left(\mathrm{Gr}\,\mathcal{A}_{\geq -i}\right) =
\left\langle
\cdots,\mathrm{D}^b(Z) \otimes_{\mathcal{O}_Z} \mathcal{A}(i-2),\mathrm{D}^b(Z) \otimes_{\mathcal{O}_Z} \mathcal{A}(i-1),\mathrm{D}^b(Z) \otimes_{\mathcal{O}_Z} \mathcal{A}(i),\
\mathcal{T}_{-i}
\right\rangle,
\]
the residual category \( \mathcal{T}_{-i} \), together with its projection to \(\mathrm{D}^b\left(\mathrm{QGr}\left(\operatorname{gr}_P(S)\right)\right)\)
defines an equivalence of triangulated categories:
\begin{equation}\label{Orlovequ}
\mathcal{T}_{-i} \simeq \mathcal{A}_i.
\end{equation}
The key idea in Orlov's proof lies in the application of Serre dualizing functor together with a weak semi-orthogonal decomposition of \( \mathrm{D}^b\left(\mathrm{Gr}\,\mathcal{A}\right) \). The essential argument, informally speaking, is based on a ``division" principle: the behavior of the projective resolution on the ``tail" (i.e., in lower degrees) of modules on \( \mathrm{Gr}\,\mathcal{A}_{\geq -i} \) and \( \mathrm{QGr}\,\mathcal{A} \) coincides.\\

If we consider the standard resolution of \( P^r \mathcal{O}_{B} \) on \( S \):
\[
\begin{tikzcd}[column sep=normal]
\cdots \arrow[r] &
\mathbf{V}^r_{3} \arrow[r, "\mathbf{F}_{3}"] &
\mathbf{V}^r_2 \arrow[r, "\mathbf{F}_{2}"] &
\mathbf{V}^r_1 \arrow[r, "\mathbf{F}_1"] &
\mathbf{V}^r_0 \arrow[r, "\mathbf{F}_0"] &
P^r \arrow[r, "\mathbf{F}_{-1}"] &
P^r \mathcal{O}_{B} \arrow[r] &
0
\end{tikzcd}
\]
it induces via taking leading terms the following exact sequence:
\[
\begin{tikzcd}[column sep=normal]
&\cdots
\arrow[r] &
\mathbf{V}^r_2 / \mathbf{V}^{r+1}_2 \arrow[r, "L(\mathbf{F}_2)"] &
\mathbf{V}^r_1 / \mathbf{V}^{r+1}_1 \arrow[r, "L(\mathbf{F}_1)"] &
\mathbf{V}^r_0 / \mathbf{V}^{r+1}_0 \arrow[r, "L(\mathbf{F}_0)"] &
P^r / P^{r+1} \arrow[r, "L(\mathbf{F}_{-1})"] &
P^r\cdot\mathcal{O}_{B} / P^{r+1} \arrow[r] &
0
\end{tikzcd}
\]
taking the direct sum over all \( r \), we obtain the following exact sequence in \( \mathrm{D}^b(\mathrm{Gr}\, \mathcal{A}) \) which we denote as:
\[
\begin{tikzcd}[column sep=normal]
\cdots \arrow[r] &
\mathcal{U}_2 \arrow[r, "L(\mathcal{F}_2)"] &
\mathcal{U}_1 \arrow[r, "L(\mathcal{F}_1)"] &
\mathcal{U}_0 \arrow[r, "L(\mathcal{F}_0)"] &
\mathcal{A} \arrow[r, "L(\mathcal{F}_{-1})"] &
L(B^+)\cdot\mathcal{A} \arrow[r] &
0
\end{tikzcd}
\]
If \( M \) is a Cohen–Macaulay module over \( S \) with Bourbaki cycle \(B\), we consider its standard resolution:
\[
\begin{tikzcd}[column sep=normal]
\cdots \arrow[r] &
\mathbf{U}^r_4 \arrow[r, "\mathbf{\bar{F}}_4"] &
\mathbf{U}^r_3 \arrow[r, "\mathbf{\bar{F}}_3"] &
\mathbf{U}^r_2 \arrow[r, "\mathbf{\bar{F}}_2"] &
\mathbf{U}^r_1 \arrow[r, "\mathbf{\bar{F}}_1"] &
\mathbf{U}^r_0 \arrow[r, "\mathbf{\bar{F}}_0"] &
\mathbf{M}^r \arrow[r] &
0
\end{tikzcd}
\]
it induces via taking leading terms the following exact sequence:

\[
\begin{tikzcd}[column sep=normal]
\cdots \arrow[r] &
\mathbf{U}^r_3 / \mathbf{U}^{r+1}_3 \arrow[r, "L(\mathbf{\bar{F}}_3)"] &
\mathbf{U}^r_2 / \mathbf{U}^{r+1}_2 \arrow[r, "L(\mathbf{\bar{F}}_2)"] &
\mathbf{U}^r_1 / \mathbf{U}^{r+1}_1 \arrow[r, "L(\mathbf{\bar{F}}_1)"] &
\mathbf{U}^r_0 / \mathbf{U}^{r+1}_0 \arrow[r, "L(\mathbf{\bar{F}}_0)"] &
\mathbf{M}^r/\mathbf{M}^{r+1}\arrow[r] &
0
\end{tikzcd}
\]
similarly we obtain the following exact sequence in \( \mathrm{D}^b(\mathrm{Gr}\, \mathcal{A}) \):
\[
\begin{tikzcd}[column sep=normal]
\cdots \arrow[r] &
\mathcal{U}_4 \arrow[r, "L(\mathcal{\bar{F}}_4)"] &
\mathcal{U}_3 \arrow[r, "L(\mathcal{\bar{F}}_3)"] &
\mathcal{U}_2 \arrow[r, "L(\mathcal{\bar{F}}_2)"] &
\mathcal{U}_1 \arrow[r, "L(\mathcal{\bar{F}}_1)"] &
\mathcal{U}_0 \arrow[r, "L(\mathcal{\bar{F}}_0)"] &
\operatorname{gr}\mathcal{M} \arrow[r] &
0
\end{tikzcd}
\]
where we may assume that
\[
\mathcal{U}_0 := \bigoplus_{i=1}^m \mathcal{A}(a^{(0)}_{i}), \quad
\mathcal{U}_1 := \bigoplus_{i=1}^m \mathcal{A}(a^{(1)}_{i}).
\]
Moreover, we have a forward extension of the exact sequence above, which remains acyclic:
\[
\begin{tikzcd}[column sep=normal]
\cdots \arrow[r] &
\mathcal{U}_4 \arrow[r, "L(\mathcal{\bar{F}}_4)"] &
\mathcal{U}_3 \arrow[r, "L(\mathcal{\bar{F}}_3)"] &
\mathcal{U}_2 \arrow[r, "L(\mathcal{\bar{F}}_2)"] &
\mathcal{U}_1 \arrow[r, "L(\mathcal{\bar{F}}_1)"] &
\mathcal{U}_0 \arrow[r, "L(\mathcal{\bar{F}}_0)"] &
\mathcal{U}_{-1} \arrow[r, "L(\mathcal{\bar{F}}_{-1})"] &
\mathcal{U}_{-2} \arrow[r, "L(\mathcal{\bar{F}}_{-2})"] &
\cdots
\end{tikzcd}
\]
For sufficiently large \( j \), we have
\[
\mathcal{U}_{-j} = \bigoplus_{i=1}^m \mathcal{A}(a^{(-j)}_{i})
\quad \text{with } a^{(-j)}_{i} > 0.
\]
so we may consider the truncation of the above complex within
\(\mathrm{D}^b\left(\mathrm{Gr}_{\geq -i}\,\mathcal{A}\right)\), that is the subcomplex whose terms are
\[
\mathcal{U}^0_{j} := \bigoplus_{i=1}^m \mathcal{A}(a^{(j)}_{i}) \quad \text{with } a^{(j)}_{i} \leq 0,
\]
and where all maps are the corresponding restrictions from the original complex. We denote this truncated complex as \(\left\lfloor \operatorname{gr} \mathcal{M} \right\rfloor_{\geq 0}\) and the left complex after truncation by \(\left\lceil \operatorname{gr} \mathcal{M} \right\rceil_{<0}\). Similarly, for any integer \( k \), we define
\[
\mathcal{U}^k_{j} := \bigoplus_{i=1}^m \mathcal{A}(a^{(j)}_{i}) \quad \text{with } a^{(j)}_{i} \leq -k,
\]
as it is a truncation of an complex by direct summands and denote this truncated complex by a  \(\left\lfloor \operatorname{gr} \mathcal{M} \right\rfloor_{\geq k}\) and left complex after truncation by \( \left\lceil \operatorname{gr} \mathcal{M} \right\rceil_{\leq k} \). In particular, we have an isomorphism in the derived category \( \mathrm{D}^b\left(\mathrm{Gr}\,\mathcal{A}\right) \):
\[
\left\lfloor \operatorname{gr} \mathcal{M} \right\rfloor_{\geq k} \simeq \left\lceil \operatorname{gr} \mathcal{M} \right\rceil_{<k}[1].
\]
\begin{lemma}\label{MFlemma}
    We have
    \[
        \left\lfloor \operatorname{gr} \mathcal{M} \right\rfloor_{\geq k} \in \mathcal{T}_k,
    \]
    and thus its projection belongs to $\mathcal{A}_{-k}$ for any integer $k$.
\end{lemma}

\begin{proof}
We consider the naive truncation for \( \left\lfloor \operatorname{gr} \mathcal{M} \right\rfloor_{\geq k} \):
\[
P_{<n}\left\lfloor \operatorname{gr} \mathcal{M} \right\rfloor_{\geq k} \longrightarrow \left\lfloor \operatorname{gr} \mathcal{M} \right\rfloor_{\geq k} \longrightarrow \sigma_{\geq n}\left\lfloor \operatorname{gr} \mathcal{M} \right\rfloor_{\geq k},
\]
for \( n \gg 0 \). Then the truncated complex $\sigma_{\geq n}\left\lfloor \operatorname{gr} \mathcal{M} \right\rfloor_{\geq k}$
is a subcomplex of the kernel of morphisms between $\mathcal{U}^{k}_{n-1}$ and $\mathcal{U}^{k}_{n}$
such that all direct summands $\mathcal{A}(a_i^{(n-1)})$ (respectively $\mathcal{A}(a_i^{(n)})$)
satisfy $a_i^{(n-1)} < -k$ (respectively $a_i^{(n)} < -k$). On the other hand, the truncation \( P_{<n}\left\lfloor \operatorname{gr} \mathcal{M} \right\rfloor_{\geq k} \) is a finite-length complex whose terms lie in \( \mathrm{Gr}_{\geq k}\mathcal{A} \). Therefore, we conclude that \(\left\lfloor \operatorname{gr} \mathcal{M} \right\rfloor_{\geq k} \in \mathrm{D}^b\left(\mathrm{Gr}_{\geq k}\mathcal{A}\right)\).\\

We consider the naive truncation of \( \left\lceil \operatorname{gr} \mathcal{M} \right\rceil_{<k} \) for \( n \ll 0 \):
\[
\sigma_{\leq n} \left\lceil \operatorname{gr} \mathcal{M} \right\rceil_{<k}
\longrightarrow \left\lceil \operatorname{gr} \mathcal{M} \right\rceil_{<k}
\longrightarrow P_{>n} \left\lceil \operatorname{gr} \mathcal{M} \right\rceil_{<k}.
\]
The truncation \( P_{>n} \left\lceil \operatorname{gr} \mathcal{M} \right\rceil_{<k} \) is a bounded complex whose terms are direct sums of \( \mathcal{A}(a) \) with \( a \geq -k \). Therefore, for any \( b < -k \), we have
\[
\operatorname{Ext}^*\left( P_{>n} \left\lceil \operatorname{gr} \mathcal{M} \right\rceil_{<k},\, \mathcal{A}(b) \right) = 0.
\]
On the other hand, let us denote \( N := \sigma_{\leq n} \left\lceil \operatorname{gr} \mathcal{M} \right\rceil_{<k} \) then \( N \) admits a short exact sequence of graded \( \mathcal{A} \)-modules:
\[
0 \longrightarrow \operatorname{Syz}(N) \longrightarrow \mathcal{U}_{j}^{n-1}= \bigoplus_{i=1}^m \mathcal{A}\left( a^{(n-1)}_i \right) \longrightarrow N \longrightarrow 0,
\]
where \( a^{(n-1)}_i > -k \) for all \( i \). Applying \( \operatorname{Hom}(-,\mathcal{A}(b)) \) for any \( b < -k \), and using the long exact sequence, we obtain:
\[
0 \longrightarrow \operatorname{Hom}(N, \mathcal{A}(b)) \longrightarrow \operatorname{Hom}(\mathcal{U}_j^{n-1}, \mathcal{A}(b)) = 0\longrightarrow \cdots.
\]
Hence \( \operatorname{Hom}(N, \mathcal{A}(b)) = 0 \) for all \( b < -k \). Furthermore, since the projective resolution of \( N \) has dual projective resolution which is also acyclic, we conclude that \( \operatorname{Ext}^{>0}(N, \mathcal{A}(b)) = 0 \) for all \( b < -k \).
Combining both sides of the truncation, we conclude:
\[
\operatorname{Ext}^*\left( \left\lceil \operatorname{gr} \mathcal{M} \right\rceil_{<k},\, \mathcal{A}(b) \right) = 0 \quad \text{for all } b < -k
\]
Finally, observe that \( Z \) is defined by a regular local ring, so every object in \( \mathrm{D}^b(Z) \) admits a finite projective (hence free) resolution. The desired result then follows by equivalence (\ref{Orlovequ}).

\end{proof}

We have the following geometric description:
\[
\mathrm{D}^b\left(\mathrm{Coh}(E_Y)\right) =
\left\langle
\pi^*_{E_Y}\mathrm{D}^b(\mathrm{Coh}(Z)) (i - a + 1),\
\dots,\
\pi^*_{E_Y}\mathrm{D}^b(\mathrm{Coh}(Z)) (i),\
\mathcal{A}_i
\right\rangle.
\]
Since all the components in the above semi-orthogonal decomposition are admissible, we may apply mutation to obtain the following equivalent decomposition:
\[
\mathrm{D}^b\left(E_Y\right) =
\left\langle
\pi^*_{E_Y} \mathrm{D}^b(Z)(i - a + 1 + b),\
\dots,\
\pi^*_{E_Y} \mathrm{D}^b(Z)(i),\
\mathcal{A}_i,\
\pi^*_{E_Y} \mathrm{D}^b(Z)(i+1),\
\dots,\
\pi^*_{E_Y} \mathrm{D}^b(Z)(i+b)
\right\rangle.
\]
When we consider $\mathcal{O}_{E_Y}(k)$ as a pre-tilting object with tilting algebra $\mathcal{O}_{Z}$, we have the following standard description:
\[
\mathrm{D}^b\left(E_Y\right) =
\left\langle
\mathcal{O}_{E_Y}(i - a + 1 + b),\
\dots,\
\mathcal{O}_{E_Y}(i),\
\mathcal{A}_i,\
\mathcal{O}_{E_Y}(i+1),\
\dots,\
\mathcal{O}_{E_Y}(i+b)
\right\rangle.
\]
Moreover, for any integer \( i \), we have \( \left\lfloor \operatorname{gr} \mathcal{M} \right\rfloor_{\geq i} \in \mathcal{A}_i \). Under the assumption that the Gorenstein parameter of $E_Y$ is positive, which means $c-d>0$, we will use the following semi-orthogonal decomposition:
\[
\mathrm{D}^b\left(E_Y\right) =
\left\langle
\mathcal{A}_{-1},\,\mathcal{O}_{E_Y},\
\mathcal{O}_{E_Y}(1),\
\dots,\
\mathcal{O}_{E_Y}(c-d-1)
\right\rangle.
\]
and its dual decomposition:
\[
\mathrm{D}^b\left(E_Y\right) =
\left\langle
\mathcal{A}_{-1},\,\
\mathcal{D}^Y_{c-d-1},\
\dots,\
\mathcal{D}^Y_{1},\
\mathcal{D}^Y_{0}
\right\rangle.
\]
where we define the dual object as mutation of pre-tilting objects, i.e.
\[
\mathcal D^Y_{i}:=\mathds{L}_{\langle\mathcal{O}_{E_Y},...,\mathcal{O}_{E_Y}(i-1)\rangle}\mathcal{O}_{E_Y}(i), \quad 0\leq i\leq c-d-1.
\]
Next, we compute some simple examples. Suppose our cycle \( B \) is the permissible center \( Z \). According to the Koszul resolution of \( \mathcal{O}_{Z} \) on  \(X\), each term in the resolution has \( \operatorname{ord}_P = 1 \). Hence, we obtain a resolution of \( P^r \mathcal{O}_{Z} \) in which each term of morphisms has \( \operatorname{ord}_P = 1 \), we have the standard resolution of \( P^r \mathcal{O}_{Z} \) over \( X \):
\[
\begin{tikzcd}
0 \arrow[r] & P^{r-c} \otimes_R \wedge^c V_c \arrow[r, "F_{c-1}"] & \cdots \arrow[r, "F_2"] & P^{r-2} \otimes_R \wedge^2 V_c \arrow[r, "F_1"] & P^{r-1} \otimes_R V_c \arrow[r, "F_0"] & P^r \arrow[r] & P^r \mathcal{O}_{Z} \arrow[r] & 0,
\end{tikzcd}
\]
where \( V_c \) is the vector space of codimension $P$, note that if \( r > 0 \), then \( P^r \mathcal{O}_{Z} = 0 \), since \(P^r\) vanishes on its support.
By taking the direct sum over all \( r \geq 0 \) of the exact sequences above, we obtain an exact sequence in the category of \( \mathrm{Gr}\big(\mathcal{R}(P)\big) \)-modules:
\[
\begin{tikzcd}
0 \arrow[r] &
\bigoplus_{r \geq 0} P^{r-c} \otimes_R \wedge^c V_c \arrow[r, "F_{c-1}"] & \cdots \arrow[r] &
\bigoplus_{r \geq 0} P^{r-1} \otimes_R V_c \arrow[r, "F_0"] &
\bigoplus_{r \geq 0} P^r \arrow[r] &
\mathcal{O}_{Z} \arrow[r] & 0,
\end{tikzcd}
\]
Then by projecting the above exact sequence to the quotient category
\(\mathrm{QGr}\big(\mathcal{R}(P)\big) \simeq \mathrm{Coh}(X^+)\),
we obtain the following exact sequence:
\[
\begin{tikzcd}
0 \arrow[r] &
\wedge^c V_c \otimes \mathcal{O}_{X^+}(cE) \arrow[r, "F_{c-1}"] &
\cdots \arrow[r] &
V_c \otimes \mathcal{O}_{X^+}(E) \arrow[r, "F_0"] &
\mathcal{O}_{X^+} \arrow[r] & 0,
\end{tikzcd}
\]
Restricting the above sequence to \( \mathrm{Coh}(E) \), we obtain the following exact sequence:
\begin{equation}\label{eulerexact}
\begin{tikzcd}
0 \arrow[r] &
\wedge^c V_c \otimes \mathcal{O}_{E}(-c) \arrow[r, "F_{c-1}"] &
\cdots \arrow[r] &
V_c \otimes \mathcal{O}_{E}(-1) \arrow[r, "F_0"] &
\mathcal{O}_{E} \arrow[r] & 0.
\end{tikzcd}
\end{equation}
Specifically, we consider the following commutative diagram in complex category of \(R\) module:
\[
\begin{tikzcd}
0 \arrow[r] &
P^{r-c} \otimes \wedge^c V_c \arrow[r, "F_{c-1}"] &
\cdots \arrow[r, "F_2"] &
P^{r-2} \otimes \wedge^2 V_c \arrow[r, "F_1"] &
P^{r-1} \otimes V_c \arrow[r, "F_0"] &
P^r \arrow[r] &
P^r \mathcal{O}_{Z} \arrow[r] &
0 \\
0 \arrow[r] &
P^{r-c+1} \otimes \wedge^c V_c \arrow[r, "F_{c-1}"] \arrow[u, hook] &
\cdots \arrow[r, "F_2"] &
P^{r-1} \otimes \wedge^2 V_c \arrow[r, "F_1"] \arrow[u, hook] &
P^{r} \otimes V_c \arrow[r, "F_0"] \arrow[u, hook] &
P^r \arrow[r] \arrow[u, equal] &
P^r \mathcal{O}_{Z} \arrow[r] \arrow[u, equal] &
0 \\
0 \arrow[r] &
P^{r-c+2} \otimes \wedge^c V_c \arrow[r, "F_{c-1}"] \arrow[u, hook] &
\cdots \arrow[r, "F_2"] &
P^r \otimes \wedge^2 V_c \arrow[r, "F_1"] \arrow[u, hook] &
P^r \otimes V_c \arrow[r, "F_0"] \arrow[u, equal] &
P^r \arrow[r] \arrow[u, equal] &
P^r \mathcal{O}_{Z} \arrow[r] \arrow[u, equal] &
0 \\
\vdots  &
\vdots \arrow[u] &
& \vdots \arrow[u] &
\vdots \arrow[u] &
\vdots \arrow[u] &
\vdots \arrow[u] &
\vdots \\
0 \arrow[r] &
P^{r} \otimes \wedge^c V_c \arrow[r, "F_{c-1}"] \arrow[u, hook] &
\cdots \arrow[r, "F_2"] &
P^r \otimes \wedge^2 V_c \arrow[r, "F_1"] \arrow[u, equal] &
P^r \otimes V_c \arrow[r, "F_0"] \arrow[u, equal] &
P^r \arrow[r] \arrow[u, equal] &
P^r \mathcal{O}_{Z} \arrow[r] \arrow[u, equal] &
0
\end{tikzcd}
\]
Note that each row in the diagram is a complex, and the first row is exact. For any \( r \geq 0 \), we take the direct sum and then project to \( \mathrm{Coh}(X^+) \). We obtain the following commutative diagram:
\[
\begin{tikzcd}
0 \arrow[r] &
\wedge^c V_c \otimes \mathcal{O}_{R^+}(cE) \arrow[r, "F_{c-1}"] &
\cdots \arrow[r, "F_2"] &
\wedge^2 V_c \otimes \mathcal{O}_{R^+}(2E) \arrow[r, "F_1"] &
V_c \otimes \mathcal{O}_{R^+}(E) \arrow[r, "F_0"] &
\mathcal{O}_{R^+} \arrow[r] &
0 \\
0 \arrow[r] &
\wedge^c V_c \otimes \mathcal{O}_{R^+}((c-1)E) \arrow[r, "F_{c-1}"] \arrow[u, hook] &
\cdots \arrow[r, "F_2"] &
\wedge^2 V_c \otimes \mathcal{O}_{R^+}(E) \arrow[r, "F_1"] \arrow[u, hook] &
V_c \otimes \mathcal{O}_{R^+} \arrow[r, "F_0"] \arrow[u, hook] &
\mathcal{O}_{R^+} \arrow[r] \arrow[u, equal] &
0 \\
0 \arrow[r] &
\wedge^c V_c \otimes \mathcal{O}_{R^+}((c-2)E) \arrow[r, "F_{c-1}"] \arrow[u, hook] &
\cdots \arrow[r, "F_2"] &
\wedge^2 V_c \otimes \mathcal{O}_{R^+} \arrow[r, "F_1"] \arrow[u, hook] &
V_c \otimes \mathcal{O}_{R^+} \arrow[r, "F_0"] \arrow[u, equal] &
\mathcal{O}_{R^+} \arrow[r] \arrow[u, equal] &
0 \\
\vdots  &
\vdots \arrow[u] &
& \vdots \arrow[u] &
\vdots \arrow[u] &
\vdots \arrow[u] &
\\
0 \arrow[r] &
\wedge^c V_c \otimes \mathcal{O}_{R^+} \arrow[r, "F_{c-1}"] \arrow[u, hook] &
\cdots \arrow[r, "F_2"] &
\wedge^2 V_c \otimes \mathcal{O}_{R^+} \arrow[r, "F_1"] \arrow[u, equal] &
V_c \otimes \mathcal{O}_{R^+} \arrow[r, "F_0"] \arrow[u, equal] &
\mathcal{O}_{R^+} \arrow[r] \arrow[u, equal] &
0
\end{tikzcd}
\]
Note that each row in the diagram is a complex, where the first row is exact and the last row is quasi-isomorphic to the derived pull-back \(\pi^*\mathcal{O}_{Z}\), specially the quotient of two adjacent rows \( i \) and \( i+1 \) has the following form:
\[
\begin{tikzcd}
0 \arrow[r] &
\mathcal{O}_E(-c+i-1) \otimes \wedge^c V_c \arrow[r, "L(F_{c-1})"] &
\cdots \arrow[r] &
\mathcal{O}_E(-2) \otimes \wedge^{i+1} V_c \arrow[r, "L(F_i)"] &
\mathcal{O}_E(-1) \otimes \wedge^i V_c \arrow[r] &
0
\end{tikzcd}
\]
Comparing with the Euler exact sequence (\ref{eulerexact}), we have that it is quasi-isomorphic to:
\[
\begin{tikzcd}
0 \arrow[r] &
\wedge^{i-1} V_c \otimes \mathcal{O}_{E}\arrow[r, "L(F_{i-2})"] &
\cdots \arrow[r] &
V_c \otimes \mathcal{O}_{E}(i-2) \arrow[r, "L(F_0)"] &
\mathcal{O}_{E}(i-1) \arrow[r] &
0,
\end{tikzcd}
\]
We can easily see that it is isomorphic to the $(i-1)$-th dual exceptional (pre-tilting) collection $\mathcal{D}_{i-1}$ on $E=\mathbb{P}_Z(V_c)$, and we have the following classical result:
The derived pushforward $\pi^*\mathcal{O}_{Z}$ has a unique right Postnikov filtration in $\mathrm{D}^b(X^+)$:
$$\begin{tikzcd}[column sep=0.5em]
 & \pi^*\mathcal{O}_{Z} \arrow{rr}&& E_{c-1}\arrow{dl}\arrow{rr}&& E_{c-2}  \arrow{dl} \\
&& k^+_{*}\mathcal{D}_{c-1}[1]\arrow[ul,dashed,"\Delta"]&& k^+_{*}\mathcal{D}_{c-2}[1]\arrow[ul,dashed,"\Delta"]
\end{tikzcd}
...\begin{tikzcd}[column sep=0.5em]
 &   E_{2}\arrow{rr}&&  E_{1}\arrow{rr}\arrow{dl}&& 0\arrow{dl} \\
 && k^+_{*}\mathcal{D}_{1}[1]\arrow[ul,dashed,"\Delta"]  && k^+_{*}\mathcal{D}_{0}[1]\arrow[ul,dashed,"\Delta"]
\end{tikzcd}$$
More familiarly, if we reverse the direction in which $r - i$ increases in the previous double complex—i.e., if we increase upwards—we obtain a similar left Postnikov filtration in $\mathrm{D}^b(X^+)$.\\

Third, we consider the previously constructed standard resolution of \( P^r \mathcal{O}_{Z} \) over \( Y \), which is induced by restriction of the totalization of the following double factorization over \( X \):
\[
\begin{tikzcd}[column sep=small, row sep=small]
 & \vdots \arrow[d] & \vdots \arrow[d] \\
 P^{r-2-2d} \otimes \wedge^2 V_c \arrow[r, "F_{1}"] \arrow[d, "-K_2"] &
 P^{r-1-2d} \otimes V_c \arrow[r, "F_{0}"] \arrow[d, "K_1"] &
 P^{r-2d} \arrow[ddl, dashed, "-K_0^{(1)}"] \arrow[d, "-K_0"] \\
 P^{r-3-d} \otimes \wedge^3 V_c \arrow[r, "F_{2}"] \arrow[d, "-K_3"] &
 P^{r-2-d} \otimes \wedge^2 V_c \arrow[r, "F_{1}"] \arrow[d, "K_2"] &
 P^{r-1-d} \otimes V_c \arrow[r, "F_{0}"] \arrow[d, "-K_1"] &
 P^{r-d} \arrow[d, "K_0"] \\
 \cdots \arrow[r] &
 P^{r-3} \otimes \wedge^3 V_c \arrow[r, "F_{2}"] &
 P^{r-2} \otimes \wedge^2 V_c \arrow[r, "F_{1}"] &
 P^{r-1} \otimes V_c \arrow[r, "F_{0}"] &
 P^{r}
\end{tikzcd}
\]
In this construction, we specifically note that:
\[ K^{(a)}_i = 0 \quad \text{for all } i \text{ and } a > 0
\]
we refer to \cite[Lemma 2.6]{Dyckerhoff2011},. Moreover, our standard resolution carries a natural filtration \( F^0(P^r \mathcal{O}_{Z}) \) induced by the factorization on $R$, which gives rise to a complex over $S$ through
\[
\begin{tikzcd}[column sep=small, row sep=small]
 & \vdots \arrow[d] & \vdots \arrow[d] \\
 P^{r-1-2d} \otimes \wedge^2 V_c \arrow[r, "F_{1}"] \arrow[d, "-K_2"] &
 P^{r-2d} \otimes V_c \arrow[r, "F_{0}"] \arrow[d, "K_1"] &
 P^{r+1-2d}  \arrow[d, "-K_0"] \\
 P^{r-2-d} \otimes \wedge^3 V_c \arrow[r, "F_{2}"] \arrow[d, "-K_3"] &
 P^{r-1-d} \otimes \wedge^2 V_c \arrow[r, "F_{1}"] \arrow[d, "K_2"] &
 P^{r-d} \otimes V_c \arrow[r, "F_{0}"] \arrow[d, "-K_1"] &
 P^{r+1-d} \arrow[d, "K_0"] \\
 \cdots \arrow[r] &
 P^{r-2} \otimes \wedge^3 V_c \arrow[r, "F_{2}"] &
 P^{r-1} \otimes \wedge^2 V_c \arrow[r, "F_{1}"] &
 P^{r} \otimes V_c \arrow[r, "F_{0}"] &
 P^{r}
\end{tikzcd}
\]
The quotients of the two complexes over $S$ above are given by the following  factorization:
\[
\begin{tikzcd}[column sep=small, row sep=small]
 & \vdots \arrow[d] & \vdots \arrow[d] \\
 P^{r-2-2d}/P^{r-1-2d} \otimes \wedge^2 V_c \arrow[r, "L(F_{1})"] \arrow[d, "-L(K_2)"] &
 P^{r-1-2d}/P^{r-2d} \otimes V_c \arrow[r, "L(F_{0})"] \arrow[d, "L(K_1)"] &
 P^{r-2d}/P^{r+1-2d} \arrow[d, "-L(K_0)"] \\
 P^{r-3-d}/P^{r-2-d} \otimes \wedge^3 V_c \arrow[r, "L(F_{2})"] \arrow[d, "-L(K_3)"] &
 P^{r-2-d}/P^{r-1-d} \otimes \wedge^2 V_c \arrow[r, "L(F_{1})"] \arrow[d, "L(K_2)"] &
 P^{r-1-d}/P^{r-d} \otimes V_c \arrow[r, "L(F_{0})"] \arrow[d, "-L(K_1)"] &
 P^{r-d}/P^{r+1-d} \arrow[d, "L(K_0)"] \\
 \cdots \arrow[r] &
 P^{r-3}/P^{r-2} \otimes \wedge^3 V_c \arrow[r, "L(F_{2})"] &
 P^{r-2}/P^{r-1} \otimes \wedge^2 V_c \arrow[r, "L(F_{1})"] &
 P^{r-1}/P^{r} \otimes V_c \arrow[r] &
 0
\end{tikzcd}
\]
according to the previous construction, we know that it is isomorphic to $P^{r}/P^{r+1}$ over $S$. Using the same logic, we take the direct sum over all $r\geq 0$ and consider the natural projection to $\mathrm{D}^b(Y^+)$: we have morphism between \( F^0(\mathcal{O}_{\widetilde{Z}})\) which is also a trivial object:
\[
\begin{tikzcd}[column sep=small, row sep=small]
 & \vdots \arrow[d] & \vdots \arrow[d] \\
 \mathcal{O}_{Y^+}((2+2d)E) \otimes \wedge^2 V_c \arrow[r, "F_{1}"] \arrow[d, "-K_2"] &
 \mathcal{O}_{Y^+}((1+2d)E) \otimes V_c \arrow[r, "F_{0}"] \arrow[d, "K_1"] &
 \mathcal{O}_{Y^+}(2dE) \arrow[d, "-K_0"] \\
 \mathcal{O}_{Y^+}((3+d)E) \otimes \wedge^3 V_c \arrow[r, "F_{2}"] \arrow[d, "-K_3"] &
 \mathcal{O}_{Y^+}((2+d)E) \otimes \wedge^2 V_c \arrow[r, "F_{1}"] \arrow[d, "K_2"] &
 \mathcal{O}_{Y^+}((1+d)E) \otimes V_c \arrow[r, "F_{0}"] \arrow[d, "-K_1"] &
 \mathcal{O}_{Y^+}(dE) \arrow[d, "K_0"] \\
 \cdots \arrow[r] &
 \mathcal{O}_{Y^+}(3E) \otimes \wedge^3 V_c \arrow[r, "F_{2}"] &
 \mathcal{O}_{Y^+}(2E) \otimes \wedge^2 V_c \arrow[r, "F_{1}"] &
 \mathcal{O}_{Y^+}(E) \otimes V_c \arrow[r, "F_{0}"] &
 \mathcal{O}_{Y^+}
\end{tikzcd}
\]
and \(F^1(\mathcal{O}_{\widetilde{Z}})\):
\[
\begin{tikzcd}[column sep=small, row sep=small]
 & \vdots \arrow[d] & \vdots \arrow[d] \\
 \mathcal{O}_{Y^+}((1+2d)E) \otimes \wedge^2 V_c \arrow[r, "F_{1}"] \arrow[d, "-K_2"] &
 \mathcal{O}_{Y^+}(2dE) \otimes V_c \arrow[r, "F_{0}"] \arrow[d, "K_1"] &
 \mathcal{O}_{Y^+}((-1+2d)E) \arrow[d, "-K_0"] \\
 \mathcal{O}_{Y^+}((2+d)E) \otimes \wedge^3 V_c \arrow[r, "F_{2}"] \arrow[d, "-K_3"] &
 \mathcal{O}_{Y^+}((1+d)E) \otimes \wedge^2 V_c \arrow[r, "F_{1}"] \arrow[d, "K_2"] &
 \mathcal{O}_{Y^+}(dE) \otimes V_c \arrow[r, "F_{0}"] \arrow[d, "-K_1"] &
 \mathcal{O}_{Y^+}((-1+d)E) \arrow[d, "K_0"] \\
 \cdots \arrow[r] &
 \mathcal{O}_{Y^+}(2E) \otimes \wedge^3 V_c \arrow[r, "F_{2}"] &
 \mathcal{O}_{Y^+}(E) \otimes \wedge^2 V_c \arrow[r, "F_{1}"] &
 \mathcal{O}_{Y^+} \otimes V_c \arrow[r, "F_{0}"] &
 \mathcal{O}_{Y^+}
\end{tikzcd}
\]
such that their quotient  \(F^0/F^1(\mathcal{O}_{\widetilde{Z}})\) is:
\[
\begin{tikzcd}[column sep=small, row sep=small]
 & \vdots \arrow[d] & \vdots \arrow[d] \\
 \mathcal{O}_{E_Y}(-(2+2d)) \otimes \wedge^2 V_c \arrow[r, "L(F_{1})"] \arrow[d, "-L(K_2)"] &
 \mathcal{O}_{E_Y}(-(1+2d)) \otimes V_c \arrow[r, "L(F_{0})"] \arrow[d, "L(K_1)"] &
 \mathcal{O}_{E_Y}(-2d) \arrow[d, "-L(K_0)"] \\
 \mathcal{O}_{E_Y}(-(3+d)) \otimes \wedge^3 V_c \arrow[r, "L(F_{2})"] \arrow[d, "-L(K_3)"] &
 \mathcal{O}_{E_Y}(-(2+d)) \otimes \wedge^2 V_c \arrow[r, "L(F_{1})"] \arrow[d, "L(K_2)"] &
 \mathcal{O}_{E_Y}(-(1+d)) \otimes V_c \arrow[r, "L(F_{0})"] \arrow[d, "-L(K_1)"] &
 \mathcal{O}_{E_Y}(-d) \arrow[d, "L(K_0)"] \\
 \cdots \arrow[r] &
 \mathcal{O}_{E_Y}(-3) \otimes \wedge^3 V_c \arrow[r, "L(F_{2})"] &
 \mathcal{O}_{E_Y}(-2) \otimes \wedge^2 V_c \arrow[r, "L(F_{1})"] &
 \mathcal{O}_{E_Y}(-1) \otimes V_c \arrow[r] &
 0
\end{tikzcd}
\]
which is quasi-isomorphic to \(\mathcal{D}_{0}\), that's just \(\mathcal{O}_{E_Y}\). The complex is naturally induced by the following semi-orthogonal decomposition
\[
\mathrm{D}^b(E_Y) =
\left\langle
\mathcal{A}_{-(c-d)-1},\,
\mathcal{O}_{E_Y}(-(c-d)),\,
\dots,\,
\mathcal{O}_{E_Y}(-2),\,
\mathcal{O}_{E_Y}(-1)
\right\rangle.
\]
For  $\mathcal{O}_{E_Y}$, we have:
\begin{corollary}
For any $1 \leq k \leq c-d$, the mutation
$\mathbb{L}_{\left\langle \mathcal{O}_{E_Y}(-k),\dots,\mathcal{O}_{E_Y}(-1) \right\rangle}\,{\mathcal{O}_{E_Y}}$
is exactly the truncation of the above complex at graded degree $\geq k+1$. When $k = c-d$, this mutation belongs to $\mathcal{A}_{-(c-d)-1}$ by Lemma \ref{MFlemma}. More generally, this conclusion holds for all elements with an appropriate twist.
\end{corollary}
We iteratively apply above method to obtain a natural filtration \(F^s(\mathcal{O}_{\widetilde{Z}})\) for any positive integers $s$ as follows:
\[
\begin{tikzcd}[column sep=0.5em]
 & F^s(\mathcal{O}_{\widetilde{Z}})\arrow{rr}&& F^{s-1}\arrow{dl}\arrow{rr}&& F^{s-2}  \arrow{dl} \\
&& i^+_{*}\mathcal{D}^Y_{s}[1]\arrow[ul,dashed,"\Delta"]&& i^+_{*}\mathcal{D}^Y_{s-1}[1]\arrow[ul,dashed,"\Delta"]
\end{tikzcd}
\quad\cdots\quad
\begin{tikzcd}[column sep=0.5em]
 &   F^1\arrow{rr}&&  F^0\arrow{rr}\arrow{dl}&& 0\arrow{dl} \\
 && i^+_{*}\mathcal{D}^Y_{1}[1]\arrow[ul,dashed,"\Delta"]  &&i^+_{*}\mathcal{D}^Y_{0}[1]\arrow[ul,dashed,"\Delta"]
\end{tikzcd}
\]
In particular, for any \( i \leq c - d - 1 \), we have
\[
\mathcal{D}^Y_i
\]
is precisely the \( i \)-th dual exceptional (or pre-tilting) object that appeared earlier in the definition of \( \mathrm{D}^b(E_Y) \) as a hypersurface fibration. On the other hand, for \( i \geq c - d \), we have
\[
\mathcal{D}^Y_i \in \mathcal{A}_{-1}
\]
by the previous corollary. This is a purely combinatorial and straightforward verification, which we leave to the reader.\\
\begin{lemma}\label{ortha-1o}
Under the assumption \( c - d > 1 \), for any object \( N \in \mathcal{A}_{-1} \) we have
\[
\operatorname{Ext}^*_{Y^+}\left(i^+_*N,\ i^+_*\mathcal{O}_{E_Y}(-i)\right) = 0
\quad \text{for all } 1 \leq i \leq c - d - 1.
\]
\end{lemma}
\begin{proof}
The proof is similar. Consider the excess distinguished triangle:
\[
\mathcal{O}_{E_Y}(-i) \longrightarrow i^{+!}i^+_*\mathcal{O}_{E_Y}(-i) \longrightarrow \mathcal{O}_{E_Y}(-i)(E_Y)[-1].
\]
Also consider the semi-orthogonal decomposition of \( \mathrm{D}^b(E_Y) \):
\[
\mathrm{D}^b(E_Y) =
\left\langle
\mathcal{O}_{E_Y}(-(c-d)),\,
\dots,\,
\mathcal{O}_{E_Y}(-1),\,
\mathcal{A}_{-1}
\right\rangle
\]
Using the long exact sequence in cohomology induced from the distinguished triangle, we obtain the result.
\end{proof}
\begin{lemma}\label{localpushspan}
For all $s \geq c - d - 1$ and $1 \leq i \leq c - d - 1$, we have
\[
\operatorname{Ext}_{Y^+}^*\big(F^s(\mathcal{O}_{\widetilde{Z}}), i^+_*\mathcal{O}_{E_Y}(-i)\big) = 0.
\]
\end{lemma}
\begin{proof}
We need to show that for any integer $k$, $\operatorname{Ext}_{Y^+}^k\big(F^s, i^+_*\mathcal{O}_{E_{Y}}(-i)\big)=0$. Consider the following filtration for any $s'\gg s$:
\[
\begin{tikzcd}[column sep=0.5em]
 & \pi^*_{Y}\mathcal{O}_{Z} \arrow{rr}&& F^{s'}\arrow{rr}\arrow{dl} && F^s(\mathcal{O}_{\widetilde{Z}})\arrow{dl} \\
 && \tau^1\arrow[ul,dashed,"\Delta"] && \tau^2\arrow[ul,dashed,"\Delta"]
\end{tikzcd}
\]
Since $\tau^2$ is generated by a sequence of $i^+_* N$ where $N$ belongs to $\mathcal{A}_{-1}$, and according to the previous lemma, we have $\operatorname{Ext}_{Y^+}^*(\tau^2, i^+_*\mathcal{O}_{E_Y}(-i))=0$. Therefore, we obtain
\[
\operatorname{Ext}_{Y^+}^k\big(F^s, i^+_*\mathcal{O}_{E_Y}(-i)\big)=\operatorname{Ext}_{Y^+}^k\big(F^{s'}, i^+_*\mathcal{O}_{E_Y}(-i)\big).
\]
On the other hand, since
\[
\operatorname{Ext}_{Y^+}^*\big(\pi^*_{Y}\mathcal{O}_{Z}, i^+_*\mathcal{O}_{E_Y}(-i)\big)=\operatorname{Ext}_{Y}^*\big(\mathcal{O}_{Z}, \pi_{Y*}i^+_*\mathcal{O}_{E_Y}(-i)\big)=0,
\]
and for $s'\gg s$ we have
\[
\operatorname{Ext}_{Y^+}^k\big(\tau^1, i^+_*\mathcal{O}_{E_Y}(-i)\big)=0,
\]
because the upper bound of the degree support of $\tau^1$ decreases as $s'$ increases. We finish the proof.
\end{proof}
If we consider natural embedding of complex induced by:
\[
\begin{tikzcd}[column sep=small, row sep=small]
 & \vdots \arrow[d] & \vdots \arrow[d] \\
 P^{r} \otimes \wedge^2 V_c \arrow[r, "F_{1}"] \arrow[d, "-K_2"] &
 P^{r} \otimes V_c \arrow[r, "F_{0}"] \arrow[d, "K_1"] &
 P^{r} \arrow[d, "-K_0"] \\
 P^{r} \otimes \wedge^3 V_c \arrow[r, "F_{2}"] \arrow[d, "-K_3"] &
 P^{r} \otimes \wedge^2 V_c \arrow[r, "F_{1}"] \arrow[d, "K_2"] &
 P^{r} \otimes V_c \arrow[r, "F_{0}"] \arrow[d, "-K_1"] &
 P^{r} \arrow[d, "K_0"] \\
 \cdots \arrow[r] &
 P^{r} \otimes \wedge^3 V_c \arrow[r, "F_{2}"] &
 P^{r} \otimes \wedge^2 V_c \arrow[r, "F_{1}"] &
 P^{r} \otimes V_c \arrow[r, "F_{0}"] &
 P^{r}
\end{tikzcd}
\]
to our filtration, it induces a natural morphism of \(\pi^*_{Y}\mathcal{O}_{Z}\) to \(F^s(\mathcal{O}_{\widetilde{Z}})\)  in  \( \mathrm{D}^-(Y^+) \) so we can extend the filtration:
\[
\begin{tikzcd}[column sep=0.5em]
\pi^*_{Y}\mathcal{O}_{Z}\arrow{rr} && F^s(\mathcal{O}_{\widetilde{Z}})\arrow{dl}\arrow{rr}&& F^{s-1}\arrow{dl}\arrow{rr}&& F^{s-2}  \arrow{dl} \\
&\tau_s\arrow[ul,dashed,"\Delta"]&& i^+_{*}\mathcal{D}^Y_{s}[1]\arrow[ul,dashed,"\Delta"]&& i^+_{*}\mathcal{D}^Y_{s-1}[1]\arrow[ul,dashed,"\Delta"]
\end{tikzcd}
\quad\cdots\quad
\begin{tikzcd}[column sep=0.5em]
 &   F^1\arrow{rr}&&  F^0\arrow{rr}\arrow{dl}&& 0\arrow{dl} \\
 && i^+_{*}\mathcal{D}^Y_{1}[1]\arrow[ul,dashed,"\Delta"]  &&i^+_{*}\mathcal{D}^Y_{0}[1]\arrow[ul,dashed,"\Delta"]
\end{tikzcd}
\]
We observe that when \( s \) is sufficiently large, the cohomology \( \mathcal{H}^i(\tau_s) \) vanishes for all \( i > -s + N \), where \( N \) is a fixed integer determined by the parameters \( c \) and \( d \). Therefore, we may view
\[
F^s(\mathcal{O}_{\widetilde{Z}})
\]
as an approximation to \( \pi^*_{Y} \mathcal{O}_{Z} \). In fact, under suitable defintion we have:
\[
\pi^*_{Y} \mathcal{O}_{Z} = \varprojlim_{s \to \infty} F^s(\mathcal{O}_{\widetilde{Z}}).
\]
We have the following immediate corollary:
\begin{corollary}\label{pushfoFs}
\begin{enumerate}
\item For any integer $s$, we have:
\[
\pi_{Y*}\left[F^s(\mathcal{O}_{\widetilde{Z}})\right] = \pi_{E_{Y}*}(\mathcal{D}^Y_0) = \mathcal{O}_Z
\]

\item There exist positive integers $p,q$ determined by $(c,d)$ such that $\tau_s$ has the following periodicity:
\[
\tau_{s+p} = \tau_s [q]
\]

Moreover, we have a distinguished triangle in $\mathrm{D}^-(Y^+)$:
\[
\tau_{s+N} \to G_{s,N} \to \tau_s {\to} \tau_{s+N}[1]
\]
where $G_{s,N}$ belongs to the full saturated subcategory of $\mathrm{D}^b(Y^+)$ generated by $i^+_*\mathcal{A}_{-1}$ for any $s$.
\end{enumerate}
\end{corollary}

For any Cohen–Macaulay module $M$ over $S$, we consider its standard resolution:
\[
\begin{tikzcd}[column sep=normal]
\cdots \arrow[r] &
\mathbf{U}^r_4 \arrow[r, "\mathbf{\bar{F}}_4"] &
\mathbf{U}^r_3 \arrow[r, "\mathbf{\bar{F}}_3"] &
\mathbf{U}^r_2 \arrow[r, "\mathbf{\bar{F}}_2"] &
\mathbf{U}^r_1 \arrow[r, "\mathbf{\bar{F}}_1"] &
\mathbf{U}^r_0 \arrow[r, "\mathbf{\bar{F}}_0"] &
\mathbf{M}^r \arrow[r] &
0
\end{tikzcd}
\]
Through the right projective resolution (or by the periodicity property for hypersurfaces), we obtain the following extended complex:
\[
\begin{tikzcd}[column sep=normal]
\cdots \arrow[r] &
\mathbf{U}^r_{2} \arrow[r, "A"] &
\mathbf{U}^r_{1} \arrow[r, "B"] &
\mathbf{U}^r_{0} \arrow[r, "A"] &
\mathbf{U}^r_{-1} \arrow[r, "B"] &
\mathbf{U}^r_{-2} \arrow[r, "A"] &
\cdots
\end{tikzcd}
\]
where we assume for any integers \(k\):
\[
\mathbf{U}^r_{k} := \bigoplus_{i=1}^m P^{r-a^{(k)}_{i}},
\]
for a series of integers \(a^{(k)}_{i}\).
We construct a subcomplex whose terms are given by
\[
\mathbf{U}'^r_{k} := \bigoplus_{i=1}^m P^{r - a'^{(k)}_i},
\]
where
\[
a'^{(k)}_i :=
\begin{cases}
a^{(k)}_i & \text{if } a^{(k)}_i \geq 0, \\
0 & \text{if } a^{(k)}_i < 0.
\end{cases}
\]
and where all maps are the corresponding restrictions from the original complex, by the definition of the degrees it is evident that this is indeed a complex. Next, we consider applying the naive truncation \( \sigma^{\leq 2n} \) to the above complex for integer \( n \). We require that \( n \) is sufficiently large so that
\[
a'^{(2n)}_i = 0 \quad \text{for all } i.
\]
Then, by applying the shift \([2n]\), we turn the truncated complex into a complex resembling \( M \), and we denote it as \(\left\lfloor M \right\rfloor_{\leq r}\):
\[
\left\lfloor M \right\rfloor_{\leq r}:=\sigma^{\leq 2n} \mathbf{U}'^r_{\bullet}[2n],
\]
which we may regard as an approximation to \( M \). It has the following form:
\[
\begin{tikzcd}[column sep=normal]
\cdots \arrow[r] &
\displaystyle\bigoplus_{i=1}^m P^{r - a'^{(2n+2)}_i} \arrow[r, "A"] &
\displaystyle\bigoplus_{i=1}^m P^{r - a'^{(2n+1)}_i} \arrow[r, "B"] &
\displaystyle\bigoplus_{i=1}^m P^{r} \arrow[r, "\sim"] &
\left\lfloor M \right\rfloor_{\leq r}
\end{tikzcd}
\]
with $a'^{(k)}_i\geq0$.
Similarly, we consider the following complex:
\[
\begin{tikzcd}[column sep=normal]
\cdots \arrow[r] &
\displaystyle\bigoplus_{i=1}^m P^{r - a''^{(2n+2)}_i} \arrow[r, "A"] &
\displaystyle\bigoplus_{i=1}^m P^{r - a''^{(2n+1)}_i} \arrow[r, "B"] &
\displaystyle\bigoplus_{i=1}^m P^{r} \arrow[r, "\sim"] &
\left\lfloor M \right\rfloor^{(1)}_{\leq r}
\end{tikzcd}
\]
where
\[
a''^{(k)}_i :=
\begin{cases}
a'^{(k)}_i - 1 & \text{if } a'^{(k)}_i - 1 \geq 0, \\
0 & \text{if } a'^{(k)}_i - 1 < 0.
\end{cases}
\]
They have a natural embedding via the following commutative diagram with quotient:
\[
\begin{tikzcd}[column sep=normal]
\cdots \arrow[r] &
\displaystyle\bigoplus_{i=1}^m P^{r - a''^{(2n+2)}_i} \arrow[r, "A"] \arrow[d, hook] &
\displaystyle\bigoplus_{i=1}^m P^{r - a''^{(2n+1)}_i} \arrow[r, "B"] \arrow[d, hook] &
\displaystyle\bigoplus_{i=1}^m P^r \arrow[r, "\sim"] \arrow[d, equal] &
\left\lfloor M \right\rfloor^{(1)}_{\leq r} \arrow[d, hook] \\
\cdots \arrow[r] &
\displaystyle\bigoplus_{i=1}^m P^{r - a'^{(2n+2)}_i} \arrow[r, "A"] \arrow[d, two heads] &
\displaystyle\bigoplus_{i=1}^m P^{r - a'^{(2n+1)}_i} \arrow[r, "B"] \arrow[d, two heads] &
\displaystyle\bigoplus_{i=1}^m P^r \arrow[r, "\sim"] \arrow[d, equal] &
\left\lfloor M \right\rfloor_{\leq r} \arrow[d, two heads] \\
\cdots \arrow[r] &
\displaystyle\bigoplus_{i=1}^m  P^{r - a'^{(2n+2)}_i}/ P^{r - a''^{(2n+2)}_i} \arrow[r, "L(A)"] &
\displaystyle\bigoplus_{i=1}^m P^{r - a'^{(2n+1)}_i}/ P^{r - a''^{(2n+1)}_i}  \arrow[r] &
0\arrow[r, two heads] &
\left\lfloor M \right\rfloor_{\leq r}/\left\lfloor M \right\rfloor^{(1)}_{\leq r}
\end{tikzcd}
\]
Using our previous operation, for any \( r \geq 0 \), we take the direct sum and then project to $\mathrm{D}^b(Y^+)$. We obtain the following commutative diagram:
\[
\begin{tikzcd}[column sep=normal]
\cdots \arrow[r] &
\displaystyle\bigoplus_{i=1}^m \mathcal{O}_{Y^+}(a''^{(2n+2)}_i E) \arrow[r, "A"] \arrow[d, hook] &
\displaystyle\bigoplus_{i=1}^m \mathcal{O}_{Y^+}(a''^{(2n+1)}_i E) \arrow[r, "B"] \arrow[d, hook] &
\displaystyle\bigoplus_{i=1}^m \mathcal{O}_{Y^+} \arrow[r, "\sim"] \arrow[d, equal] &
\mathcal{F}^1(\mathcal{R}(M)) \arrow[d, hook] \\
\cdots \arrow[r] &
\displaystyle\bigoplus_{i=1}^m \mathcal{O}_{Y^+}(a'^{(2n+2)}_i E) \arrow[r, "A"] \arrow[d, two heads] &
\displaystyle\bigoplus_{i=1}^m \mathcal{O}_{Y^+}(a'^{(2n+1)}_i E) \arrow[r, "B"] \arrow[d, two heads] &
\displaystyle\bigoplus_{i=1}^m \mathcal{O}_{Y^+} \arrow[r, "\sim"] \arrow[d, equal] &
\mathcal{R}(M) \arrow[d, two heads] \\
\cdots \arrow[r] &
\displaystyle\bigoplus_{\substack{i=1\\a'^{(2n+2)}_i>0}}^m\mathcal{O}_{E_Y}(-a'^{(2n+2)}_i) \arrow[r, "L(A)"] &
\displaystyle\bigoplus_{\substack{i=1\\a'^{(2n+1)}_i>0}}^m \mathcal{O}_{E_Y}(-a'^{(2n+1)}_i) \arrow[r] &
0 \arrow[r,"\sim"] &
\operatorname{gr}^1(M)
\end{tikzcd}
\]
we readily observe that $\operatorname{gr}^1(M)$ forms a graded matrix factorization truncated at graded degree $\geq 1$, which belongs to $\mathcal{A}_{-1}$ by Lemma \ref{MFlemma}.\\

We iteratively apply above method to obtain a natural filtration \(\mathcal{F}^s(\mathcal{R}(M))\) for any positive integers $s$ as follows:
\[
\begin{tikzcd}[column sep=0.5em]
 & \mathcal{F}^s(\mathcal{R}(M)) \arrow{rr}&& \mathcal{F}^{s-1}\arrow{dl}\arrow{rr}&& \mathcal{F}^{s-2}  \arrow{dl} \\
&&  i^+_{*}\operatorname{gr}^{s}(M)\arrow[ul,dashed,"\Delta"]&&  i^+_{*}\operatorname{gr}^{s-1}\arrow[ul,dashed,"\Delta"]
\end{tikzcd}
\quad\cdots\quad
\begin{tikzcd}[column sep=0.5em]
 &   \mathcal{F}^2\arrow{rr}&&  \mathcal{F}^1\arrow{rr}\arrow{dl}&& \mathcal{R}(M)\arrow{dl} \\
 && i^+_{*}\operatorname{gr}^2\arrow[ul,dashed,"\Delta"]  &&i^+_{*}\operatorname{gr}^1(M)\arrow[ul,dashed,"\Delta"]
\end{tikzcd}
\]
for any $i$, we have $\operatorname{gr}^i(M) \in \mathcal{A}_{-1}$. Naturally, we obtain a map from $\pi_Y^*M$ to $\mathcal{F}^s(\mathcal{R}(M))$ which extends canonically to the following diagrams:
\[
\begin{tikzcd}[column sep=0.5em]
 & \pi_Y^*M \arrow{rr} && \mathcal{F}^{s} \arrow{dl} \arrow{rr} && \mathcal{F}^{s-1} \arrow{dl} \\
 && \tau^{s} \arrow[ul, dashed, "\Delta"] && i^+_{*}\operatorname{gr}^{s}(M) \arrow[ul, dashed, "\Delta"]
\end{tikzcd}
\quad\cdots\quad
\begin{tikzcd}[column sep=0.5em]
 & \mathcal{F}^2 \arrow{rr} && \mathcal{F}^1 \arrow{rr} \arrow{dl} && \mathcal{R}(M) \arrow{dl} \\
 && i^+_{*}\operatorname{gr}^2(M) \arrow[ul, dashed, "\Delta"] && i^+_{*}\operatorname{gr}^1(M) \arrow[ul, dashed, "\Delta"]
\end{tikzcd}
\]
We observe that when \( s \) is sufficiently large, the cohomology \( \mathcal{H}^i(\tau_s) \) vanishes for all \( i > -s + N \), where \( N \) is a fixed integer determined by \( M \) and our choice of standard resolution. Therefore, we may view
\[
\mathcal{F}^s(\mathcal{R}(M))
\]
as an approximation to \( \pi^*_{Y}M \). In fact, under a suitable definition, we have:
\[
\pi^*_{Y}M = \varprojlim_{s \to \infty} \mathcal{F}^s(\mathcal{R}(M)).
\]
There is an immediate corollary:
\begin{corollary}
There exist positive integers \( p \) and \( q \), determined by \( M \) and our choice of standard resolution, such that the associated graded pieces \( \mathrm{gr}^s(M) \) satisfy the following periodicity:
\[
\mathrm{gr}^{s+p}(M) \cong \mathrm{gr}^s(M)[q].
\]
for any $s$.
\end{corollary}
We would like to perform a computation analogous to \( \pi_{Y*}\pi_Y^*M \) in the derived category. However, using only elementary methods, we can define
\[
\pi_Y^* \colon \mathrm{D}^-(Y) \to \mathrm{D}^-(Y^+), \quad \text{and} \quad \pi_{Y*} \colon \mathrm{D}^+(Y^+) \to \mathrm{D}^+(Y),
\]
so the composition \( \pi_{Y*}\pi_Y^*M \) requires more advanced techniques. To avoid certain technical difficulties, we instead consider the Verdier pullback:
\[
\pi_Y^! \colon \mathrm{D}^+(Y) \to \mathrm{D}^+(Y^+),
\]
and aim to compute \( \pi_Y^!M \).
The functor \( \pi_Y^! \) is defined via the composition of a closed immersion followed by a smooth morphism. Assuming that we are given the dualizing divisors \( K_Y \) and \( K_{Y^+} \) for \( Y \) and \( Y^+ \), respectively, we have:
\[
\pi_Y^! M = \mathcal{S}_{Y^+} \circ \pi_Y^* \circ \mathcal{S}^{-1}_Y(M),
\]
where \( \mathcal{S}_{Y^+} \) and \( \mathcal{S}_Y \) denote the dualizing functor associated to \( Y^+ \) and \( Y \), respectively. Explicitly, these are given by:
\[
\mathcal{S}_Y(M) := \mathbf{R}\!\operatorname{\mathcal{H}om}_{\mathcal{O}_Y}(M, \mathcal{O}_Y(K_Y)[n]), \qquad
\mathcal{S}_{Y^+}(N) := \mathbf{R}\!\operatorname{\mathcal{H}om}_{\mathcal{O}_{Y^+}}(N, \mathcal{O}_{Y^+}(K_{Y^+})[n]).
\]
Since \( Y \) and \( Y^+ \) are birational, their shifts cancel out, we will therefore omit shifts in what follows.\\

Recall \( M \) is a Cohen–Macaulay module, its Serre dual
\[
M^T := \mathcal{S}^{-1}_Y(M)= \mathcal{S}_Y(M)
\]
is also a CM module over \( Y \). Then under the action of \( \pi_Y^* \), and using the computations developed earlier, we obtain the following filtration diagram:
\[
\begin{tikzcd}[column sep=0.5em]
 & \pi_Y^*M^T \arrow{rr} && \mathcal{F}^{s} \arrow{dl} \arrow{rr} && \mathcal{F}^{s-1} \arrow{dl} \\
 && \tau^{s} \arrow[ul, dashed, "\Delta"] && i^+_{*}\operatorname{gr}^{s}(M^T) \arrow[ul, dashed, "\Delta"]
\end{tikzcd}
\quad\cdots\quad
\begin{tikzcd}[column sep=0.5em]
 & \mathcal{F}^2 \arrow{rr} && \mathcal{F}^1 \arrow{rr} \arrow{dl} && \mathcal{R}(M^T) \arrow{dl} \\
 && i^+_{*}\operatorname{gr}^2(M^T) \arrow[ul, dashed, "\Delta"] && i^+_{*}\operatorname{gr}^1(M^T) \arrow[ul, dashed, "\Delta"]
\end{tikzcd}
\]
Applying the functor $\mathcal{S}_{Y^+}$ to the above filtration, we obtain the following mirror-symmetric diagram:
\[
\begin{tikzcd}[column sep=0.1em]
 & \mathcal{S}_{Y^+}(\mathcal{R}(M^T)) \arrow{rr} &&\mathcal{S}_{Y^+}(\mathcal{F}^{1}) \arrow{dl} \\
 && \mathcal{S}_{Y^+}(i^+_{*}\operatorname{gr}^1(M^T))[1] \arrow[ul, dashed, "\Delta"]
\end{tikzcd}
\cdots
\begin{tikzcd}[column sep=0.001em]
 & \mathcal{S}_{Y^+}(\mathcal{F}^{s-1}) \arrow{rr} && \mathcal{S}_{Y^+}(\mathcal{F}^s) \arrow{rr} \arrow{dl} && \pi_Y^!M \arrow{dl} \\
 && \mathcal{S}_{Y^+}(i^+_{*}\operatorname{gr}^{s}(M^T))[1]  \arrow[ul, dashed, "\Delta"] &&  \mathcal{S}_{Y^+}(\tau^{s})[1] \arrow[ul, dashed, "\Delta"]
\end{tikzcd}
\]
Noticing we have:
\[
\begin{aligned}
\mathcal{S}_{Y^+}(i^+_*N)
&:= \mathbf{R}\!\operatorname{\mathcal{H}om}_{\mathcal{O}_{Y^+}}\left(i^+_*N, \mathcal{O}_{Y^+}(K_{Y^+})\right)[1] \\
&= i^+_*\mathbf{R}\!\operatorname{\mathcal{H}om}_{\mathcal{O}_{E_Y}}\left(N, (i^+)^!\mathcal{O}_{Y^+}(K_{Y^+})\right)[1] \\
&= i^+_*\mathbf{R}\!\operatorname{\mathcal{H}om}_{\mathcal{O}_{E_Y}}\left(N, \mathcal{O}_{E_Y}(K_{E_Y})\right) \\
&= i^+_*\mathbf{R}\!\operatorname{\mathcal{H}om}_{\mathcal{O}_{E_Y}}\left(N, \mathcal{O}_{E_Y}(-(c-d))\right)
\end{aligned}
\]
In particular, if \( N \in \mathcal{A}_{-1} \), then by Lemma \ref{MFlemma}:
\[
\mathbf{R}\!\operatorname{\mathcal{H}om}_{\mathcal{O}_{E_Y}}(N, \mathcal{O}_{E_Y}) \in \mathcal{A}_{0},
\]
so
\[
\mathcal{S}_{Y^+}(i^+_*N) \in \mathcal{A}_{-(c-d)}.
\]
Now consider the following semiorthogonal decomposition of \( \mathrm{D}^b(E_Y) \):
\[
\mathrm{D}^b(E_Y) =
\left\langle
\mathcal{A}_{-(c-d)},\,
\mathcal{O}_{E_Y}(-(c-d)+1),\,
\dots,\,
\mathcal{O}_{E_Y}(-1),\,
\mathcal{O}_{E_Y}
\right\rangle.
\]
It follows that:
\[
\pi_{Y*}\left(\mathcal{S}_{Y^+}(i^+_*N)\right) = 0.
\]
On the other hand, we observe that when \( s \) is sufficiently large, the cohomology
\[
\mathcal{H}^i\big(\mathcal{S}_{Y^+}(\tau_s)\big)
\]
vanishes for all \( i < s + N \), where \( N \) is an integer determined by the choice of standard resolution of \( M^T \).

\begin{corollary}
Let \( M \) be a Cohen–Macaulay module over \( Y \). Given a choice of standard resolution, we can construct an object
\[
\mathcal{R}^!(M) := \mathcal{S}_{Y^+}(\mathcal{R}(M^T)),
\]
which by construction belongs to $\mathrm{D}^b(Y^+)$, such that
\[
\pi_{Y*}(\mathcal{R}^!(M)) \cong M.
\]
\end{corollary}

\begin{proof}
Consider the distinguished triangle in \( \mathrm{D}^+(Y^+) \):
\[
\mathcal{R}^!(M) \longrightarrow \pi_Y^!M \longrightarrow \tau \longrightarrow \mathcal{R}^!(M)[1].
\]
Applying \( \pi_{Y*} \), we obtain a triangle in \( \mathrm{D}^+(Y) \):
\[
\pi_{Y*} \mathcal{R}^!(M) \longrightarrow \pi_{Y*} \pi_Y^! M \longrightarrow \pi_{Y*} \tau \longrightarrow \pi_{Y*} \mathcal{R}^!(M)[1].
\]
By Grothendieck–Verdier duality and the fact that \( \pi_{Y*} \mathcal{O}_{Y^+} = \mathcal{O}_Y \) by Lemma \ref{ratlemma}, we have
\[
\pi_{Y*} \pi_Y^! M \cong M.
\]
To conclude the proof, it remains to show \( \pi_{Y*} \tau \simeq 0 \). This is equivalent to showing that
\[
\mathcal{H}^i(\pi_{Y*} \tau) = 0 \quad \text{for all integers } i.
\]
Fix any integer \( i \), and choose a sufficiently large \( s \gg i \). Then we consider the truncated triangle
\[
\mathcal{R}^!(M) \longrightarrow \mathcal{S}_{Y^+}(\mathcal{F}^s) \longrightarrow \tau^s \longrightarrow \mathcal{R}^!(M)[1],
\]
which implies
\[
\mathcal{H}^i(\pi_{Y*} \tau) = \mathcal{H}^i(\pi_{Y*} \tau^s).
\]
By the previous computation, we have
\[
\pi_{Y*} \tau^s = 0.
\]
Hence \( \mathcal{H}^i(\pi_{Y*} \tau) = 0 \) for this \( i \), and the proof is complete.
\end{proof}

\begin{corollary}\label{pushforwardlemma}
Let \( M \) be a Cohen–Macaulay module over \( Y \). Given a choice of standard resolution, we can construct an object
\[
\widetilde{\mathcal{R}}^*(M) := \mathcal{R}^!(M)(-K_{Y^+/Y})=\mathcal{S}_{Y^+}(\mathcal{R}(M^T)) \otimes_{\mathcal{O}_{Y^+}} \mathcal{O}(-K_{Y^+/Y}),
\]
which by construction belongs to $\mathrm{D}^b(Y^+)$, such that
\[
\pi_{Y*}(\widetilde{\mathcal{R}}^*(M)) \cong M.
\]
\end{corollary}

\begin{proof}
By Grothendieck–Verdier duality, we have:
\[
\pi_* \pi^!M \otimes \mathcal{O}_{Y^+}(-K_{Y^+/Y})
\simeq \pi_* \mathbf{R}\!\operatorname{\mathcal{H}om}_{\mathcal{O}_{Y^+}}\left( \mathcal{O}_{Y^+}(K_{Y^+/Y}), \pi^!M \right)
\simeq \mathbf{R}\!\operatorname{\mathcal{H}om}_{\mathcal{O}_Y}\left( \pi_*\mathcal{O}_{Y^+}(K_{Y^+/Y}), M \right),
\]
We know that the relative canonical divisor is
\[
K_{Y^+/Y} \sim (c - d - 1) E_Y.
\]
Under the assumption that \( c - d - 1 \geq 0 \), and from cohomology computations on hypersurfaces, we have:
\[
\operatorname{Ext}^*(\mathcal{O}_{E_Y}, \mathcal{O}_{E_Y}(i)) = 0 \quad \text{for} \quad -(c - d) + 1 \leq i \leq -1,
\]
and
\[
\operatorname{Ext}^*(\mathcal{O}_{E_Y}, \mathcal{O}_{E_Y}) = \mathcal{O}_Z.
\]
Now, consider the short exact sequence on \( Y^+ \) for \( j \geq 1 \):
\[
0 \longrightarrow \mathcal{O}_{Y^+}((j - 1)E_Y) \longrightarrow \mathcal{O}_{Y^+}(jE_Y) \longrightarrow \mathcal{O}_{E_Y}(-j) \longrightarrow 0.
\]
and \( \pi_{Y*} \mathcal{O}_{Y^+} \cong \mathcal{O}_Y \) by Lemma \ref{ratlemma}, it follows by induction that:
\begin{equation}\label{eq:pushforward-relative-dualizing-Y}
\pi_* \mathcal{O}_{Y^+}(K_{Y^+/Y}) \cong \mathcal{O}_Y.
\end{equation}
Therefore, by Grothendieck–Verdier duality, we conclude:
\[
\pi_* \left( \pi^! M \otimes \mathcal{O}_{Y^+}(-K_{Y^+/Y}) \right) \cong M.
\]
On the other hand, for any object \( N \), if
\[
\mathcal{S}_{Y^+}(i^+_*N) \in \mathcal{A}_{-(c-d)},
\]
then since
\[
\mathcal{O}_{Y^+}(K_{Y^+/Y})|_{E_Y} \simeq \mathcal{O}_{E_Y}(-(c-d)+1),
\]
we have
\[
\mathcal{S}_{Y^+}(i^+_*N) \otimes_{\mathcal{O}_{Y^+}} \mathcal{O}_{Y^+}(-K_{Y^+/Y}) \in \mathcal{A}_{-1}.
\]
In particular, this implies
\[
\pi_{Y*} \left( \mathcal{S}_{Y^+}(i^+_*N) \otimes_{\mathcal{O}_{Y^+}}  \mathcal{O}_{Y^+}(-K_{Y^+/Y}) \right) = 0.
\]
Putting everything together and applying exactly the same argument as before, we obtain the desired conclusion.
\end{proof}

\begin{corollary} \label{constructionlemma-Y}
Let \( F \) be any object in \( \mathrm{D}^b(Y) \). Given a choice of standard resolution, we can construct objects
\[
\mathcal{R}\pi^!_Y(F) \quad \text{and} \quad \widetilde{\mathcal{R}}\pi^*_Y(F)
\]
which belong to \( \mathrm{D}^b(Y^+) \), such that
\[
\pi_{Y*}\mathcal{R}\pi^!_Y(F) \cong\pi_{Y*}\widetilde{\mathcal{R}}\pi^*_Y(F) \cong F.
\]
\end{corollary}

\begin{proof}
For any object \( F \in \mathrm{D}^{b}(Y) \), there exists a distinguished triangle by Lemma \ref{CMapprox}:
\[
P_{F} \longrightarrow F \longrightarrow M_{F}[n_{F}] \in \Delta,
\]
where \( M_{F} \) is a Cohen–Macaulay \( Y\)-module and \( P_{F} \in \mathrm{D}^{\mathrm{Perf}}(Y) \). Then we have the following distinguished triangle in \( \mathrm{D}^b(Y^+) \):
\[
\pi_Y^! P_{F} \longrightarrow \mathcal{R}\pi^!_Y(F)\longrightarrow \mathcal{R}^!(M_{F})[n_{F}]\in \Delta,
\]
by considering the morphism space:
\[
\operatorname{Hom}_{\mathrm{D}^b(Y^+)}(\mathcal{R}^!(M_{F})[n_{F}], \pi_Y^! P_{F}[1])
\simeq \operatorname{Hom}_{\mathrm{D}^b(Y)}(\pi_{Y*}\mathcal{R}^!(M_{F})[n_{F}], P_{F}[1])
\simeq \operatorname{Hom}_{\mathrm{D}^b(Y)}(M_{F}[n_{F}], P_{F}[1]).
\]
Therefore, the morphism in the first triangle uniquely determines the morphism in the second triangle by adjunction.\\

We define the mapping cone of this morphism to be \( \widetilde{\mathcal{R}}\pi^!_Y(F) \). Applying the pushforward functor \( \pi_{Y*} \) to the second triangle, and using the property \( \pi_{Y*} \widetilde{\mathcal{R}}^!(M_F) \simeq M_F \), and the projection formula for \( \pi_Y^*P_F \), we obtain:
\[
\pi_{Y*} \widetilde{\mathcal{R}}\pi^!_Y(F) \simeq F.
\]
A similar argument holds for \(\widetilde{\mathcal{R}}\pi^* _Y(F)\) which completes the construction and proof.
\end{proof}

Recall that we would like to perform a computation analogous to \( \pi_{Y*} \pi_Y^* M \) in the derived category. In complete analogy with the above, we define:
\[
\pi_{Y!} \mathcal{F} := \mathcal{S}_{Y} \circ \pi_{Y*} \circ \mathcal{S}_{Y^+}(\mathcal{F}),
\]
which gives a well-defined functor
\[
\pi_{Y!} \colon \mathrm{D}^-(Y) \to \mathrm{D}^-(Y^+).
\]
According to \eqref{eq:pushforward-relative-dualizing-Y} and noting that duality behaves well under pushforward, we have
\[
\pi_{Y*}\mathbf{R}\!\operatorname{\mathcal{H}om}_{\mathcal{O}_{Y^+}}(\pi_Y^*M, \mathcal{O}_{Y^+}(K_{Y^+}))
= \mathbf{R}\!\operatorname{\mathcal{H}om}_{\mathcal{O}_Y}(M, \pi_{Y*} \mathcal{O}_{Y^+}(K_{Y^+}))
= \mathbf{R}\!\operatorname{\mathcal{H}om}_{\mathcal{O}_Y}(M, \mathcal{O}_Y(K_Y)).
\]
Therefore, we obtain
\[
\pi_{Y!} \pi_Y^* M = \mathcal{S}_Y \circ \mathcal{S}_Y (M) = M.
\]

By a similar argument as above, we have the following:

\begin{corollary}
Let \( M \) be a Cohen–Macaulay module over \( Y \). Given a choice of standard resolution, we can define a functorial lift of \( M \) to \( \mathrm{D}^b(Y^+) \) via the previously constructed complex
\[
\mathcal{R}(M)
\]
which belongs to \( \mathrm{D}^b(Y^+) \), and satisfies
\[
\pi_{Y!}(\mathcal{R}(M)) \cong M.
\]
\end{corollary}

Let us observe the behavior of \( \pi_{Y!}(\mathcal{R}(M)) \). By applying Grothendieck--Verdier duality, one obtains:
\[
\begin{aligned}
\pi_{Y*}\mathbf{R}\!\operatorname{\mathcal{H}om}_{\mathcal{O}_{Y^+}}(\mathcal{R}(M), \mathcal{O}_{Y^+}(K_{Y^+}))
&= \pi_{Y*}\mathbf{R}\!\operatorname{\mathcal{H}om}_{\mathcal{O}_{Y^+}}(\mathcal{R}(M), \pi_Y^! \mathcal{O}_Y(K_Y)) \\
&\simeq \mathbf{R}\!\operatorname{\mathcal{H}om}_{\mathcal{O}_Y}(\pi_{Y*} \mathcal{R}(M), \mathcal{O}_Y(K_Y)).
\end{aligned}
\]
This leads to the identity:
\[
\pi_{Y!}(\mathcal{R}(M)) \simeq \mathcal{S}_Y \circ \mathcal{S}_Y(\pi_{Y*}\mathcal{R}(M)) \simeq \pi_{Y*}\mathcal{R}(M).
\]
As a consequence, we arrive at the following corollary.

\begin{corollary}
Let \( M \) be a Cohen–Macaulay module over \( Y \). Given a choice of standard resolution, we can define a functorial lift of \( M \) to \( \mathrm{D}^b(Y^+) \) via the previously constructed complex
\[
\mathcal{R}(M)
\]
which belongs to \( \mathrm{D}^b(Y^+) \), and satisfies
\[
\pi_{Y*}(\mathcal{R}(M)) \cong M.
\]
\end{corollary}

\begin{corollary} \label{constructionlemma}
Let \( F \) be any object in \( \mathrm{D}^b(Y) \). Given a choice of standard resolution, we can construct objects
\[
\mathcal{R}\pi^*_Y(F) \quad \text{and} \quad \widetilde{\mathcal{R}}\pi^!_Y(F):=\mathcal{R}\pi^*_Y(F)(K_{Y^+/Y})
\]
which belong to \( \mathrm{D}^b(Y^+) \), such that
\[
\pi_{Y*}\mathcal{R}\pi^*_Y(F) \cong\pi_{Y*}\widetilde{\mathcal{R}}\pi^!_Y(F) \cong F.
\]
\end{corollary}

\subsection{Local adjunction}
We continue with the notation and assumptions introduced in the previous subsection and consider the following cohomology computation.

\begin{lemma}\label{ortho}
Let \( i^+ : E_Y \hookrightarrow Y^+ \) denote the natural inclusion. Then under the assumption \( c - d > 1 \), we have:
\[
\operatorname{Ext}^*_{\mathrm{D}^b(Y^+)}\left(i^+_*\mathcal{O}_{E_Y}(-i),\,i^+_*\mathcal{O}_{E_Y}(-j)\right) \cong \mathcal{O}_Z\cdot\delta_{i,j}.
\]
for any \(1\leq i\leq j\leq c-d-1\). In particular, the object \( i^+_*\mathcal{O}_{E_Y}(-i) \in \mathrm{D}^b(Y^+) \) is an exceptional (pre-tilting) object in \( \mathrm{D}^b(Y^+) \).
\end{lemma}

\begin{proof}
Consider the following excess distinguished triangle:
\[
\mathcal{O}_{E_Y}(-i)(-E_Y)[1] \longrightarrow i^{+*}i^+_*\mathcal{O}_{E_Y}(-i) \longrightarrow \mathcal{O}_{E_Y}(-i).
\]
By assumption, \( E_Y \subset E \) is a hypersurface with Gorenstein parameter \( a > 1 \). It follows that
\[
\operatorname{Ext}^*\left(\mathcal{O}_{E_Y}(-i)(-E_Y),\, \mathcal{O}_{E_Y}(-j)\right) = 0,
\quad \text{and} \quad
\operatorname{Ext}^*\left(\mathcal{O}_{E_Y}(-i),\, \mathcal{O}_{E_Y}(-j)\right) = \mathcal{O}_Z\cdot\delta_{i,j}.
\]
for any \(i\) and \(j\) in the range. Applying the long exact sequence in cohomology associated to the distinguished triangle yields the result.
\end{proof}

\begin{corollary}\label{sodofker}
In \(\mathrm{D}^b(Y^+)\) there exists a full saturated subcategory generated by the following semi-orthogonal decomposition:
\[
\left\langle
i^+_*\mathcal{O}_{E_Y}(-(c-d)+1),\, \ldots,\, i^+_*\mathcal{O}_{E_Y}(-1),\, \left\langle i^+_*\mathcal{A}_{-1} \right\rangle
\right\rangle,
\]
where \(\left\langle i^+_*\mathcal{A}_{-1} \right\rangle\) denotes the full saturated subcategory of \(\mathrm{D}^b(Y^+)\) generated by objects of the form \(i^*_+(N)\) for \(N \in \mathcal{A}_{-1}\).
\end{corollary}

\begin{definition}
we define
\[
[\mathrm{D}^b(Y^+)] := \mathrm{D}^b(Y^+) \big/ \langle i^+_*\mathcal{A}_{-1} \rangle
\]
to be the Verdier quotient of \(\mathrm{D}^b(Y^+)\) along the full saturated subcategory \(\langle i^+_*\mathcal{A}_{-1} \rangle\).
\end{definition}

\begin{lemma}\label{semiorthofserre}
The collection
\[
\left\langle
i^+_*\mathcal{O}_{E_Y}(-(c-d)+1),\, \ldots,\, i^+_*\mathcal{O}_{E_Y}(-1)
\right\rangle
\]
forms a sequence of exceptional (pre-tilting) objects in $[\mathrm{D}^b(Y^+)]$ that constitutes a semi-orthogonal decomposition.
\end{lemma}

\begin{proof}
We only need to consider $\operatorname{Ext}^*_{[\mathrm{D}^b(Y^+)]}\left(i^+_*\mathcal{O}_{E_Y}(-i),\,i^+_*\mathcal{O}_{E_Y}(-j)\right) \cong \mathcal{O}_Z\cdot\delta_{i,j}$ for any $j \leq i$. Consider the roof diagram:
\[
\begin{tikzcd}[column sep=small, row sep=small]
 & F^+ \arrow[dl, "t"'] \arrow[dr, "f"] & \\
i^+_*\mathcal{O}_{E_Y}(-i) & & i^+_*\mathcal{O}_{E_Y}(-j),
\end{tikzcd}
\]
where $\operatorname{Cone}(t) \in \langle i^+_*\mathcal{A}_{-1} \rangle$. By Lemma \ref{ortha-1o}, we have
\[
\operatorname{Ext}^*_{[\mathrm{D}^b(Y^+)]}\left(i^+_*\mathcal{O}_{E_Y}(-i),\,i^+_*\mathcal{O}_{E_Y}(-j)\right) \cong \operatorname{Ext}^*_{\mathrm{D}^b(Y^+)}\left(i^+_*\mathcal{O}_{E_Y}(-i),\,i^+_*\mathcal{O}_{E_Y}(-j)\right).
\]
According to Lemma \ref{ortho}, the latter equals $\mathcal{O}_Z\cdot\delta_{i,j}$.
\end{proof}

\begin{lemma}\label{keylemmakernal}
Assume \( c - d > 0 \) and \(s\geq (c-d)-1\). Then for any \( G \in \mathrm{D}^b(Y^+) \) there is a natural isomorphism
\[
\operatorname{Ext}^*_{[\mathrm{D}^b(Y^+)]}\left(F^s(\mathcal{O}_{\widetilde{Z}}), G\right) \cong \operatorname{Ext}^*_{\mathrm{D}^b(Y)}\left(\mathcal{O}_Z,\, \pi_{Y*}G\right),
\]
and this isomorphism is induced by the derived pushforward functor \( \pi_{Y*} \colon \mathrm{D}^b(Y^+) \to \mathrm{D}^b(Y) \).
\end{lemma}

\begin{proof}
We first consider an arbitrary roof:
\[
\begin{tikzcd}[column sep=small, row sep=small]
 & F^+ \arrow[dl, "t"'] \arrow[dr, "f"] & \\
F^s & & G
\end{tikzcd}
\]
such that \( \operatorname{Cone}(t) \in\langle i^+_*\mathcal{A}_{-1} \rangle\). Then applying the derived pushforward \( \pi_{Y*} \) to \( f \), we get a morphism
\[
\pi_{Y*}(f): \pi_{Y*}(F^+) \to \pi_{Y*}(G).
\]
Since \( \pi_{Y*}(F^+) \cong \pi_{Y*}(F^s(\mathcal{O}_{\widetilde{Z}})) \cong \mathcal{O}_{\widetilde{Z}} \) by Corollary \ref{pushfoFs}, this gives a morphism from \( \mathcal{O}_{\widetilde{Z}} \) to \( \pi_{Y*}(G) \). Then we proceed to show that it induces an isomorphism by an argument analogous to that in \cite[Proof of Theorem 2.14]{BondalKapranovSchechtman2018}.\\

The surjectivity of \(\pi_{Y*}\) is straightforward. Note that for any morphism
\[
g: \mathcal{O}_Z \longrightarrow \pi_{Y*}G \quad \text{in} \quad \mathrm{D}^b(Y),
\]
by adjunction, there exists a morphism
\[
f := \pi_Y^*\mathcal{O}_Z \longrightarrow G \quad \text{in} \quad \mathrm{D}^-(Y^+)
\]
such that \(\pi_{Y*}(f) = g\).
On the other hand, consider the filtration of \(\pi_Y^*\mathcal{O}_Z\) constructed previously:
\[
\begin{tikzcd}[column sep=0.5em]
 & \pi^*_{Y}\mathcal{O}_{Z}\arrow{rr}&& F^{s'} \arrow{rr} \arrow{dl} && F^s(\mathcal{O}_{\widetilde{Z}}) \arrow{dl} \\
 && \tau^1 \arrow[ul,dashed,"\Delta"] && \tau^2 \arrow[ul,dashed,"\Delta"]
\end{tikzcd}
\]
If we choose \(s' \gg s\), then the cohomological degrees supporting \(\tau^1\) are sufficiently low, so we have an isomorphism
\[
\operatorname{Hom}_{\mathrm{D}^-(Y^+)}(\pi_Y^*\mathcal{O}_Z, G) \cong \operatorname{Hom}_{\mathrm{D}^-(Y^+)}(F^{s'}(\mathcal{O}_{\widetilde{Z}}), G).
\]
Under this isomorphism, \(f\) induces a bounded morphism, which we also denote by \(f\). Therefore, we obtain a roof
\[
\begin{tikzcd}[column sep=small, row sep=small]
 & F^{s'} \arrow[dl, "t"'] \arrow[dr, "f"] & \\
F^s & & G
\end{tikzcd}
\]
such that \(\mathrm{Cone}(t) = \tau^2\) lies in the subcategory \(\langle i^+_*\mathcal{A}_{-1} \rangle\), and
\[
\pi_{Y*}(f \circ t^{-1}) = g.
\]
For the injectivity of \(\pi_{Y*}\), suppose we are given a roof
\[
\begin{tikzcd}[column sep=small, row sep=small]
 & F^+ \arrow[dl, "t"'] \arrow[dr, "f"] & \\
F^s & & G
\end{tikzcd}
\]
such that \( \pi_{Y*}(f) = 0 \). By adjunction, this implies that the induced morphism
\[
f: \pi_Y^* \pi_{Y*}F^+ \longrightarrow G
\]
in \( \mathrm{D}^-(Y^+) \) is the zero map.
Noting that \( \pi_Y^* \pi_{Y*}F^+ \cong \pi_Y^* \mathcal{O}_Z \), we consider the following filtration:
\[
\begin{tikzcd}[column sep=0.5em]
 & \pi_Y^* \pi_{Y*}F^+ \arrow{rr} && F^{s'} \arrow{rr} \arrow{dl} && F^s(\mathcal{O}_{\widetilde{Z}}) \arrow{dl} \\
 && \tau^1 \arrow[ul,dashed,"\Delta"] && \tau^2 \arrow[ul,dashed,"\Delta"]
\end{tikzcd}
\]
If we choose \( s' \gg s \), then the cohomological degrees supporting \( \tau^1 \) are sufficiently low. Hence we obtain the isomorphism
\[
\operatorname{Hom}_{\mathrm{D}^-(Y^+)}\left( \pi_Y^* \pi_{Y*}F^+, G \right)
\cong \operatorname{Hom}_{\mathrm{D}^-(Y^+)}\left( F^{s'}(\mathcal{O}_{\widetilde{Z}}), G \right),
\]
under which the morphism \( f \) corresponds to a morphism from \( F^{s'}(\mathcal{O}_{\widetilde{Z}}) \) to \( G \), which is necessarily zero.
On the other hand, consider the unit morphism in \( \mathrm{D}^-(Y^+) \) induced by adjunction:
\[
f: \pi_Y^* \pi_{Y*}F^+ \longrightarrow F^+.
\]
Again, using the same type of filtration, we choose \( s' \gg s \) so that
\[
\operatorname{Hom}_{\mathrm{D}^-(Y^+)}\left( \pi_Y^* \pi_{Y*}F^+, F^+ \right)
\cong \operatorname{Hom}_{\mathrm{D}^-(Y^+)}\left( F^{s'}(\mathcal{O}_{\widetilde{Z}}), F^+ \right),
\]
and this isomorphism induces a morphism from \( F^{s'}(\mathcal{O}_{\widetilde{Z}}) \) to \( F^+ \).
Combining these, we have a morphism from
\[
\begin{tikzcd}[column sep=0.5em, row sep=small]
 & F^{s'} \arrow[dl, "t'"'] \arrow[dr, "0"] & \\
F^s & & G
\end{tikzcd}
\qquad\text{to}\qquad
\begin{tikzcd}[column sep=0.5em, row sep=small]
 & F^+ \arrow[dl, "t"'] \arrow[dr, "f"] & \\
F^s & & G
\end{tikzcd}
\]
In fact, according to the definition of equivalent morphisms in the Verdier quotient, this shows that $f\circ t^{-1} = 0$ in $[\mathrm{D}^b(Y^+)]$.
\end{proof}

\begin{lemma}\label{lemmasupport}
Let \( \mathrm{D}^b_{E_Y}(Y^+) \) denote the full triangulated subcategory of \( \mathrm{D}^b(Y^+) \) consisting of objects supported on \( E_Y \). Define the Verdier quotient category:
\[
\left[ \mathrm{D}^b_{E_Y}(Y^+) \right] := \mathrm{D}^b_{E_Y}(Y^+) \big/ \left\langle i^+_*\mathcal{A}_{-1} \right\rangle.
\]
Then this quotient category admits a semi-orthogonal decomposition:
\[
\left[ \mathrm{D}^b_{E_Y}(Y^+) \right] = \left\langle
i^+_*\mathcal{O}_{E_Y}(-(c-d)+1),\, \ldots,\, i^+_*\mathcal{O}_{E_Y}(-1),\left\langle F^s\left( \mathcal{O}_{\widetilde{Z}} \right) \right\rangle
\right\rangle.
\]
where \(s=:c-d-1\) and \( \left\langle F^s\left( \mathcal{O}_{\widetilde{Z}} \right) \right\rangle \) denotes the  full triangulated subcategory of \( \left[ \mathrm{D}^b_{E_Y}(Y^+) \right] \) generated by the image of \( F^s\left( \mathcal{O}_{\widetilde{Z}} \right) \) under the natural projection functor.
\end{lemma}

\begin{proof}
By Lemmas~\ref{ortha-1o}, \ref{localpushspan} and \ref{semiorthofserre}, the orthogonality of the above semi-orthogonal decomposition is established. It remains to show that it is generating.\\

Let \( G \in \mathrm{D}^b_{E_Y}(Y^+) \). Since \( G \) belongs to the bounded derived category with support on \( E_Y \), its cohomology sheaves \( \mathcal{H}^i(G) \) are coherent sheaves supported on \( E_Y \), and hence \( G \) lies in the extension closure of these cohomology sheaves.
On the other hand, by the assumption \( c-d > 0 \) we know that the sheaves
\[
i^+_*\mathcal{O}_{E_Y}(-(c-d)+1),\, \ldots,\, i^+_*\mathcal{O}_{E_Y}(-1),\, i^+_*\mathcal{O}_{E_Y}
\]
form a collection of exceptional (or pre-tilting) objects which generate \( \mathrm{D}^b(E_Y) \). Hence, \( \mathcal{H}^i(G) \) lies in the triangulated subcategory generated by them.
Next, recall the filtration of the object \( F^s(\mathcal{O}_{\widetilde{Z}}) \), which fits into a sequence of distinguished triangles:
\[
\begin{tikzcd}[column sep=0.5em]
 & F^s(\mathcal{O}_{\widetilde{Z}}) \arrow{rr} && F^{s-1} \arrow{dl} \arrow{rr} && F^{s-2} \arrow{dl} \\
&& i^+_*\mathcal{D}^Y_s[1] \arrow[ul, dashed, "\Delta"] && i^+_*\mathcal{D}^Y_{s-1}[1] \arrow[ul, dashed, "\Delta"]
\end{tikzcd}
\quad\cdots\quad
\begin{tikzcd}[column sep=0.5em]
 & F^1 \arrow{rr} && F^0 \arrow{rr} \arrow{dl} && 0 \arrow{dl} \\
&& i^+_*\mathcal{D}^Y_1[1] \arrow[ul, dashed, "\Delta"] && i^+_*\mathcal{D}^Y_0[1] \arrow[ul, dashed, "\Delta"]
\end{tikzcd}
\]
Moreover, up to the Verdier quotient by \( \left\langle i^+_*\mathcal{A}_{-1} \right\rangle \), each \( \mathcal{D}^Y_i \in \mathrm{Coh}(E_Y) \) (for \( 1 \leq i \leq c-d-1 \)) belongs to the triangulated subcategory generated by
\[
\left\langle
i^+_*\mathcal{O}_{E_Y}(-(c-d)+1),\, \ldots,\, i^+_*\mathcal{O}_{E_Y}(-1)
\right\rangle.
\]
Therefore, we conclude that \( G \) lies in the triangulated subcategory
\[
\left\langle
i^+_*\mathcal{O}_{E_Y}(-(c-d)+1),\, \ldots,\, i^+_*\mathcal{O}_{E_Y}(-1),\, \left\langle F^s\left( \mathcal{O}_{\widetilde{Z}} \right) \right\rangle
\right\rangle
\]
modulo \( \left\langle i^+_*\mathcal{A}_{-1} \right\rangle \) since \(i^+_*\mathcal{D}^Y_0\simeq i^+_*\mathcal{O}_{E_Y}\), which completes the proof.
\end{proof}

\begin{prop}\label{proppushG0}
Let \( G \in \mathrm{D}^b(Y^+) \) be an object such that its derived pushforward satisfies \( \pi_{Y*} G = 0 \). Then \( G \) belongs to the triangulated subcategory with the following semi-orthogonal decomposition:
\[
\left\langle
i^+_*\mathcal{O}_{E_Y}(-(c-d)+1),\, \ldots,\, i^+_*\mathcal{O}_{E_Y}(-1),\, \left\langle i^+_*\mathcal{A}_{-1} \right\rangle
\right\rangle.
\]
which is equivalent to:
\[
\left\langle
\left\langle i^+_*\mathcal{A}_{-(c-d)} \right\rangle,\,
i^+_*\mathcal{O}_{E_Y}(-(c-d)+1),\, \ldots,\, i^+_*\mathcal{O}_{E_Y}(-1)
\right\rangle.
\]
\end{prop}

\begin{proof}
Consider the admissible subcategory
\[
\left\langle
i^+_*\mathcal{O}_{E_Y}(-(c-d)+1),\, \ldots,\, i^+_*\mathcal{O}_{E_Y}(-1)
\right\rangle.
\]
Since the derived pushforward \( \pi_{Y*} \) vanishes on this subcategory, we may replace \( G \) with its right mutation
\[
\mathbb{R}_{\left\langle
i^+_*\mathcal{O}_{E_Y}(-(c-d)+1),\, \ldots,\, i^+_*\mathcal{O}_{E_Y}(-1)
\right\rangle} G,
\]
without loss of generality.
Now assume \( \pi_{Y*} G = 0 \). Due to the birational nature of the morphism \( \pi_Y \colon Y^+ \to Y \), we have an isomorphism \( Y^+ \setminus E_Y \simeq Y \setminus Z \), so the support of \( G \) must be contained in \( E_Y \).
By Lemma~\ref{lemmasupport}, it follows that \( G \) is generated, up to a quotient of \( \left\langle i^+_*\mathcal{A}_{-1} \right\rangle \), by iterated extensions of \( F^s\left( \mathcal{O}_{\widetilde{Z}} \right) \).
On the other hand, applying Lemma~\ref{keylemmakernal}, under the assumption \( \pi_{Y*} G = 0 \), we obtain
\[
\operatorname{Ext}^*_{[\mathrm{D}^b(Y^+)]}\left(F^s(\mathcal{O}_{\widetilde{Z}}), G\right)
\cong \operatorname{Ext}^*_{\mathrm{D}^b(Y)}\left(\mathcal{O}_Z,\, \pi_{Y*}G\right) = 0.
\]
Using the long exact sequences induced from the iterated extensions generating \( G \), we conclude that
\[
\operatorname{Ext}^*_{[\mathrm{D}^b(Y^+)]}(G, G) = 0.
\]
This implies that \( G \simeq 0 \) in the Verdier quotient \( [\mathrm{D}^b(Y^+)] \), so \( G \) belongs to the subcategory \( \left\langle i^+_*\mathcal{A}_{-1} \right\rangle \).
Moreover, by Corollary~\ref{sodofker}, the above decomposition is a well-defined semi-orthogonal decomposition.
\end{proof}

By Lemma~\ref{semiorthofserre}, we obtain the following consequence:
\begin{corollary}\label{cor[D]kernal}
Let \( G \in [\mathrm{D}^b(Y^+)] \) be an object such that its derived pushforward vanishes: \( \pi_{Y*} G = 0 \). Then \( G \) belongs to the triangulated subcategory
\[
\left\langle
i^+_*\mathcal{O}_{E_Y}(-(c-d)+1),\, \ldots,\, i^+_*\mathcal{O}_{E_Y}(-1)
\right\rangle.
\]
In particular, this gives an admissible subcategory of \([\mathrm{D}^b_{E_Y}(Y^+)]\) generated by a sequence of pre-tilting objects.
\end{corollary}

\subsection{Global adjunction}
We recall that in Corollary \ref{constructionlemma}, for any object \( F \in \mathrm{D}^b(Y) \), we constructed an object \( \mathcal{R}\pi^!_Y(F) \in \mathrm{D}^b(Y^+) \). However, this construction involves a certain ambiguity: it depends on the choice of a Bourbaki sequence associated to the Cohen–Macaulay approximation \( M^T_F \) of \( F \).

To reflect this, we denote by \( \mathcal{R}^{(1)}\pi^!_Y(F) \) and \( \mathcal{R}^{(2)}\pi^!_Y(F) \) two constructions corresponding to different choices of Bourbaki sequences or even some other freedom. It is clear that both satisfy the conclusions of Corollary~\ref{constructionlemma}.
\begin{definition}
we define
\[
[\mathrm{D}^b(Y^+)]^! := \mathrm{D}^b(Y^+) \big/ \langle i^+_*\mathcal{A}_{-(c-d)} \rangle
\]
to be the Verdier quotient of \(\mathrm{D}^b(Y^+)\) along the full saturated subcategory \(\langle i^+_*\mathcal{A}_{-(c-d)} \rangle\).
\end{definition}
We now state the following lemma which can be regarded as a global counterpart of Lemma~\ref{keylemmakernal}:
\begin{lemma}\label{keylemmaglobal}
For any object \( G \in \mathrm{D}^b(Y^+) \), there exists a natural isomorphism
\[
\operatorname{Ext}^*_{[\mathrm{D}^b(Y^+)]^!}\left(G, \mathcal{R}^{(1)}\pi^!_Y(F)\right) \cong \operatorname{Ext}^*_{\mathrm{D}^b(Y)}\left(\pi_{Y*}G,\,F\right),
\]
and this isomorphism is induced by the derived pushforward functor \( \pi_{Y*} \colon \mathrm{D}^b(Y^+) \to \mathrm{D}^b(Y) \).
\end{lemma}
\begin{proof}
We first consider an arbitrary anti roof:
\[
\begin{tikzcd}[column sep=small, row sep=small]
G\arrow[dr, "f"] & &\mathcal{R}^{(1)}\pi^!_Y(F)\arrow[dl, "t"'] \\
& \mathcal{F}^-&
\end{tikzcd}
\]
such that \( \operatorname{Cone}(t) \in\langle i^+_*\mathcal{A}_{-1} \rangle\). Then applying the derived pushforward \( \pi_{Y*} \) to \( f \), we get a morphism
\[
\pi_{Y*}(f): \pi_{Y*}(G) \to \pi_{Y*}(\mathcal{F}^-).
\]
Since \( \pi_{Y*}(\mathcal{F}^+) \cong F \) by Corollary \ref{pushforwardlemma}, this gives a morphism from \( F \) to \( \pi_{Y*}(G) \). Then we proceed to show that it induces an isomorphism by an argument analogous to that in \cite[Proof of Theorem 2.14]{BondalKapranovSchechtman2018}.\\

The surjectivity of \(\pi_{Y*}\) is straightforward. Note that for any morphism
\[
g: \pi_{Y*}G \longrightarrow F \quad \text{in} \quad \mathrm{D}^b(Y),
\]
by adjunction, there exists a morphism
\[
f := G\longrightarrow \pi_Y^!F \quad \text{in} \quad \mathrm{D}^+(Y^+)
\]
such that \(\pi_{Y*}(f) = g\).
On the other hand, consider the filtration of \(\pi_Y^!F\) constructed previously:
\[
\begin{tikzcd}[column sep=0.5em]
 & \mathcal{R}^{(1)}\pi^!_Y(F)\arrow{rr}&& \mathcal{F}^{s'} \arrow{rr} \arrow{dl} && \pi_Y^!F  \arrow{dl} \\
 && \tau_!^1 \arrow[ul,dashed,"\Delta"] && \tau_!^2 \arrow[ul,dashed,"\Delta"]
\end{tikzcd}
\]
If we choose \(s' \gg s\), then the cohomological degrees supporting \(\tau_!^2\) are sufficiently high, so we have an isomorphism
\[
\operatorname{Hom}_{\mathrm{D}^+(Y^+)}(G,\pi_Y^!F) \cong \operatorname{Hom}_{\mathrm{D}^+(Y^+)}(G,\mathcal{F}^{s'}).
\]
Under this isomorphism, \(f\) induces a bounded morphism, which we also denote by \(f\). Therefore, we obtain an anti roof
\[
\begin{tikzcd}[column sep=small, row sep=small]
G\arrow[dr, "f"] & &\mathcal{R}^{(1)}\pi^!_Y(F)\arrow[dl, "t"'] \\
& \mathcal{F}^{s'} &
\end{tikzcd}
\]
such that \(\mathrm{Cone}(t) = \tau^1_!\) lies in the subcategory \(\langle i^+_*\mathcal{A}_{-(c-d)} \rangle\), and
\[
\pi_{Y*}(t^{-1}\circ f) = g.
\]

For the injectivity of \(\pi_{Y*}\), suppose we are given a roof
\[
\begin{tikzcd}[column sep=small, row sep=small]
G\arrow[dr, "f"] & &\mathcal{R}^{(1)}\pi^!_Y(F)\arrow[dl, "t"'] \\
& \mathcal{F}^-&
\end{tikzcd}
\]
such that \( \pi_{Y*}(f) = 0 \). By adjunction, this implies that the induced morphism
\[
f: G\longrightarrow \pi_Y^! \pi_{Y*}\mathcal{F}^-
\]
in \( \mathrm{D}^+(Y^+) \) is the zero map.
Noting that \( \pi_Y^! \pi_{Y*}\mathcal{F}^- \cong \pi_Y^! F \), we consider the following filtration:
\[
\begin{tikzcd}[column sep=0.5em]
 & \mathcal{R}^{(1)}\pi^!_Y(F)\arrow{rr}&& \mathcal{F}^{s'} \arrow{rr} \arrow{dl} && \pi_Y^!F  \arrow{dl} \\
 && \tau_!^1 \arrow[ul,dashed,"\Delta"] && \tau_!^2 \arrow[ul,dashed,"\Delta"]
\end{tikzcd}
\]
If we choose \( s' \gg s \), then the cohomological degrees supporting \( \tau^1_! \) are sufficiently high. Hence we obtain the isomorphism
\[
\operatorname{Hom}_{\mathrm{D}^+(Y^+)}\left(G,  \pi_Y^! \pi_{Y*}\mathcal{F}^- \right)
\cong \operatorname{Hom}_{\mathrm{D}^+(Y^+)}\left(G,\mathcal{F}^{s'}\right),
\]
under which the morphism \( f \) corresponds to a morphism from \( G\) to \( \mathcal{F}^{s'} \), which is necessarily zero.
On the other hand, consider the unit morphism in \( \mathrm{D}^+(Y^+) \) induced by adjunction:
\[
f: \mathcal{F}^-\longrightarrow \pi_Y^! \pi_{Y*}\mathcal{F}^-.
\]
Again, using the same type of filtration, we choose \( s' \gg s \) so that
\[
\operatorname{Hom}_{\mathrm{D}^+(Y^+)}\left(\mathcal{F}^-, \pi_Y^! \pi_{Y*}\mathcal{F}^-\right)
\cong \operatorname{Hom}_{\mathrm{D}^+(Y^+)}\left(\mathcal{F}^-, \mathcal{F}^{s'}\right),
\]
and this isomorphism induces a morphism from \( \mathcal{F}^- \) to \(\mathcal{F}^{s'}\).
Combining these, we have a morphism from
\[
\begin{tikzcd}[column sep=small, row sep=small]
G\arrow[dr, "f"] & &\mathcal{R}^{(1)}\pi^!_Y(F)\arrow[dl, "t"'] \\
& \mathcal{F}^-&
\end{tikzcd}
\qquad\text{to}\qquad
\begin{tikzcd}[column sep=small, row sep=small]
G\arrow[dr, "0"] & &\mathcal{R}^{(1)}\pi^!_Y(F)\arrow[dl, "t'"'] \\
& \mathcal{F}^{s'}&
\end{tikzcd}
\]
In fact, according to the definition of equivalent morphisms in the Verdier quotient, this shows that $ t^{-1}\circ f = 0$ in $[\mathrm{D}^b(Y^+)]^!$.
\end{proof}

In the special case of Lemma~\ref{keylemmaglobal}, we consider the following isomorphisms:
\[
\operatorname{Ext}^*_{[\mathrm{D}^b(Y^+)]^!}\left(\mathcal{R}^{(2)}\pi^!_Y(F),\, \mathcal{R}^{(1)}\pi^!_Y(F)\right)
\cong
\operatorname{Ext}^*_{\mathrm{D}^b(Y)}\left(\pi_{Y*}\mathcal{R}^{(1)}\pi^!_Y(F),\,F\right)
\cong
\operatorname{Ext}^*_{\mathrm{D}^b(Y)}\left(F,\,F\right),
\]
where the last isomorphism follows from Corollary~\ref{constructionlemma}. In particular, we obtain a morphism corresponding to the identity map \( \mathrm{id}_F \), given by the following anti-roof:
\[
\begin{tikzcd}[column sep=small, row sep=small]
\mathcal{R}^{(2)}\pi^!_Y(F)\arrow[dr, "f"] & &\mathcal{R}^{(1)}\pi^!_Y(F)\arrow[dl, "t"'] \\
& \mathcal{F}^- &
\end{tikzcd}
\]
We now consider the distinguished triangle associated to this morphism:
\[
\mathcal{R}^{(2)}\pi^!_Y(F) \xrightarrow{f} \mathcal{F}^- \to G := \operatorname{Cone}(f) \to \mathcal{R}^{(2)}\pi^!_Y(F)[1].
\]
By applying the derived pushforward, we have \( \pi_{Y*}G = 0 \). Then by Proposition~\ref{proppushG0} and the observation that for any \( 1 \leq i \leq (c - d) - 1 \) and \( j = 1, 2 \), we have
\[
\operatorname{Ext}^*_{[\mathrm{D}^b(Y^+)]^!}
\left(i^+_*\mathcal{O}_{E_Y}(-i),\, \mathcal{R}^{(j)}\pi^!_Y(F)\right)
\cong
\operatorname{Ext}^*_{\mathrm{D}^b(Y)}
\left(\pi_{Y*}i^+_*\mathcal{O}_{E_Y}(-i),\, F\right)
\cong
0,
\]
where the vanishing follows from the fact that \( \pi_{Y*}i^+_*\mathcal{O}_{E_Y}(-i) = 0 \) in this range. We conclude that
\[
\operatorname{Ext}^*_{[\mathrm{D}^b(Y^+)]^!}
\left(i^+_*\mathcal{O}_{E_Y}(-i),\, G\right)
\cong
0,
\]
for any \( 1 \leq i \leq (c - d) - 1 \), so
\[
G \in \left\langle i^+_*\mathcal{A}_{-(c-d)} \right\rangle.
\]
Hence using the anti-roof above, we find that \( \mathcal{R}^{(1)}\pi^!_Y(F) \cong \mathcal{R}^{(2)}\pi^!_Y(F) \) in the category \( [\mathrm{D}^b(Y^+)]^! \). Therefore, we use the unified notation \( \mathcal{R}\pi^!_Y(F) \) to denote their common projection in the quotient category.\\

If we are given a distinguished triangle \( F \xrightarrow{g} G \to H \) in \( \mathrm{D}^b(Y) \), then applying the functor \( \pi_Y^! \), we obtain a distinguished triangle
\[
\pi_Y^!F \xrightarrow{\pi_Y^!g} \pi_Y^!G \to \pi_Y^!H
\]
in \( \mathrm{D}^+(Y^+) \).
Recall that we have the following filtration on \( \pi^!_Y F \) (or on \( \pi^!_Y G \)):
\[
\begin{tikzcd}[column sep=0.1em]
 & \mathcal{R}\pi^!_Y(F) \arrow{rr} && \mathcal{F}_!^{1}(F) \arrow{dl} \\
 && \operatorname{gr}_!^1 \arrow[ul, dashed, "\Delta"]
\end{tikzcd}
\cdots
\begin{tikzcd}[column sep=0.01em]
 & \mathcal{F}_!^{s-1}(F) \arrow{rr} && \mathcal{F}_!^{s}(F) \arrow{rr} \arrow{dl} && \pi_Y^!F \arrow{dl} \\
 && \operatorname{gr}_!^{s} \arrow[ul, dashed, "\Delta"] && \tau_!^{s} \arrow[ul, dashed, "\Delta"]
\end{tikzcd}
\]
According to a similar discussion as before, we can choose integers \( s' \gg s'' \gg 0 \) such that there exists an induced morphism
\[
\mathcal{F}_!^{s''}(F) \longrightarrow \mathcal{F}_!^{s'}(G)
\]
such that it fits in the following commutative diagram:
\[
\begin{tikzcd}
\mathcal{F}_!^{s''}(F) \arrow[r] \arrow[d] & \mathcal{F}_!^{s'}(G) \arrow[d] \\
\pi_Y^!F \arrow[r] & \pi_Y^!G
\end{tikzcd}
\]
which also represents the morphism \( f \colon \mathcal{R}\pi^!_Y(F) \to \mathcal{R}\pi^!_Y(G) \) in \( [\mathrm{D}^b(Y^+)]^! \) corresponding to a given morphism \( g \) in \( \mathrm{D}^b(Y) \) via adjunction. Therefore, we obtain the following distinguished triangle in \( \mathrm{D}^b(Y^+) \):
\[
\mathcal{F}_!^{s''}(F) \longrightarrow \mathcal{F}_!^{s'}(G) \longrightarrow \mathcal{H}_1 := \mathrm{Cone}(f),
\]
such that the image of \( \mathcal{H}_1 \) in \( [\mathrm{D}^b(Y^+)]^! \) is isomorphic to \( \operatorname{Cone}(f)\). Then we can complete the above commutative diagram
\[
\begin{tikzcd}
\mathcal{F}_!^{s''}(F) \arrow[r] \arrow[d] & \mathcal{F}_!^{s'}(G) \arrow[r] \arrow[d] & \mathcal{H}_1 \arrow[r] \arrow[d,"h"] & \mathcal{F}_!^{s''}(F)[1]\arrow[d] \\
\pi_Y^!F \arrow[r, "\pi_Y^!g"] & \pi_Y^!G \arrow[r] & \pi_Y^!H \arrow[r] & \pi_Y^!F[1]
\end{tikzcd}
\]
Correspondingly, we consider the filtration on \( \pi_Y^!H \) as above. Then for sufficiently large \( s \gg 0 \), we obtain the following anti-roof:
\[
\begin{tikzcd}[column sep=small, row sep=small]
\mathcal{H}_1 \arrow[dr, "h"] & & \mathcal{R}\pi^!_Y(H) \arrow[dl, "t"'] \\
& \mathcal{H}^-:=\mathcal{F}_!^{s}(H) &
\end{tikzcd}
\]
Moreover by our setting, it is not difficult to see that the distinguished triangle
\[
\mathcal{H}_1 \xrightarrow{h} \mathcal{H}^- \to \operatorname{Cone}(h)
\]
satisfies \( \pi_{Y*}(\operatorname{Cone}(h)) = 0 \). Similarly:
\[
\operatorname{Ext}^*_{[\mathrm{D}^b(Y^+)]^!}
\left(i^+_*\mathcal{O}_{E_Y}(-i),\, \operatorname{Cone}(h)\right)
\cong
0,
\]
for any \( 1 \leq i \leq (c - d) - 1 \), so by Proposition \ref{proppushG0}
\[
\operatorname{Cone}(h) \in \left\langle i^+_*\mathcal{A}_{-(c-d)} \right\rangle.
\]
We have \( \mathcal{H}_1 \simeq \mathcal{R}\pi^!_Y(H) \) in \( [\mathrm{D}^b(Y^+)]^! \). Finally, we obtain the following distinguished triangle in the quotient category \( [\mathrm{D}^b(Y^+)]^! \):
\[
\mathcal{R}\pi^!_Y(F) \longrightarrow \mathcal{R}\pi^!_Y(G) \longrightarrow \mathcal{R}\pi^!_Y(H),
\]
while the last morphism is also induced by adjunction.\\

Finally, we verify that the correspondence in Lemma~\ref{keylemmaglobal} is indeed natural with respect to the above morphisms. Thus we have:
\begin{prop}\label{maniprop}
Based on all the constructions above, if \( c - d > 0 \) we obtain a well-defined, left admissible, and fully-faithful exact triangulated functor:
\[
\mathcal{R}\pi^!_Y \colon \mathrm{D}^b(Y) \longrightarrow [\mathrm{D}^b(Y^+)]^!
\]
from the bounded derived category of \( Y \) to the Verdier quotient \( [\mathrm{D}^b(Y^+)]^! \).
\end{prop}
Up to a twist by \( -K_{Y^+/Y} \), we have the following:
\begin{prop}\label{maniprop1}
If \( c - d > 0 \), then we obtain a well-defined, left admissible, and fully faithful exact triangulated functor
\[
\widetilde{\mathcal{R}}\pi^*_Y \colon \mathrm{D}^b(Y) \longrightarrow [\mathrm{D}^b(Y^+)]
\]
from the bounded derived category of \( Y \) to the Verdier quotient \( [\mathrm{D}^b(Y^+)] \).
\end{prop}
\begin{theorem}\label{thmsod}
If \( c - d > 0 \),  we have the following semi-orthogonal decompositions:
\[
[\mathrm{D}^b(Y^+)]^!
= \left\langle
\operatorname{Im} \mathcal{R}\pi^!_Y,\
i^+_*\mathcal{O}_{E_Y}(-(c{-}d)+1),\
\ldots,\
i^+_*\mathcal{O}_{E_Y}(-2),\
i^+_*\mathcal{O}_{E_Y}(-1)
\right\rangle.
\]
or dually,
\[
[\mathrm{D}^b(Y^+)]
= \left\langle \operatorname{Im} \widetilde{\mathcal{R}}\pi^*_Y,\
i^+_*\mathcal{O}_{E_Y},\
i^+_*\mathcal{O}_{E_Y}(1),\
\ldots,\
i^+_*\mathcal{O}_{E_Y}((c{-}d){-}2) \right\rangle.
\]
\end{theorem}
\begin{proof}
Lemma~\ref{semiorthofserre} and Lemma~\ref{keylemmaglobal}  imply that the following category
\[
\mathcal{D} :=
\left\langle
\operatorname{Im} \mathcal{R}\pi^!_Y,\
i^+_*\mathcal{O}_{E_Y}(-(c{-}d)+1),\
\ldots,\
i^+_*\mathcal{O}_{E_Y}(-2),\
i^+_*\mathcal{O}_{E_Y}(-1)
\right\rangle
\]
forms a semi-orthogonal decomposition of a left admissible subcategory $\mathcal{D} \subset [\mathrm{D}^b(Y^+)]^!$.\\

For the generated properties, for any object in \(\mathrm{D}^b(Y^+)\) we have a natural morphism \( G \to \mathcal{R}\pi^!_Y\pi_{Y*}(G)\) in \([\mathrm{D}^b(Y^+)]^!\) by adjunction, which induces a distinguished triangle
\[
G \longrightarrow \mathcal{R}\pi^!_Y\pi_{Y*}(G) \longrightarrow H
\]
in [\( \mathrm{D}^b(Y^+)]^! \). By construction, the object \( H \) satisfies \( \pi_{Y*}H=0 \). Therefore, by Proposition~\ref{proppushG0} we conclude that \( H \in \mathcal{D} \), and hence \( G \in \mathcal{D} \).
\end{proof}

According to Lemma~\ref{keylemmaglobal} , we have
\[
\operatorname{Ext}^*_{[\mathrm{D}^b(Y^+)]^!}
\left(F^s(\mathcal{O}_{\widetilde{Z}}),\ \mathcal{R}\pi^!_Y(\mathcal{O}_Z)\right)
\cong
\operatorname{Ext}^*_{\mathrm{D}^b(Y)}(\mathcal{O}_Z,\ \mathcal{O}_Z)
\]
for any \( s \geq c-d-1 \). Let us now consider the identity morphism \( \operatorname{id}_{\mathcal{O}_Z} \), which corresponds to a roof of the form
\[
\begin{tikzcd}[column sep=small, row sep=small]
& \mathcal{F}^+ \arrow[dl, "t"'] \arrow[dr, "f"] & \\
F^s(\mathcal{O}_{\widetilde{Z}}) & & \mathcal{R}\pi^!_Y(\mathcal{O}_Z)
\end{tikzcd}
\]
in the Verdier quotient. Therefore, \( F^s(\mathcal{O}_{\widetilde{Z}}) \) is isomorphic to \( \mathcal{R}\pi^!_Y(\mathcal{O}_Z) \) up to the quotient by the kernel of \( \pi_{Y*} \), which is by Lemma~\ref{proppushG0} generated by
\[
\left\langle
\left\langle i^+_*\mathcal{A}_{-(c-d)} \right\rangle,\
i^+_*\mathcal{O}_{E_Y}(-(c{-}d)+1),\
\ldots,\
i^+_*\mathcal{O}_{E_Y}(-1)
\right\rangle.
\]
by the orthogonality of the pretilting objects with \( F^s(\mathcal{O}_{\widetilde{Z}}) \) and \( \mathcal{\widetilde{R}}\pi^*_Y(\mathcal{O}_Z) \) established in Lemma \ref{localpushspan}, we obtain an isomorphism
\[
F^s(\mathcal{O}_{\widetilde{Z}}) \simeq \mathcal{\widetilde{R}}\pi^*_Y(\mathcal{O}_Z) \quad \text{in} \quad {[\mathrm{D}^b(Y^+)]}.
\]
Naturally the Lemma \ref{keylemmaglobal} has its twin lemma, we denote by \( \mathcal{R}^{(1)}\pi^*_Y(F) \) any construction under a Bourbaki sequence choice:
\begin{lemma}\label{keylemmagglobal2}
For any object \( G \in \mathrm{D}^b(Y^+) \), there exists a natural isomorphism
\[
\operatorname{Ext}^*_{[\mathrm{D}^b(Y^+)]}\left(\mathcal{R}^{(1)}\pi^*_Y(F), G\right) \cong \operatorname{Ext}^*_{\mathrm{D}^b(Y)}\left(F,\,\pi_{Y*}G\right),
\]
and this isomorphism is induced by the derived pushforward functor \( \pi_{Y*} \colon \mathrm{D}^b(Y^+) \to \mathrm{D}^b(Y) \).
\end{lemma}
In analogy with the previous discussion, and by Lemma~\ref{keylemmagglobal2}, we see that if \( s \geq c-d-1 \), \( F^s(\mathcal{O}_{\widetilde{Z}}) \) is isomorphic to \( \mathcal{R}\pi^*_Y(\mathcal{O}_Z) \) up to the quotient by the kernel of \( \pi_{Y*} \),
\[
\left\langle
i^+_*\mathcal{O}_{E_Y}(-(c{-}d)+1),\, \ldots,\, i^+_*\mathcal{O}_{E_Y}(-1),\, \left\langle i^+_*\mathcal{A}_{-1} \right\rangle
\right\rangle.
\]
Moreover, by adjunction we have
\[
\operatorname{Ext}^*_{[\mathrm{D}^b(Y^+)]}\left(\mathcal{R}^{(1)}\pi^*_Y(F),\, i^+_*\mathcal{O}_{E_Y}(-i)\right)
= \operatorname{Ext}^*_{[\mathrm{D}^b(Y^+)]}\left(F^s(\mathcal{O}_{\widetilde{Z}}),\, i^+_*\mathcal{O}_{E_Y}(-i)\right)
= 0
\]
for any \( 1 \leq i \leq c{-}d{-}1 \). Therefore, we conclude that \( F^s(\mathcal{O}_{\widetilde{Z}}) \) is isomorphic to \( \mathcal{R}\pi^*_Y(\mathcal{O}_Z) \) in the Verdier quotient category \( [\mathrm{D}^b(Y^+)] \).\\

There are the following twin properties:
\begin{prop}\label{prop:output}
Based on all the constructions above, if \( c - d > 0 \) we obtain a well-defined, right admissible, and fully-faithful exact triangulated functor:
\[
\mathcal{R}\pi^*_Y \colon \mathrm{D}^b(Y) \longrightarrow [\mathrm{D}^b(Y^+)]
\]
from the bounded derived category of \( Y \) to the Verdier quotient \( [\mathrm{D}^b(Y^+)] \).
\end{prop}
Up to a twist by \( K_{Y^+/Y} \), we have the following:
\begin{prop}\label{prop:output1}
If \( c - d > 0 \), then we obtain a well-defined, right admissible, and fully faithful exact triangulated functor
\[
\widetilde{\mathcal{R}}\pi^!_Y \colon \mathrm{D}^b(Y) \longrightarrow [\mathrm{D}^b(Y^+)]^!
\]
from the bounded derived category of \( Y \) to the Verdier quotient \( [\mathrm{D}^b(Y^+)]^! \).
\end{prop}

\begin{remark}
The symbol \( \mathcal{R} \) originates from the construction of the Rees module. When \( P = \mathfrak{m} \) and \( d = 2 \), we have for any Cohen–Macaulay module \(M\):
\[
\mathcal{R}(M)\simeq\operatorname{Rees}(M)\simeq\mathbf{L}_0\pi^*_Y M / \mathrm{Tor}.
\]
Hence, we refer to the above construction as a \emph{derived Rees (type) pullback}.
\end{remark}

Finally, we also have a mirror-symmetric counterpart of Theorem~\ref{thmsod}.
\begin{theorem} \label{thm:output}
If \( c - d > 0 \),  we have the following semi-orthogonal decompositions:
\[
[\mathrm{D}^b(Y^+)]
= \left\langle
i^+_*\mathcal{O}_{E_Y}(-(c{-}d)+1),\
\ldots,\
i^+_*\mathcal{O}_{E_Y}(-2),\
i^+_*\mathcal{O}_{E_Y}(-1),\
\operatorname{Im} \mathcal{R}\pi^*_Y
\right\rangle.
\]
or dually,
\[
[\mathrm{D}^b(Y^+)]^!
= \left\langle
i^+_*\mathcal{O}_{E_Y},\
i^+_*\mathcal{O}_{E_Y}(1),\
\ldots,\
i^+_*\mathcal{O}_{E_Y}((c{-}d){-}2),\
\operatorname{Im} \widetilde{\mathcal{R}}\pi^!_Y \right\rangle,
\]
\end{theorem}

In analogy with Verdier--Artin theory, our main results are formulated under the assumption that $R$ is a formal or analytic local ring. The primary reason for this restriction is that we rely on standard basis techniques, which are classically well-defined only in the formal or analytic setting. In particular, our standard resolution is constructed only at the formal level. Nevertheless, when $R$ is an \textbf{algebraic} regular local ring and $M$ an algebraic module, the associated Rees transform, for instance $\mathcal{R}^!(M)$, exhibits behaviour that is also algebraic.
\[
\begin{tikzcd}
\widehat{Y^+_{E_Y}}\arrow[r] \arrow[d] & Y^+ \arrow[d] \\
\widehat{Y} \arrow[r] & Y
\end{tikzcd}
\]
Consider the distinguished triangle in \(\mathrm{D}^+\bigl(\widehat{Y^+_{E_Y}}\bigr)\) we constructed before:
\[
\mathcal{R}^!(\widehat{M})
\longrightarrow
\pi_{\widehat{Y}}^!\widehat{M}
\longrightarrow
\tau
\longrightarrow
\mathcal{R}^!(\widehat{M})[1].
\]
Here $\pi_{\widehat{Y}}^!\widehat{M}$ is algebraic from the commutative diagram, and via the identification
\[
\mathrm{Coh}_{E_Y}(Y^+) = \mathrm{Coh}_{E_{Y}}\bigl(\widehat{Y^+_{E_Y}}\bigr),
\]
the object $\tau$ is algebraic as well if we choose a good well--order in (\ref{orderchoice}) \footnote{The rigorous  proof should be given in the complete intersection case.}. Moreover, any morphism from $\pi_{\widehat{Y}}^!\widehat{M}$ to a sheaf supported on $E_Y$ is algebraic. Consider the cohomological truncation of $\pi_{\widehat{Y}}^! M$ in the standard $t$-structure on $\mathrm{D}^b_{E_Y}(\widehat{Y}^+)$. There is a canonical distinguished triangle
\begin{equation}\label{truncational}
\tau^{\le N-1} \pi_{\widehat{Y}}^! \widehat{M}
\longrightarrow
\pi_{\widehat{Y}}^! \widehat{M}
\longrightarrow
\tau^{\ge N} \pi_{\widehat{Y}}^! \widehat{M}
\longrightarrow
(\tau^{\le N-1} \pi_{\widehat{Y}}^! \widehat{M})[1],
\end{equation}
where $\tau^{\le N-1}$ and $\tau^{\ge N}$ denote the standard cohomological truncation functors for a sufficiently large integer $N$. Such that we have the following commutative triangles:
\[
\begin{tikzcd}[row sep=small, column sep=small]
\tau^{\le N-1} \mathcal{R}^!(\widehat{M}) \simeq \mathcal{R}^!(\widehat{M}) \arrow[r] \arrow[d]
& \tau^{\le N-1} \pi_{\widehat{Y}}^! \widehat{M} \arrow[r] \arrow[d]
& \tau^{\le N-1} \tau \arrow[r] \arrow[d]
& (\tau^{\le N-1} \mathcal{R}^!(\widehat{M}))[1] \arrow[d] \\
\mathcal{R}^!(\widehat{M}) \arrow[r] \arrow[d]
& \pi_{\widehat{Y}}^! \widehat{M} \arrow[r] \arrow[d]
& \tau \arrow[r] \arrow[d]
& \mathcal{R}^!(\widehat{M})[1] \arrow[d] \\
\tau^{\ge N} \mathcal{R}^!(\widehat{M}) \simeq 0 \arrow[r]
& \tau^{\ge N} \pi_{\widehat{Y}}^! \widehat{M} \arrow[r]
& \tau^{\ge N} \tau \arrow[r]
& (\tau^{\ge N} \mathcal{R}^!(\widehat{M}))[1]
\end{tikzcd}
\]
We need to show that the morphism
\[
\tau^{\le N-1} \pi_{\widehat{Y}}^! \widehat{M}
\longrightarrow
\tau^{\le N-1} \tau
\]
arising above is the completion of some algebraic morphism. If this is the case, then this suggests that
\[
\mathcal{R}^!(\widehat{M})
\]
itself should be algebraic since the cohomological truncation preserves algebraicity of \(\pi_{\widehat{Y}}^! \widehat{M}\) and \(\tau\).
\begin{lemma}
There is a natural isomorphism
\[
\operatorname{Ext}^*_{\mathrm{D}^+(\widehat{Y}^+_{E_Y})}\Bigl(
\tau^{\le N-1} \pi_{\widehat{Y}}^! \widehat{M},
\tau^{\le N-1} \tau
\Bigr)
\simeq
\operatorname{Ext}^*_{\mathrm{D}^+(Y^+)}\Bigl(
\tau^{\le N-1} \pi_{Y}^! M,
\tau^{\le N-1} \tau
\Bigr).
\]
\end{lemma}
\begin{proof}
Consider the distinguished triangle (\ref{truncational}) in $\mathrm{D}^+(\widehat{Y}^+_{E_Y})$ arising from the cohomological truncation of $\pi_{\widehat{Y}}^! \widehat{M}$, by considering the long exact sequence in $\operatorname{Ext}^*$ induced by this triangle, it suffices to show that
\begin{equation}\label{eq:ext-high-trunc}
\operatorname{Ext}^*_{\mathrm{D}^+(\widehat{Y}^+_{E_Y})}\Bigl(
\tau^{\ge N} \pi_{\widehat{Y}}^! \widehat{M},
\tau^{\le N-1} \tau
\Bigr)
\simeq
\operatorname{Ext}^*_{\mathrm{D}^+(Y^+)}\Bigl(
\tau^{\ge N} \pi_{Y}^! M,
\tau^{\le N-1} \tau
\Bigr).
\end{equation}
and
\begin{equation}\label{eq:ext-low-trunc}
\operatorname{Ext}^*_{\mathrm{D}^+(\widehat{Y}^+_{E_Y})}\Bigl(
\pi_{\widehat{Y}}^! \widehat{M},
\tau^{\le N-1} \tau
\Bigr)
\simeq
\operatorname{Ext}^*_{\mathrm{D}^+(Y^+)}\Bigl(
\pi_{Y}^! M,
\tau^{\le N-1} \tau
\Bigr),
\end{equation}
For \eqref{eq:ext-high-trunc}, we note that there is a natural identification of derived categories with support on $E_Y$:
\[
\mathrm{D}^+\big( \mathrm{Coh}_{E_Y}(Y^+)\big)
\simeq
\mathrm{D}^+ \big(\mathrm{Coh}_{E_{Y}}\bigl(\widehat{Y^+_{E_Y}}\bigr)\big),
\]
and moreover, we have
\[
\tau^{\ge N} \pi_Y^! M
\simeq
\tau^{\ge N} \pi_{\widehat{Y}}^! \widehat{M}
\simeq
\tau^{\ge N} \tau .
\]
which all support on with support on $E_Y$. For \eqref{eq:ext-low-trunc}, we note that by construction
\[
\tau^{\le N-1} \in \mathrm{D}^b_{E_Y}(Y^+),
\]
i.e., it lies in the bounded derived category of coherent sheaves supported on $E_Y$.
Then, using induction on cohomology and the following natural identifications of Ext groups for any sheaf $F$ support on $E_Y$:
\[
\begin{aligned}
\operatorname{Ext}_{Y^+}^*\bigl(\pi_Y^! M, i^+_* F\bigr)
&\simeq \operatorname{Ext}_{E_Y}^*\bigl((i^+)^* \pi_Y^! M, F\bigr) \\
&\simeq \operatorname{Ext}_{E_Y}^*\bigl((i^+)^* \pi_{\widehat{Y}}^! \widehat{M}, F\bigr) \\
&\simeq \operatorname{Ext}_{\widehat{Y^+_{E_Y}}}^*\bigl(\pi_{\widehat{Y}}^! \widehat{M}, i^+_* F\bigr),
\end{aligned}
\]
Since \(i^+\) is a closed embedding of a Cartier divisor, the isomorphism in the first step is justified.
Therefore, we obtain the desired natural isomorphism of Ext groups.
\end{proof}
Hence, we denote by \(\mathcal{R}^!(M)\) the algebraic counterpart of \(\mathcal{R}^!(\widehat{M})\). Similarly, for the other constructed functors, we have analogous algebraic descriptions. Thus, we have reached the desired outcome that we anticipated from the outset.
\begin{theorem}\label{thm:algebraic-version}
Our Propositions~\ref{maniprop}, \ref{maniprop1}, \ref{prop:output}, \ref{prop:output1} and Theorems~\ref{thmsod}, \ref{thm:output} continue to hold if we assume that our regular local ring is \textbf{algebraic}, in our primary setting~(\ref{expla}).
\end{theorem}
\begin{remark}
In the work of Kuznetsov and Lunts~\cite[Definition~6.1]{KuznetsovLunts2015}, the notion of a
\emph{non-rational center} is introduced as follows.
Let $f \colon X \to Y$ be a proper birational morphism.
A subscheme $Z \subset Y$ is called a non-rational center of $Y$ with respect to $f$
if the canonical morphism
\[
\mathcal I_Z \longrightarrow Rf_*\bigl(\mathcal I_{f^{-1}(Z)}\bigr)
\]
is an isomorphism. In view of earlier work of Mori~\cite{Mori1982}, where the proof of (3.32.2) is given over the base field $\mathds{C}$,
it follows that in the setting relevant to us—namely, when $f$ is the blow-up of a complete intersection
and the exceptional divisor is Fano—for any non-negative integer $n$, the $n$-th infinitesimal neighborhood
$Z^{(n)}$ is a non-rational center of $Y$.
\end{remark}
\begin{remark}
\begin{enumerate}
\item However, this method is not a direct construction, as an alternative, one may pass to a polynomial presentation of the algebraic regular local ring, for example
\[
\mathrm{k}[x_i]_{(x_i)}/I,
\]
and perform a division procedure analogous to Grauert division for algebraic complete intersections. This requires a sequence of refinements allowing one to extract the tail of the periodic division into identity components and to justify the process rigorously, with the aim of producing an asymptotically periodic presentation.
\item A further aspect of the argument that could be improved concerns the use of monomial-refined division. In the present approach we rely primarily on order division, using monomial division only to guarantee that over a regular local ring we obtain a finite resolution whose syzygies satisfy constraints on their leading terms. If one works solely with order division, such finiteness properties can no longer be expected, and this may be problematic. An alternative direction would be to identify suitable key varieties capturing the relevant geometric structure.
\item The complete intersection case requires a more systematic treatment. It is likely that an appropriate dg-algebra formalism would simplify the above constructions and diagrams \cite{Avramov2002}.
\end{enumerate}
\end{remark}

\section{Acknowledgements}

I would like to thank Professor Ryo Takahashi, Shinya Kumashiro, Xiaowu Chen, Yuki Mizuno, and Tomoki Yoshida for helpful discussions. I am especially grateful to my supervisor, Professor Yasunari Nagai, for his constant support and guidance throughout the development of this work.

\appendix
\section{Standard resolution on complete intersections}\label{subsec:ci}
The above argument can be further extended to the case where \( J \) is a complete intersection. Let \( (R, \mathfrak{m}) \) be a regular local ring, and let \( J \subset R \) be an ideal such that the quotient ring \( S := R/J \) is Gorenstein. Suppose \( P \) is a prime ideal in \( R \) containing \( J \) with the property that the residue ring \( Z := R/P \) is regular. We consider the situation where \( J \) is generated by a regular sequence of length \(2\), for instance, \( J = (w_1, w_2) \). We have the following commutative diagram:
\[
\begin{tikzcd}[row sep=large, column sep=large]
\operatorname{Spec} R/P \arrow[d, "k_2"'] \arrow[dr, "k_1"] \arrow[drr, "k", bend left=15] & & \\
\operatorname{Spec} R/J \arrow[r, "j_2"'] & \operatorname{Spec} R/(w_1) \arrow[r, "j_1"'] & \operatorname{Spec} R \\
\operatorname{Spec} R/I \arrow[u, "i_2"'] \arrow[ur, "i_1"'] \arrow[urr, "i", bend right=15]
\end{tikzcd}
\]
where all morphisms are induced by the corresponding quotient maps, and we use the same notation for the canonical morphism between corresponding rings and their analytic (formal) completion.

Our goal is to construct the standard resolution of $i_{2*}\mathcal{O}_B := i_{2*}(\operatorname{Spec} R/I)$ over $\operatorname{Spec} R/J$. First, we work over the analytic (formal) completion $\widehat{R}$ of $R$. From previous discussions, considering either $i_*\mathcal{O}_{\widehat{B}}$ or $i_{1*}\mathcal{O}_{\widehat{B}}$, we obtain the following standard resolution for \(\widehat{P}^{r}\mathcal{O}_{\widehat{B}}\) in  \(\widehat{R}/(w_1)\):

\[
\begin{tikzcd}[column sep=small, row sep=small]
 &\vdots \arrow[d] & \vdots \arrow[d] \\
 \bigoplus^{m_{1}}_{l=1} \widehat{P}^{r-\mathrm{a}^{(1)}_l-2d}(-2\widehat{S}_1) \arrow[r, "F_{1}"] \arrow[d, "-K_2"] & \bigoplus^{m_{0}}_{l=1} \widehat{P}^{r-\mathrm{a}^{(0)}_l-2d}(-2\widehat{S}_1) \arrow[r, "F_{0}"] \arrow[d, "K_1"] &  \widehat{P}^{r-2d}(-2\widehat{S}_1) \arrow[ddl,dashed, "-K^{(1)}_0"]\arrow[d, "-K_0"] \\
 \cdots\arrow[r, "F_{2}"] & \bigoplus^{m_{1}}_{l=1} \widehat{P}^{r-\mathrm{a}^{(1)}_l-d}(-\widehat{S}_1) \arrow[r, "F_{1}"] \arrow[d, "K_2"] & \bigoplus^{m_{0}}_{l=1} \widehat{P}^{r-\mathrm{a}^{(0)}_l-d}(-\widehat{S}_1) \arrow[r, "F_{0}"] \arrow[d, "-K_1"] &  \widehat{P}^{r-d}(-\widehat{S}_1) \arrow[d, "K_0"] \\
 \cdots & \bigoplus^{m_{2}}_{l=1} \widehat{P}^{r-\mathrm{a}^{(2)}_l} \arrow[r, "F_{2}"] & \bigoplus^{m_{1}}_{l=1} \widehat{P}^{r-\mathrm{a}^{(1)}_l} \arrow[r, "F_{1}"] & \bigoplus^{m_{0}}_{l=1} \widehat{P}^{r-\mathrm{a}^{(0)}_l} \arrow[r, "F_{0}"] & \widehat{P}^{r}
\end{tikzcd}
\]
where $r$ is any integer and $S_1$ is the hypersurface defined by $w_1$.  Furthermore, by the property of the regular sequence, $w_2$ defines a Cartier hypersurface $S_2$ in $\operatorname{Spec} R/(w_1)$. Therefore, by Corollary \ref{excesshyper} $i_{2*}\mathcal{O}_{B_s}$ admits a filtration with graded pieces given by:
\[
i_{2*}\mathcal{O}_B \longrightarrow i_{2*}\mathcal{O}_B(-sS_2)[2s] \longrightarrow i_{2*}\mathcal{O}_{B_s} \longrightarrow i_{2*}\mathcal{O}_B[1],
\]
and $i_{2*}\mathcal{O}_{B_s}$ admits the following filtration:
\[
\begin{tikzcd}[column sep=small, row sep=small]
 & 0 \arrow[rr] && i_{2*}\mathcal{O}_{B_1} \arrow[dl] \arrow[rr] && i_{2*}\mathcal{O}_{B_2} \arrow[dl] \\
 & & j_{2}^{*}j_{2*}i_{2*}\mathcal{O}_B[1] \arrow[ul, dashed, "\Delta"]
 && j_{2}^{*}j_{2*}i_{2*}\mathcal{O}_B(-S_2)[3] \arrow[ul, dashed, "\Delta"] &
\end{tikzcd}
\]
\[
\cdots
\]
\[
\begin{tikzcd}[column sep=small, row sep=small]
 & i_{2*}\mathcal{O}_{B_{r-1}} \arrow[rr] && i_{2*}\mathcal{O}_{B_s} \arrow[dl] \\
 & & j_{2}^{*}j_{2*}i_{2*}\mathcal{O}_B(-(s-1)S_2)[2s-1] \arrow[ul, dashed, "\Delta"]
\end{tikzcd}
\]
In the conventional setting, we consider sufficiently large integers $s$ and view the filtration structure of $i_{2*}\mathcal{O}_{B_s}$ as a filtration of $i_{2*}\mathcal{O}_{B}$. Crucially, like Lemma \ref{lemmaexccheck} we need to examine the morphisms within this filtration:
\[
\mathrm{Ext}^1 \left( j_2^{*}i_{1*} \mathcal{O}_{B}(- (n-1) S_2) [2n-1], j_2^{*}i_{1*} \mathcal{O}_{B}(- (m-1) S_2) [2m-1] \right)
=
\begin{cases}
\mathcal{O}_{B} & \text{for } n = m+1 \\
0 & \text{for } n > m+1
\end{cases}
\]
then by symmetry and the construction of the excess distinguished triangle, the non-trivial contributions in the filtration arise from the unit element in:
\[
\mathrm{Hom} \left( j_2^{*}i_{1*} \mathcal{O}_{B}(- S_2)[1], j_2^{*}i_{1*} \mathcal{O}_{B} \right)
\]
Similarly, we construct the complete system of above morphism in the analytic (formal) completion $\widehat{R}$.
\begin{enumerate}
\item As above, we have the standard resolution of $i_{2*}\mathcal{O}_{\widehat{B}}$ induced by the totalization of the following double factorization on $\widehat{R}$,
\[
\begin{tikzcd}[column sep=small, row sep=small]
& & \vdots \arrow[d] & & \vdots \arrow[d] \\
\cdots \arrow[r, "F_{3}"] & \mathcal{O}_{\widehat{R}}^{m_2}(-2\widehat{S_1}) \arrow[r, "F_{2}"] \arrow[d, "K_3"] & \mathcal{O}_{\widehat{R}}^{m_1}(-2\widehat{S_1}) \arrow[r, "F_{1}"] \arrow[d, "-K_2"] & \mathcal{O}_{\widehat{R}}^{m_0}(-2\widehat{S_1}) \arrow[ddl,dashed, "K^{(1)}_1"]\arrow[r, "F_{0}"] \arrow[d, "K_1"] & \mathcal{O}_{\widehat{R}}(-2\widehat{S_1})  \arrow[ddl,dashed, "-K^{(1)}_0"]\arrow[d, "-K_0"] \\
& \cdots \arrow[r, "F_{3}"] & \mathcal{O}_{\widehat{R}}^{m_2}(-\widehat{S_1}) \arrow[r, "F_{2}"] \arrow[d, "-K_3"] & \mathcal{O}_{\widehat{R}}^{m_1}(-\widehat{S_1}) \arrow[r, "F_{1}"] \arrow[d, "K_2"] & \mathcal{O}_{\widehat{R}}^{m_0}(-\widehat{S_1}) \arrow[r, "F_{0}"] \arrow[d, "-K_1"] & \mathcal{O}_{\widehat{R}}(-\widehat{S_1}) \arrow[d, "K_0"] \\
& & \cdots \arrow[r, "F_{3}"] & \mathcal{O}_{\widehat{R}}^{m_2} \arrow[r, "F_{2}"] & \mathcal{O}_{\widehat{R}}^{m_1} \arrow[r, "F_{1}"] & \mathcal{O}_{\widehat{R}}^{m_0} \arrow[r, "F_{0}"] & \mathcal{O}_{\widehat{R}}
\end{tikzcd}
\]In particular, for every integer $i$, the differentials satisfy the following relations:
\[
K_i \circ F_i + F_{i+1} \circ K_{i+1} = w_1 \cdot \text{id}
\]
\[
\sum_{\substack{a_1+a_2=a-1}} K_{i+2a_2+1}^{(a_1)} \circ K_i^{(a_2)} + F_{i+2a} \circ K_i^{(a)} + K_{i-1}^{(a)} \circ F_{i-1} = 0
\]
Then applying the $j_2^*$ pullback to this resolution, we obtain a complex for $j_2^* i_{2*}\mathcal{O}_{\widehat{B}}$ while preserving the names of all morphisms. Consider the morphism
\[
\mathrm{Hom} \left( j_2^{*} i_{1*} \mathcal{O}_{\widehat{B}}(- \widehat{S_2})[1],\ j_2^{*} i_{1*} \mathcal{O}_{\widehat{B}} \right),
\]
which arises from  two induced complexes described above. To define these morphisms, we begin by constructing the building blocks of a map between these two factorizations on $\widehat{R}$:
\[
\begin{tikzcd}[column sep=small, row sep=small]
& \vdots \arrow[d] & & \vdots \arrow[d] \\
\cdots \arrow[r, "F_{2}"] & \mathcal{O}_{\widehat{R}}^{m_1}(-2\widehat{S}_1-\widehat{S}_2) \arrow[r, "F_{1}"] \arrow[d, "-K_2"] & \mathcal{O}_{\widehat{R}}^{m_0}(-2\widehat{S}_1-\widehat{S}_2) \arrow[ddl,dashed, "K^{(1)}_1"]\arrow[r, "F_{0}"] \arrow[d, "K_1"] & \mathcal{O}_{\widehat{R}}(-2\widehat{S}_1-\widehat{S}_2)  \arrow[ddl,dashed, "-K^{(1)}_0"]\arrow[d, "-K_0"] \\
 \cdots \arrow[r, "F_{3}"] & \mathcal{O}_{\widehat{R}}^{m_2}(-\widehat{S}_1-\widehat{S}_2) \arrow[r, "F_{2}"] \arrow[d, "-K_3"] & \mathcal{O}_{\widehat{R}}^{m_1}(-\widehat{S}_1-\widehat{S}_2) \arrow[r, "F_{1}"] \arrow[d, "K_2"] & \mathcal{O}_{\widehat{R}}^{m_0}(-\widehat{S}_1-\widehat{S}_2) \arrow[r, "F_{0}"] \arrow[d, "-K_1"] & \mathcal{O}_{\widehat{R}}(-\widehat{S}_1-\widehat{S}_2) \arrow[d, "K_0"] \\
& \cdots \arrow[r, "F_{3}"] & \mathcal{O}_{\widehat{R}}^{m_2}(-\widehat{S}_2) \arrow[r, "F_{2}"] & \mathcal{O}_{\widehat{R}}^{m_1}(-\widehat{S}_2) \arrow[r, "F_{1}"] & \mathcal{O}_{\widehat{R}}^{m_0}(-\widehat{S}_2) \arrow[r, "F_{0}"] & \mathcal{O}_{\widehat{R}}(-\widehat{S}_2)
\end{tikzcd}
\]
and
\[
\begin{tikzcd}[column sep=small, row sep=small]
& & \vdots \arrow[d] & & \vdots \arrow[d] \\
\cdots \arrow[r, "F_{3}"] & \mathcal{O}_{\widehat{R}}^{m_2}(-2\widehat{S_1}) \arrow[r, "F_{2}"] \arrow[d, "K_3"] & \mathcal{O}_{\widehat{R}}^{m_1}(-2\widehat{S_1}) \arrow[r, "F_{1}"] \arrow[d, "-K_2"] & \mathcal{O}_{\widehat{R}}^{m_0}(-2\widehat{S_1}) \arrow[ddl,dashed, "K^{(1)}_1"]\arrow[r, "F_{0}"] \arrow[d, "K_1"] & \mathcal{O}_{\widehat{R}}(-2\widehat{S_1})  \arrow[ddl,dashed, "-K^{(1)}_0"]\arrow[d, "-K_0"] \\
& \cdots \arrow[r, "F_{3}"] & \mathcal{O}_{\widehat{R}}^{m_2}(-\widehat{S_1}) \arrow[r, "F_{2}"] \arrow[d, "-K_3"] & \mathcal{O}_{\widehat{R}}^{m_1}(-\widehat{S_1}) \arrow[r, "F_{1}"] \arrow[d, "K_2"] & \mathcal{O}_{\widehat{R}}^{m_0}(-\widehat{S_1}) \arrow[r, "F_{0}"] \arrow[d, "-K_1"] & \mathcal{O}_{\widehat{R}}(-\widehat{S_1}) \arrow[d, "K_0"] \\
& & \cdots \arrow[r, "F_{3}"] & \mathcal{O}_{\widehat{R}}^{m_2} \arrow[r, "F_{2}"] & \mathcal{O}_{\widehat{R}}^{m_1} \arrow[r, "F_{1}"] & \mathcal{O}_{\widehat{R}}^{m_0} \arrow[r, "F_{0}"] & \mathcal{O}_{\widehat{R}}
\end{tikzcd}
\]
We have the following natural morphisms:
\[
\begin{tikzcd}[row sep=small, column sep=1em]
\cdots & \mathcal{O}^{m_1}_{\widehat{R}}(-\widehat{S_1}) \arrow[rr, "F_1"] \arrow[dd,"K_1"]  & & \mathcal{O}^{m_0}_{\widehat{R}}(-\widehat{S_1}) \arrow[dd,"K_0"]\arrow[rr, "F_0"]  & & \mathcal{O}_{\widehat{R}}(-\widehat{S_1}) \arrow[dd,"K_0"]\\
\mathcal{O}^{m_0}_{\widehat{R}}(-\widehat{S_1}-\widehat{S_2}) \arrow[crossing over, "F_0"]{rr}\arrow[ur,"-L_1"] \arrow[dd,"-K_1"] & & \mathcal{O}_{\widehat{R}}(-\widehat{S_1}-\widehat{S_2}) \arrow[crossing over,"0"]{rr}\arrow[ur,"L_0"]\arrow[dashed,dl,"G_0"] \arrow[dd,"K_0"] & & 0 \arrow[ur,"0"]\arrow[dd,"0"]\arrow[dashed,dl,"0"]\\
& \mathcal{O}_{\widehat{R}}^{m_2} \arrow[rr]  & & \mathcal{O}_{\widehat{R}}^{m_1} \arrow[rr]  & & \mathcal{O}_{\widehat{R}}^{m_0} \arrow[rr, "F_0"]  & &  \mathcal{O}_{\widehat{R}}\\
\mathcal{O}^{m_1}_{\widehat{R}}(-\widehat{S_2})  \arrow[rr, "F_1"]\arrow[ur,"L_2"] & & \mathcal{O}^{m_0}_{\widehat{R}}(-\widehat{S_2})\arrow[rr, "F_0"]\arrow[ur,"-L_1"] & & \mathcal{O}_{\widehat{R}}(-\widehat{S_2}) \arrow[ur,"L_0"]\arrow[from=uu, crossing over]
\end{tikzcd}
\]
The horizontal plane in the commutative diagram is obtained by taking the standard resolution of $i_*\mathcal{O}_{\widehat{B}}$ and performing division along $w_2$. This construction is completely analogous to our previous one, with the following precise correspondence:
\[
\begin{tikzcd}[column sep=small, row sep=small]
& \cdots \arrow[r, "F_{3}"] & \mathcal{O}_{\widehat{R}}^{m_2}(-\widehat{S}_2) \arrow[r, "F_{2}"] \arrow[d, "-L_3"] & \mathcal{O}_{\widehat{R}}^{m_1}(-\widehat{S}_2) \arrow[r, "F_{1}"] \arrow[d, "L_2"] & \mathcal{O}_{\widehat{R}}^{m_0}(-\widehat{S}_2) \arrow[r, "F_{0}"] \arrow[d, "-L_1"] & \mathcal{O}_{\widehat{R}}(-\widehat{S}_2) \arrow[d, "L_0"] \\
& & \cdots \arrow[r, "F_{3}"] & \mathcal{O}_{\widehat{R}}^{m_2} \arrow[r, "F_{2}"] & \mathcal{O}_{\widehat{R}}^{m_1} \arrow[r, "F_{1}"] & \mathcal{O}_{\widehat{R}}^{m_0} \arrow[r, "F_{0}"] & \mathcal{O}_{\widehat{R}}
\end{tikzcd}
\]
it satisfies the following relations:
\[
L_i \circ F_i + F_{i+1} \circ L_{i+1} = w_2 \cdot \text{id}
\] and also induces the standard resolution of $\mathcal{O}_{\widehat{B}}$ on $\widehat{R}/(w_2)$.

\item We analyze the composition $F_1 \circ (K_1 \circ L_0 + L_1 \circ K_0)$ through the following steps:
\[
\begin{aligned}
F_1 \circ (K_1 \circ L_0 + L_1 \circ K_0) &= (F_1 \circ K_1) \circ L_0 + (F_1 \circ L_1) \circ K_0 \\
&= (w_1 - K_0 \circ F_0) \circ L_0 + (w_2 - L_0 \circ F_0) \circ K_0 \\
&= w_1 L_0 - K_0 \circ (F_0 \circ L_0) + w_2 K_0 - L_0 \circ (F_0 \circ K_0) \\
&= w_1 L_0 - w_2 K_0 + w_2 K_0 - w_1 L_0  \\
&= 0
\end{aligned}
\]
Recall we have:
  \[
  F_2 = \begin{pmatrix}
  f^{(2)}_1 & f^{(2)}_2 & \cdots & f^{(2)}_{m_2}
  \end{pmatrix},
  \]
  a standard basis for \( \mathrm{Ker}(F_1) \) in \( \widehat{R} \) under the natural module order, then applying the division theorem to each column of the \(K_1 \circ L_0 + L_1 \circ K_0\) yields:
  \[
  \sum_{l=1}^{m_2} f^{(2)}_l \cdot g^{(1)}_l + r^{(1)} = -\big(K_1 \circ L_0 + L_1 \circ K_0\big),
  \]
  with \( r^{(1)} = 0 \). Additionally, for each column, the initial terms satisfy:
  \[
  \mathrm{in}(K_1 \circ L_0 + L_1 \circ K_0) \leq \mathrm{in}(f^{(2)}_l \cdot g^{(1)}_l) \quad \text{for any } l.
  \]
  Thus, we obtain a matrix $G_0$ from \(\mathcal{O}_{\widehat{R}}(-\widehat{S_1}-\widehat{S_2})\) to \(\mathcal{O}^{m_1}_{\widehat{R}}\) such that:
  \[
  K_1 \circ L_0 + L_1 \circ K_0 + F_2\circ G_0=0.
  \]
\item Let $G_0$ be constructed as above. For each positive integer $i$, suppose by induction that $G_j$ has been constructed for all $j \leq i$ such that:
\begin{equation}
K_i \circ L_{i-1} + L_i \circ K_{i-1} + F_{i+1} \circ G_{i-1} + G_{i-2} \circ F_{i-2}=0.
\end{equation}
Then we examine the following equation:
\[
\begin{aligned}
F_{i+1} \circ \big(K_{i+1} \circ L_i + L_{i+1} \circ K_i + G_i \circ F_{i-1}\big)
&= (F_{i+1} \circ K_{i+1}) \circ L_i + (F_{i+1} \circ L_{i+1}) \circ K_i - (F_{i+1} \circ G_i) \circ F_{i-1} \\[4pt]
&= \big(w_1 - K_i \circ F_i\big) \circ L_i + \big(w_2 - L_i \circ F_i\big) \circ K_i + (F_{i+1} \circ G_i) \circ F_{i-1} \\[4pt]
&= w_1 \circ L_i - K_i \circ (F_i \circ L_i) + w_2 \circ K_i \\
&\quad - L_i \circ (F_i \circ K_i) + (F_{i+1} \circ G_i) \circ F_{i-1} \\[4pt]
&= w_1 \circ L_i - K_i \circ \big(w_1 - L_{i-1} \circ F_{i-1}\big) + w_2 \circ K_i \\
&\quad - L_i \circ \big(w_2 - K_{i-1} \circ F_{i-1}\big) + (F_{i+1} \circ G_i) \circ F_{i-1} \\[4pt]
&= w_1 \circ L_i - w_1 \circ K_i + K_i \circ L_{i-1} \circ F_{i-1} + w_2 \circ K_i \\
&\quad - w_2 \circ L_i + L_i \circ K_{i-1} \circ F_{i-1} + F_{i+1} \circ G_i \circ F_{i-1} \\[4pt]
&= \big(K_i \circ L_{i-1} + L_i \circ K_{i-1} + F_{i+1} \circ G_i\big) \circ F_{i-1}\\
&= -\big(G_{i-1} \circ F_{i-2}\big) \circ F_{i-1}=0
\end{aligned}
\]

Similarly we have:
\[
F_{i+2} = \begin{pmatrix}
f^{(i+2)}_1 & f^{(i+2)}_2 & \cdots & f^{(i+2)}_{m_{i+2}}
\end{pmatrix},
\]
a standard basis for \( \mathrm{Ker}(F_{i+1}) \) on \( \widehat{R} \) under the natural module order, then applying the division theorem to each column of the
\[
K_{i+1} \circ L_i + L_{i+1} \circ K_i + G_i \circ F_{i-1}
\]
yields a matrix \( G_{i+1} \colon \mathcal{O}^{m_i}_{\widehat{R}}(-\widehat{S_1} - \widehat{S_2}) \to \mathcal{O}^{m_{i+1}}_{\widehat{R}} \) such that:
\[
K_{i+1} \circ L_i + L_{i+1} \circ K_i + F_{i+2} \circ G_{i+1} + G_i \circ F_{i-1} = 0.
\]

In particular, this precisely serves as the required initial condition for the next step, hence for every integer \( k \) in the range, we can construct \( G_k \) via division as above.

\item
For notational convenience, we introduce the following new notation:
\[
K_i^{(a+1,\, 0,\,-a)} := K_i^{(a)}, \quad
K_i^{(0,\, 1,\, 0)} := L_i, \quad
K_i^{(1,\, 1,\,-1)} := G_i.
\]
We regard the superscripts of \( K \) as vectors, subject to the basic operations of a group. Then for any positive integer \( a \), we can construct matrices \( K_i^{(a,1,-a)} \) satisfying the following recursive compatibility condition:
\[
\sum_{\substack{v + w = (a,1,-a+1)}}
K_{\alpha(w,i)}^v \circ K_i^w
+ K_{i-1}^{(a,1,-a)} \circ F_{i-1}
+ F_{\alpha((a,1,-a),i)-1} \circ K_i^{(a,1,-a)}
= 0,
\]
where the sum ranges over all decompositions \( v + w = (a,1,-a+1) \) of multi-indices such that \( v, w < (a,1,-a+1) \) under the natural vector order that we already constructed, and \( \alpha(v, i) \) denotes the index shift function determined by the image of \( K^v_i \). We assume that for all multi-indices \( v \leq (a,1,-a) \) and all integers \( i' \),
as well as for \( v = (a,1,-a) \) and all \( i' < i \),
the matrices \( K_{i'}^v \) have already been constructed and satisfy the above recursion condition. Then, we claim that \( K_i^{(a,1,-a)} \) can also be constructed:
\[
\left(F_{\alpha((a,1,-a), i) - 2} \circ K_{i-1}^{(a,1,-a)}\right) \circ F_{i-1}
= -\left(\sum_{\substack{v + w = (a,1,-a+1)}}
K_{\alpha(w, i-1)}^v \circ K_{i-1}^w \right)\circ F_{i-1}
\]
\[
F_{\alpha((a-1,1,-a), i)-2} \circ
\left( \sum_{\substack{v + w = (a,1,-a+1)}}
K_{\alpha(w,i)}^v \circ K_i^w \right)
=
\sum_{\substack{v + w = (a,1,-a+1) }}
\left( F_{\alpha((a-1,1,-a), i)-2} \circ K_{\alpha(w,i)}^v \right) \circ K_i^w
\]
noticing \( v \leq (a,1,-a+1)  \), the following holds by assumption for any integer $i$
\[
\sum_{\substack{v' + w' = v+(0,0,1)}}
K_{\alpha(w',i)}^{v'} \circ K_i^{w'}
+ K_{i-1}^{v} \circ F_{i-1}
+ F_{\alpha(v,i)-1} \circ K_i^{v}
= 0,
\]
so we have
\begin{align*}
&\sum_{\substack{v + w = (a,1,-a+1) }}
\left( F_{\alpha((a-1,1,-a), i)-2} \circ K_{\alpha(w,i)}^v \right) \circ K_i^w \\
&= \sum_{\substack{v + w = (a,1,-a+1) }}
\left(
- \sum_{\substack{v' + w' = v +(0,0,1)}} K_{\alpha(w',\alpha(w,i))}^{v'} \circ K_{\alpha(w,i)}^{w'}
- K_{\alpha(w,i)-1}^v \circ F_{\alpha(w,i)-1}
\right) \circ K_i^w \\
&= -\sum_{\substack{v + w = (a,1,-a+1)  \\ v' + w' = v +(0,0,1)}}
K_{\alpha(w',\alpha(w,i))}^{v'} \circ K_{\alpha(w,i)}^{w'} \circ K_i^w \\
&\quad -\sum_{\substack{v + w = (a,1,-a+1) }}
K_{\alpha(w,i)-1}^v \circ \left(F_{\alpha(w,i)-1} \circ K_i^w\right) \\
&= -\sum_{\substack{v + w = (a,1,-a+1)  \\ v' + w' = v +(0,0,1)}}
K_{\alpha(w',\alpha(w,i))}^{v'} \circ K_{\alpha(w,i)}^{w'} \circ K_i^w \\
&\quad -\sum_{\substack{v + w = (a,1,-a+1) }}
K_{\alpha(w,i)-1}^v \circ \left(
- \sum_{\substack{v' + w' = w}}
K_{\alpha(w',i)}^{v'} \circ K_i^{w'}
- K_{i-1}^w \circ F_{i-1}
\right) \\
&= -\sum_{\substack{v + w = (a,1,-a+1)  \\ v' + w' = v +(0,0,1)}}
K_{\alpha(w',\alpha(w,i))}^{v'} \circ K_{\alpha(w,i)}^{w'} \circ K_i^w \\
&\quad + \sum_{\substack{v + w = (a,1,-a+1)  \\ v' + w' = w+(0,0,1)}}
K_{\alpha(w,i)-1}^v \circ K_{\alpha(w',i)}^{v'} \circ K_i^{w'} \\
&\quad + \sum_{\substack{v + w = (a,1,-a+1) }}
K_{\alpha(w,i)-1}^v \circ K_{i-1}^w \circ F_{i-1}=\sum_{\substack{v + w = (a,1,-a+1) }}
K_{\alpha(w,i)-1}^v \circ K_{i-1}^w \circ F_{i-1}
\end{align*}
Therefore,
\[
F_{\alpha((a,1,-a),i)-2} \left(
\sum_{\substack{v + w = (a,1,-a+1)}}
K_{\alpha(w,i)}^v \circ K_i^w
+ K_{i-1}^{(a,1,-a)} \circ F_{i-1}
\right)
= 0.
\]
By the division theorem with respect to the standard basis, we obtain a matrix
\( K_i^{(a,1,-a)} \)
such that
\[
\sum_{\substack{v + w = (a,1,-a+1) }}
K_{\alpha(w,i)}^v \circ K_i^w
+ K_{i-1}^{(a,1,-a)} \circ F_{i-1}
+ F_{\alpha((a,1,-a),i)-1} \circ K_i^{(a,1,-a)}
= 0.
\]
\item
From the commutative diagram, it is clear that for any non-negative integer \( a \) and any integer \( i \),
only finitely many of the recursively constructed matrices \( K_i^{(a,1,-a)} \) are nonzero. Moreover, together they define a well-defined non-trivial element  \(s({\mathcal{O}_{\widehat{B}}})\) in the morphism
\[
\mathrm{Hom}\left(
j_2^{*}i_{1*}\mathcal{O}_{\widehat{B}}(-\widehat{S}_2)[1],\
j_2^{*}i_{1*}\mathcal{O}_{\widehat{B}}
\right).
\]
for any element \( t \) of \( {\mathcal{O}_{\widehat{B}}} \), it induces a map:
\[
t: j_2^{*}i_{1*}\mathcal{O}_{\widehat{B}} \longrightarrow j_2^{*}i_{1*}\mathcal{O}_{\widehat{B}}.
\]
We consider the following composition:
\[
\begin{tikzcd}
j_2^{*}i_{1*}\mathcal{O}_{\widehat{B}}(-\widehat{S}_2)[1] \arrow[r,"s({\mathcal{O}_{\widehat{B}}})"] & j_2^{*}i_{1*}\mathcal{O}_{\widehat{B}} \arrow[r, "t"] & j_2^{*}i_{1*}\mathcal{O}_{\widehat{B}}.
\end{tikzcd}
\]
based on cohomology calculations, we obtain all the mappings.
\item Next, we consider the \(\mathcal{O}_{\widehat{B}}\) family:
\[
\mathrm{Hom}\left(
j_2^{*}i_{1*}\mathcal{O}_{\widehat{B}}(-2\widehat{S}_2)[2],\
\mathcal{O}_{\widehat{B}_2}
\right).
\]
Completely consistent for any non-negative integer \( a \), we can construct matrices \( K_i^{(a,2,-a-1)} \)  satisfying the following recursive compatibility condition:
\[
\sum_{\substack{v + w = (a,2,-a)}}
K_{\alpha(w,i)}^v \circ K_i^w
+ K_{i-1}^{(a,2,-a-1)} \circ F_{i-1}
+ F_{\alpha((a,2,-a-1),i)-1} \circ K_i^{(a,2,-a-1)}
= 0,
\]
purely via division. In particular, for such a recursion relation, we have
\[
K_i^{(0,2,-1)} = L_i^{(1)}
\]
since \(a=0\) and \((0,2,0)=(0,1,0)+(0,1,0)\).
\item By continuing the induction process for any \( b \geq 2 \), we obtain the \( \mathcal{O}_{\widehat{B}} \) family:
\[
\mathrm{Hom}\left(
j_2^{*}i_{1*}\mathcal{O}_{\widehat{B}}(-b\widehat{S}_2)[b],\
\mathcal{O}_{\widehat{B}_b}
\right).
\]
by considering, for any non-negative integer \( a \), the matrices \( K_i^{(a,b,-a-b+1)} \) satisfying the following recursive condition:
\[
\sum_{\substack{v + w =(a,b,-a-b+2)}}
K_{\alpha(w,i)}^v \circ K_i^w
+ K_{i-1}^{(a,b,-a-b+1)} \circ F_{i-1}
+ F_{\alpha((a,b,-a-b+1),i)-1} \circ K_i^{(a,b,-a-b+1)}
= 0,
\]
purely via division. In particular, for such a recursion relation, we have
\[
K_i^{(0,b,-b+1)} = L_i^{(b)}.
\]
Through this relation, it is not difficult to see that \( K_i^{(a,b,-a-b+1)} \), for any \( a \), gives rise to a nontrivial morphism as described above, via the commutative diagram. Moreover, by scaling and extension, we obtain all such morphisms.

\item We consider the following higher (cubic) block matrices:
\[
A := \Bigg\langle
\begin{pmatrix}
K_0^{(0,0)}:=0 & K_0^{(0,1,0)} & K_0^{(0,2,-1)} &\\
K_0^{(1,0,0)} & K_0^{(1,1,-1)} & K_0^{(1,2,-2)} & \cdots \\
K_0^{(2,0,-1)} & K_0^{(2,1,-2)} & K_0^{(2,2,-3)} & \\
 & \vdots &  & \ddots
\end{pmatrix}
\ \Bigg|\
\begin{pmatrix}
0 &  &  &  &  \\
& K_2^{(0,0)}:=F_1 & K^{(0,1,0)}_2 & K^{(0,2,-1)}_2 &  \\
& K^{(1,0,0)}_2 & K^{(1,1,-1)}_2  & K^{(1,2,-2)}_2 & \cdots \\
& K^{(2,0,-1)}_2 & K^{(2,1,-2)}_2 & K^{(2,2,-3)}_2 &  \\
&  & \vdots &  & \ddots
\end{pmatrix}
\ \Bigg|\
\cdots
\Bigg\rangle.
\]
\[
B := \Bigg\langle
\begin{pmatrix}
K_1^{(0,0)}:=F_0 & K_1^{(0,1,0)} & K_1^{(0,2,-1)} &\\
K_1^{(1,0,0)} & K_1^{(1,1,-1)} & K_1^{(1,2,-2)} & \cdots \\
K_1^{(2,0,-1)} & K_1^{(2,1,-2)} & K_1^{(2,2,-3)} & \\
 & \vdots &  & \ddots
\end{pmatrix}
\ \Bigg|\
\begin{pmatrix}
0 &  &  &  &  \\
& K_3^{(0,0)}:=F_2 & K^{(0,1,0)}_3 & K^{(0,2,-1)}_3 &  \\
& K^{(1,0,0)}_3 & K^{(1,1,-1)}_3  & K^{(1,2,-2)}_3 & \cdots \\
& K^{(2,0,-1)}_3 & K^{(2,1,-2)}_3 & K^{(2,2,-3)}_3 &  \\
&  & \vdots &  & \ddots
\end{pmatrix}
\ \Bigg|\
\cdots
\Bigg\rangle.
\]

where \( A \) (or \( B \)) is a finite cubic matrix, that is, a \( k' \)-family of matrices, denoted by
\[
A = \Braket{ A_1 \mid A_2 \mid \cdots \mid A_{k'} }.
\]
$A_i$ represents the $i$-th layer of our cubic matrix and in general we assume that \( 2k' \leq k+1 \). This cubic matrix with its sub blocks is capable of reproducing all of the above morphisms. Perhaps we could define a multiplication for higher matrices so that the pair \( (A, B) \) resembles a matrix factorization more closely. However, we only consider their projections along the corresponding directions:
\[
A_K:= \begin{pmatrix}
K_0 & F_1 & & & \\
K_0^{(1)} & K_2 & F_3 & & \\
K_0^{(2)} & K_2^{(1)} & K_4 & \ddots & \\
\vdots & \vdots & \ddots & \ddots & F_{2k'-1} \\
K_0^{(k')} & K_2^{(k'-1)} & \cdots & K_{2k'-2}^{(1)} & K_{2k'}
\end{pmatrix}, \quad
B_K: =\begin{pmatrix}
F_0 &        &        &        &        \\
K_1 & F_2   &        &        &        \\
K_1^{(1)} & K_3 & F_4    &        &        \\
\vdots & \vdots & \ddots & \ddots &        \\
K_1^{(k'-1)} & K_3^{(k'-2)} & \cdots & K_{2k'-1} & F_{2k'}
\end{pmatrix}
\]
correspond exactly to the matrix factorization of \( w_1 \).  In particular, if \( w_1 \) and \( w_2 \) additionally form a standard basis or a measure basis, then by considering weighted projections, we obtain matrix factorizations for arbitrary element in $\widehat{J}$.
\item From the construction of the excess distinguished triangle, we know that after applying the iterated totalization of the cubic complex along the morphism we constructed, we obtain the desired projective resolution of \( i_{2*}\mathcal{O}_B \) on \(\widehat{S}\). However, we can also verify directly that the above construction yields an exact sequence via the division relations, i.e. for the first several terms it involves the following complex:
    \[
\mathcal{O}^{m_1}_{\widehat{R}}(-\widehat{S_1})\oplus\mathcal{O}^{m_1}_{\widehat{R}}(-\widehat{S_2})\oplus\mathcal{O}^{m_2}_{\widehat{R}}
\xrightarrow{
\mathbf{F}_{2} :=
\left(\begin{smallmatrix}
F_{0} & 0 & 0 \\
 0& F_{0} &  0\\
L_{1} & K_{1} & F_{2}
\end{smallmatrix}\right)
}
\mathcal{O}^{m_0}_{\widehat{R}}(-\widehat{S_1})\oplus\mathcal{O}^{m_0}_{\widehat{R}}(-\widehat{S_2})\oplus\mathcal{O}^{m_1}_{\widehat{R}}
\xrightarrow{\mathbf{F}_{1} :=
\left(\begin{smallmatrix}
K_{0}&L_{0}& F_{1}\\
\end{smallmatrix}\right)}\mathcal{O}^{m_0}_{\widehat{R}}
\xrightarrow{F_0}\mathcal{O}_{\widehat{R}}
\]
If we have \( t = (x, y, z) \) such that \( \mathbf{F}_1(t) = 0 \) on \( \widehat{S} \), then there exist \( p, q \in \widehat{R} \) such that
\[
K_0 \circ x + L_0 \circ y + F_1 \circ z = p w_1 + q w_2 \quad \text{in } \widehat{R}.
\]
Applying \( F_0 \) to both sides, we obtain
\[
w_1 \circ (x - F_0 \circ p) + w_2 \circ (y - F_0 \circ q) = 0.
\]
Since \( \langle w_1, w_2 \rangle \) is a regular sequence, it follows that there exists \( t \in \widehat{R} \) such that
\[
x = F_0 \circ p + w_2 \circ t, \quad y = F_0 \circ q - w_1 \circ t.
\]
Substituting back into the original expression, we have
\[
K_0 \circ (F_0 \circ p + w_2 \circ t) + L_0 \circ (F_0 \circ q - w_1 \circ t) + F_1 \circ z = p w_1 + q w_2.
\]
Therefore,
\[
K_0 \circ w_2 \circ t - L_0 \circ w_1 \circ t + F_1 \circ z = F_1 \circ K_1 \circ p + F_1 \circ K_1 \circ q.
\]
Note that
\[
K_0 \circ w_2 \circ t - L_0 \circ w_1 \circ t = (L_0 \circ F_0 + F_1 \circ L_1) \circ K_0 \circ t - L_0 \circ w_1 \circ t = F_1 \circ L_1 \circ K_0 \circ t,
\]
so we have
\[
F_1 \circ z = F_1 \circ (K_1 \circ p + K_1 \circ q - L_1 \circ K_0 \circ t),
\]
which implies there exists \( r \in \widehat{R} \) such that
\[
z = K_1 \circ p + K_1 \circ q - L_1 \circ K_0 \circ t + F_2 \circ r.
\]
Hence,
\[
x = F_0 \circ (p + L_0 \circ t), \quad y = F_0 \circ (q - K_0 \circ t), \quad z = K_1 \circ (p + q) - L_1 \circ K_0 \circ t + F_2 \circ r.
\]
This confirms that the sequence is exact at degree $-1$. Similar arguments apply to the other terms, and note that the condition that \( w_1, w_2 \) form a regular sequence is essential. Most importantly, we want to mention is that the division construction and the cohomology computations are independent processes. That is to say, we can prove that the cubic factorization obtained via the division construction, once totalized modulo \( \widehat{J} \), yields an exact sequence by methods of induction on triangular decomposition. We denote this exact sequence by
\[
\begin{tikzcd}[column sep=normal]
\cdots \arrow[r] &
\mathbf{V}_{k+2} \arrow[r, "\mathbf{F}_{k+2}"] &
\mathbf{V}_{k+1} \arrow[r, "\mathbf{F}_{k+1}"] &
\mathbf{V}_k \arrow[r, "\mathbf{F}_{k}"] &
\cdots \arrow[r, "\mathbf{F}_1"] &
\mathbf{V}_0 \arrow[r, "\mathbf{F}_0"] &
\widehat{R} \arrow[r, "\mathbf{F}_{-1}"] &
\mathcal{O}_{\widehat{B}} \arrow[r] &
0
\end{tikzcd}
\]

\end{enumerate}

On the other hand, if we assume that $w_1, w_2$ additionally form a measure basis, i.e., we impose the extra assumption that $L(J) = \langle L(w_1), L(w_2) \rangle$ with $L(w_1), L(w_2)$ forming a regular sequence, then by repeating the above procedure for $L(B)$, $L(J)$, $L(w_1)$, and $L(w_2)$, we obtain a similar projective resolution of $\mathcal{O}_{L(\widehat{B})}$ in $\widehat{R}/L(\widehat{J})$. This resolution is induced by the following building block:
\[
\begin{tikzcd}[row sep=small, column sep=1em]
\cdots & \mathcal{O}^{m_1}_{\widehat{R}}(-*) \arrow[rr, "L(F_1)"] \arrow[dd,"L(K_1)"]  & & \mathcal{O}^{m_0}_{\widehat{R}}(-*) \arrow[dd,"L(K_0)"]\arrow[rr, "L(F_0)"]  & & \mathcal{O}_{\widehat{R}}(-*) \arrow[dd,"L(K_0)"]\\
\mathcal{O}^{m_0}_{\widehat{R}}(-2*) \arrow[crossing over, "L(F_0)"]{rr}\arrow[ur,"-L(L_1)"] \arrow[dd,"-L(K_1)"] & & \mathcal{O}_{\widehat{R}}(-2*) \arrow[crossing over,"0"]{rr}\arrow[ur,"L(L_0)"]\arrow[dashed,dl,"L(G_1)"] \arrow[dd,"L(K_0)"] & & 0 \arrow[ur,"0"]\arrow[dd,"0"]\arrow[dashed,dl,"-L(G_0)"]\\
& \mathcal{O}_{\widehat{R}}^{m_2} \arrow[rr]  & & \mathcal{O}_{\widehat{R}}^{m_1} \arrow[rr]  & & \mathcal{O}_{\widehat{R}}^{m_0} \arrow[rr, "L(F_0)"]  & &  \mathcal{O}_{\widehat{R}}\\
\mathcal{O}^{m_1}_{\widehat{R}}(-*)  \arrow[rr, "L(F_1)"]\arrow[ur,"L(L_2)"] & & \mathcal{O}^{m_0}_{\widehat{R}}(-*)\arrow[rr, "L(F_0)"]\arrow[ur,"-L(L_1)"] & & \mathcal{O}_{\widehat{R}}(-*) \arrow[ur,"L(L_0)"]\arrow[from=uu, crossing over]
\end{tikzcd}
\]
where $(*)$ denotes twist by the corresponding leading divisors. We have leading cubic matrices:

\[
L(A) := \Bigg\langle
\begin{pmatrix}
L(K_0^{(0,0)} := 0) & L(K_0^{(0,1,0)}) & L(K_0^{(0,2,-1)}) &\\
L(K_0^{(1,0,0)}) & L(K_0^{(1,1,-1)}) & L(K_0^{(1,2,-2)}) & \cdots \\
L(K_0^{(2,0,-1)}) & L(K_0^{(2,1,-2)}) & L(K_0^{(2,2,-3)}) & \\
 & \vdots &  & \ddots
\end{pmatrix}
\ \Bigg|\
\cdots
\Bigg\rangle.
\]

\[
L(B) := \Bigg\langle
\begin{pmatrix}
L(K_1^{(0,0)} := F_0) & L(K_1^{(0,1,0)}) & L(K_1^{(0,2,-1)}) &\\
L(K_1^{(1,0,0)}) & L(K_1^{(1,1,-1)}) & L(K_1^{(1,2,-2)}) & \cdots \\
L(K_1^{(2,0,-1)}) & L(K_1^{(2,1,-2)}) & L(K_1^{(2,2,-3)}) & \\
 & \vdots &  & \ddots
\end{pmatrix}
\ \Bigg|\
\cdots
\Bigg\rangle.
\]
Moreover, according to the construction via division, the matrix \(L(K_i^{(a,b,-a-b+1)})\) appearing above is precisely the leading term matrix of \(K_i^{(a,b,-a-b+1)}\), taken with respect to the natural module order on domain and codomain. Concretely, each entry satisfies:
\[
\{L(K_i^{(a,b,-a-b+1)})\}_{m,n} = \{K_i^{(a,b,-a-b+1)}\}_{m,n} \mod \widehat{P}^{\,ad_1 + bd_2 + a^{(i-1)}_m - a^{(i+2a+2b-2)}_n + 1},
\]
where \(d_i = \operatorname{ord}_{P}(w_i)\) for \(i = 1, 2\).\\

Using the exact diagram above and following the same degree cancellation method as before, we can prove that for any integer $r$, the module $\widehat{P}^r \cdot \mathcal{O}_{\widehat{B}}$ admits a standard resolution over $\widehat{S}$. This resolution is induced by the following building blocks:

\[
\begin{tikzcd}[row sep=small, column sep=0.6em]
\cdots & \bigoplus\limits^{m_1}_{j=1} \widehat{P}^{r - a^{(1)}_j - d_1} \arrow[rr, "F_1"] \arrow[dd,"K_1"]  & & \bigoplus\limits^{m_0}_{i=1} \widehat{P}^{r - a^{(0)}_i - d_1} \arrow[dd,"K_0"]\arrow[rr, "F_0"]  & & \widehat{P}^{r - d_1} \arrow[dd,"K_0"]\\
\bigoplus\limits^{m_0}_{i=1} \widehat{P}^{r - a^{(0)}_i - d_1 - d_2} \arrow[crossing over, "F_0"]{rr}\arrow[ur,"-L_1"] \arrow[dd,"-K_1"] & & \widehat{P}^{r - d_1 - d_2} \arrow[crossing over,"0"]{rr}\arrow[ur,"L_0"]\arrow[dashed,dl,"G_0"] \arrow[dd,"K_0"] & & 0 \arrow[ur,"0"]\arrow[dd,"0"]\arrow[dashed,dl,"0"]\\
& \bigoplus\limits^{m_2}_{l=1} \widehat{P}^{r - a^{(2)}_l} \arrow[rr]  & & \bigoplus\limits^{m_1}_{j=1} \widehat{P}^{r - a^{(1)}_j} \arrow[rr]  & & \bigoplus\limits^{m_0}_{i=1} \widehat{P}^{r - a^{(0)}_i} \arrow[rr, "F_0"]  & &  \widehat{P}^r\\
\bigoplus\limits^{m_1}_{j=1} \widehat{P}^{r - a^{(1)}_j - d_2}  \arrow[rr, "F_1"]\arrow[ur,"L_2"] & & \bigoplus\limits^{m_0}_{i=1} \widehat{P}^{r - a^{(0)}_i - d_2} \arrow[rr, "F_0"]\arrow[ur,"-L_1"] & & \widehat{P}^{r - d_2} \arrow[ur,"L_0"]\arrow[from=uu, crossing over]
\end{tikzcd}
\]
We have a bounded above acyclic complex on $\widehat{S}$:
\[
\begin{tikzcd}[column sep=small]
 \cdots & \bigoplus^{m_{2}}_{l=1} \widehat{P}^{r-\mathrm{a}^{(2)}_l}\oplus\cdots \arrow[r, "\mathbf{F}_{2}"] & \bigoplus^{m_{1}}_{l=1} \widehat{P}^{r-\mathrm{a}^{(1)}_l} \oplus \widehat{P}^{r-d_1}(-\widehat{S}_1) \oplus \widehat{P}^{r-d_2}(-\widehat{S}_2)\arrow[r, "\mathbf{F}_{1}"] & \bigoplus^{m_{0}}_{l=1} \widehat{P}^{r-\mathrm{a}^{(0)}_l} \arrow[r, "\mathbf{F}_{0}"]&\widehat{P}^r
\end{tikzcd}
\]
and denote this complex by $\mathbf{U}^L(\widehat{P}^r \otimes \mathcal{O}_{\widehat{B}})$. On the other hand, analogous to the case of hypersurfaces, this resolution can be extended to $\mathbb{Z}$-graded. Considering the right resolution of $\widehat{P}^r \cdot \mathcal{O}_{\widehat{B}}$, we employ the excess distinguished triangle method in the opposite direction to obtain the forward extension through the following building blocks:
\[
\begin{tikzcd}[row sep=small, column sep=1em]
\cdots & \mathcal{O}^{m_1}_{\widehat{R}}(\widehat{S_2}) \arrow[rr, "F_1"] \arrow[dd,"K_1"]  & & \mathcal{O}^{m_0}_{\widehat{R}}(\widehat{S_2}) \arrow[dd,"K_0"]\arrow[rr, "F_0"]  & & \mathcal{O}_{\widehat{R}}(\widehat{S_2}) \arrow[dd,"K_0"]\\
\mathcal{O}^{m_0}_{\widehat{R}} \arrow[crossing over, "F_0"]{rr}\arrow[ur,"-L_1"] \arrow[dd,"-K_1"] & & \mathcal{O}_{\widehat{R}} \arrow[crossing over,"0"]{rr}\arrow[ur,"L_0"]\arrow[dashed,dl,"G_0"] \arrow[dd,"K_0"] & & 0 \arrow[ur,"0"]\arrow[dd,"0"]\arrow[dashed,dl,"0"]\\
& \mathcal{O}_{\widehat{R}}^{m_2}(\widehat{S_1}+\widehat{S_2}) \arrow[rr]  & & \mathcal{O}_{\widehat{R}}^{m_1}(\widehat{S_1}+\widehat{S_2}) \arrow[rr]  & & \mathcal{O}_{\widehat{R}}^{m_0} (\widehat{S_1}+\widehat{S_2})\arrow[rr, "F_0"]  & &  \mathcal{O}_{\widehat{R}}(\widehat{S_1}+\widehat{S_2})\\
\mathcal{O}^{m_1}_{\widehat{R}}(\widehat{S_1})  \arrow[rr, "F_1"]\arrow[ur,"L_2"] & & \mathcal{O}^{m_0}_{\widehat{R}}(\widehat{S_1})\arrow[rr, "F_0"]\arrow[ur,"-L_1"] & & \mathcal{O}_{\widehat{R}}(\widehat{S_1}) \arrow[ur,"L_0"]\arrow[from=uu, crossing over]
\end{tikzcd}
\]
Then we have a bounded below acyclic complex on $\widehat{S}$ quasi-isomorphic to $\widehat{P}^r \cdot \mathcal{O}_{\widehat{B}}$:
\[
\begin{tikzcd}[column sep=small]
 \bigoplus^{m_{k}}_{l=1} \widehat{P}^{r-\mathrm{a}^{(k)}_l} \arrow[r, "\mathbf{F}'_{-k}"] & \bigoplus^{m_{k-1}}_{l=1} \widehat{P}^{r-\mathrm{a}^{(k-1)}_l} \arrow[r, "\mathbf{F}'_{-k+1}"] & \bigoplus^{m_{k}}_{l=1} \widehat{P}^{r-\mathrm{a}^{(k)}_l}(\widehat{S}_1)\oplus \bigoplus^{m_{k}}_{l=1} \widehat{P}^{r-\mathrm{a}^{(k)}_l}(\widehat{S}_2)\arrow[r, "\mathbf{F}'_{-k+2}"] & \cdots
\end{tikzcd}
\]
By gluing these two $\widehat{P}^r \cdot \mathcal{O}_{\widehat{B}}$ (standard) resolutions into an unbounded exact (standard) sequence on $\widehat{S}$, as a generalization of hypersurface matrix factorizations, we obtain the following representation:
\[
\begin{tikzcd}[column sep=normal]
\cdots \arrow[r] &
\mathbf{U}^r_{2} \arrow[r, "\mathbf{\bar{F}}_{2}"] &
\mathbf{U}^r_{1} \arrow[r, "\mathbf{\bar{F}}_{1}"] &
\mathbf{U}^r_{0} \arrow[r, "\mathbf{\bar{F}}_0"] &
\mathbf{U}^r_{-1} \arrow[r, "\mathbf{\bar{F}}_{-1}"] &
\mathbf{U}^r_{-2} \arrow[r, "\mathbf{\bar{F}}_{-2}"] &
\cdots
\end{tikzcd}
\]
where each $\mathbf{U}^r_i$ can be expressed as a direct sum of $\widehat{P}^{r+a}$ for various $a$, we typically denote this complex by $\mathbf{U}(\widehat{P}^r \cdot\mathcal{O}_{\widehat{B}})$.\\

For any Cohen-Macaulay module \( M \) on \( \widehat{S} \), and any Bourbaki exact sequence associated with \( M \):
\[
0 \longrightarrow L \longrightarrow M \longrightarrow I(M) \longrightarrow 0,
\]
we examine the corresponding Bourbaki cycle \( \widehat{B}(M) \) and consider the standard resolution $\mathbf{U}(\widehat{P}^r\,\mathcal{O}_{\widehat{B}(M)})$.
It is straightforward to see that, we can also extend the above exact sequences to obtain the standard version of the Bourbaki exact sequences for any integer \(r\):
\[
0 \longrightarrow \mathbf{L}^r \longrightarrow \mathbf{M}^r \longrightarrow \widehat{P}^r \longrightarrow \widehat{P}^r \mathcal{O}_{\widehat{B}(M)} \longrightarrow 0,
\]
\( \mathbf{M}^r \) has the following \textbf{standard resolution}:
\[
\begin{tikzcd}[column sep=normal]
\cdots \arrow[r] &
\mathbf{U}^r_4 \arrow[r, "\mathbf{\bar{F}}_4"] &
\mathbf{U}^r_3 \arrow[r, "\mathbf{\bar{F}}_3"] &
\mathbf{U}^r_2 \arrow[r, "\mathbf{\bar{F}}_2"] &
\mathbf{U}^r_1 \arrow[r, "\mathbf{\bar{F}}_1"] &
\mathbf{U}^r_0 \arrow[r, "\mathbf{\bar{F}}_0"] &
\mathbf{M}^r \arrow[r] &
0
\end{tikzcd}
\]
and \( \mathbf{L}^r \) also admits a standard resolution of the form:
\[
0 \longrightarrow \mathbf{W}^r_{k} \longrightarrow \cdots \longrightarrow \mathbf{W}^r_{1} \longrightarrow \mathbf{L}^r \to 0
\]
where the module \( \mathbf{W}_i^r \) as the quotient:
\[
\mathbf{W}_i^r := \mathbf{U}_i(\widehat{P}^r \cdot\mathcal{O}_{\widehat{B}}) / \mathbf{U}^{L}_i(\widehat{P}^r \cdot\mathcal{O}_{\widehat{B}}),
\]
for any integer \( i \), and the quotient is taken via the natural embedding of the submodules.
\begin{lemma}
According to the lemma, we have $\mathbf{M}^r$ stably equivalent to $M$ when $r \ll 0$. In particular, if $M$ has no free direct summand and we consider the reduced $\mathbf{M}^r_{\mathrm{red}}$, then $\mathbf{M}^r_{\mathrm{red}} = M$.
\end{lemma}

It is not difficult to see that the above discussion can be naturally extended to higher complete intersections through appropriate combinations. To avoid presenting a technically rigorous but ultimately vacuous structure, we state the following result without providing a complete proof:

\begin{claim}
Let $R$ be a regular local ring with maximal ideal $\mathfrak{m}$ and residue field \(\mathrm{k}\), $J$ an ideal of $R$ that admits a measure basis which is also a regular sequence. Let $P$ be a prime ideal containing $J$ such that the quotient ring $R/P$ is regular, and let $B$ be any cycle defined by an arbitrary ideal of $R$. Then, in the analytic sense, we obtain from the structure sheaf of $B$ a standard resolution in the quotient ring $R/J$.
\end{claim}
Analogously,
\begin{claim}
Under the same conditions as above, for any Cohen-Macaulay module $M$ on $S$, there exists a standard resolution of $M$ in the analytic sense.
\end{claim}

\quad\\
\author{HAO XINGBANG}

\newcommand{\Addresses}{{
  \bigskip
  \footnotesize

  HAO XINGBANG, \textsc{Waseda University, Ookubo, Shinjuku-ku, Tokyo, 169-8555, Japan.}\par\nopagebreak
  \textit{E-mail address}: \texttt{hao@fuji.waseda.jp}
}}
\Addresses

\end{document}